\documentclass[10pt]{article}

\usepackage{a4wide}
\usepackage{amssymb}
\usepackage{amsfonts}
\usepackage{amsmath}
\input xy
\xyoption{arrow} \xyoption{matrix}

\date{}

\newtheorem{proposition}{Proposition}[section]
\newtheorem{theorem}[proposition]{Theorem}
\newtheorem{lemma}[proposition]{Lemma}

\newtheorem{corollary}[proposition]{Corollary}

\def\GK{{\rm  GK}\,}

\def\Hom{{\rm Hom}}
\def\der{\partial }

\def\nFM0{{\nu }_{F,M_0}}
\def\nFN0{{\nu }_{F,N_0}}
\def\nGN0{{\nu }_{G,N_0}}

\def\N0{ {\bf N}_0 }

\def\t{\otimes}
\def\g{\gamma}
\def\v{\varphi}
\def\ra{\rightarrow}

\def\lra{\leftrightarrow}
\def\Xpm{X^{\pm }}

\def\s{\sigma}
\def\Z{\mathbb{Z}}

\def\l1{{\lambda}_1}

\def\a{\alpha}
\def\a0{ {\alpha }_0}
\def\a1{ {\alpha }_1}

\def\l{\lambda}
\def\o{\omega}

\def\nFGM0{{\nu }_{F,G,M_0}}

\def\nFN0{{\nu}_{F,N_0}}

\def\sm{{\sigma}^m}

\def\sm1{{\sigma}^{-1}}

\def\smtp1{{\sigma}^{-t+1}}

\def\o{\omega }
\def\S1{S^{-1}}

\def\Xpm1{X^{\pm 1}_1}

\def\sPM1{{\sigma }^{\pm 1}}
\def\sMP1{{\sigma }^{\mp 1 }}

\def\d{\delta}

\def\di{{\rm d.ind}}

\def\L{\Lambda}

\def\G{\Gamma}

\def\CA{{\cal A}}

\def\CD{{\cal D}}

\def\Ytm1{Y^{t-1}}
\def\Yim1{Y^{i-1}}

\def\CL{{\cal L}}

\def\CG{{\cal G}}
\def\CH{{\cal H}}

\def\Aut{{\rm Aut}}

\def\ad{{\rm ad }}
\def\dim{{\rm dim }}

\def\ker{ {\rm ker } }

\def\Ev{ {\rm Ev} }

\def\SL2Z{ {\rm SL}_2({\bf Z}) }

\def\th{ \theta }

\def\CL{{\cal L}}

\def\Gp1{ G^{1 , 1 } }
\def\P11{ P^{-1 , 1 } }
\def\Pp1{ P^{1 , 1 } }

\def\Supp{{\rm Supp}}

\def\th{\theta}

\def\nCLsr{{}^\nu\kern-2pt {\cal L}^{\sigma , \rho  }}
\def\nP{{}^\nu \kern-2pt P}
\def\nL{{}^\nu\kern-2pt L}
\def\nLL{{}^\nu\kern-2pt \Lambda}
\def\nPsr{{}^\nu\kern-2pt P^{\sigma , \rho  }}
\def\nLsr{{}^\nu\kern-2pt L^{\sigma , \rho  }}
\def\nuCL{{}^\nu\kern-2pt  {\cal L}}
\def\nCLsr{{}^\nu\kern-2pt {\cal L}^{\sigma , \rho  }}
\def\nCL1m{{}^\nu\kern-2pt {\cal L}^{-1 , 1  }}
\def\x1nu{x^\frac{1}{\nu}}
\def\xm1nu{x^{-\frac{1}{\nu}}}

\def\ob{\overline{b}}

\def\ra{\rightarrow }

\def\CI{{\cal I}}

\def\coker{{\rm coker}}

\def\CC{ {\cal C}}

\def\CH{ {\cal H}}
\def\CP{ {\cal P}}

\def\nAM0{{\nu }_{{\cal A},M_0}}
\def\nAN0{{\nu }_{{\cal A},N_0}}

\def\End{ {\rm End }}

\def\CP{ {\cal P }}
\def\det{ {\rm det }}
\def\ad{ {\rm ad }}

\def\bx{\overline{x}}
\def\by{\overline{y}}

\def\ga{\mathfrak{a}}
\def\gb{\mathfrak{b}}

\def\gp{\mathfrak{p}}
\def\gq{\mathfrak{q}}

\def\GL{{\rm GL}}
\def\SL{{\rm SL}}

\def\Spec{{\rm Spec}}

\def\Hom{{\rm Hom}}

\def\di!{\frac{\der^i}{i!}}
\def\dik!{\frac{\der^k_i}{k!}}

\def\id{{\rm id}}

\def\gldim{{\rm gldim}}

\def\Fun{{\rm Fun}}

\def\N{\mathbb{N}}

\def\0{\overline{0}}
\def\1{\overline{1}}

\def\Ln1{\L_{n,\overline{1}}}

\def\oa{\overline{a}}

\def\a1{a_{\overline{1}}}

\def\St{{\rm St}}
\def\S{\Sigma}

\def\vn1{\overrightarrow{n-1}}

\def\Sh{{\rm Sh}}

\def\hx{\widehat{x}}

\def\im{{\rm im}}

\def\mA{\mathbb{A}}
\def\mD{\mathbb{D}}
\def\mL{\mathbb{L}}

\def\Sub{{\rm Sub}}
\def\SSub{{\rm SSub}}
\def\Inc{{\rm Inc}}
\def\Min{{\rm Min}}

\def\Inn{{\rm Inn}}

\def\mS{\mathbb{S}}

\def\lann{{\rm l.ann}}
\def\rann{{\rm r.ann}}
\def\Cen{{\rm Cen}}
\def\clKdim{{\rm cl.Kdim}}

\def\hht{{\rm ht}}

\def\heta{\widehat{\eta}}
\def\mF{\mathbb{F}}

\def\mT{\mathbb{T}}
\def\ind{{\rm ind}}

\def\mG{\mathbb{G}}
\def\mE{\mathbb{E}}
\def\mX{\mathbb{X}}
\def\mU{\mathbb{U}}
\def\Frac{{\rm Frac}}
\def\Out{{\rm Out}}

\begin{document}

\author{V. V. \  Bavula 
}

\title{The group of automorphisms of the Jacobian algebra $\mA_n$ }

\maketitle

\begin{abstract}
The {\em Jacobian} algebra $\mA_n$ is obtained from the Weyl
algebra $A_n$ by inverting (not in the sense of Ore) of certain
elements ($\mA_n$  is neither a Noetherian algebra nor a domain,
$\mA_n$ contains the algebra $K\langle x_1, \ldots , x_n,
\frac{\der}{\der x_1}, \ldots ,\frac{\der}{\der x_n},  \int_1,
\ldots , \int_n\rangle $ of  polynomial  integro-differential
operators). The group of automorphisms $\mG_n$ of the Jacobian
algebra $\mA_n$ is found ($\mG_n$ is a huge group):
$$ \mG_n=S_n\ltimes (\mT^n\times \Xi_n)\ltimes  \Inn (\mA_n) \supseteq
 S_n\ltimes (\mT^n\times (\Z^n)^{(\Z )})\ltimes \underbrace{\GL_\infty (K)\ltimes\cdots \ltimes
\GL_\infty (K)}_{2^n-1 \;\; {\rm times}}, $$
$$ \mG_1\simeq  (\mT^1\times \Z^{(\Z )} ) \ltimes \GL_\infty (K),$$
 where  $S_n$
is the symmetric group, $\mT^n$ is the $n$-dimensional algebraic
torus, $\Xi_n\simeq \Z^n$ is a group given explicitly, $\Inn
(\mA_n)$ is the group of inner automorphisms of $\mA_n$ (which is
huge), $\GL_\infty (K)$ is the group of invertible infinite
dimensional matrices, and $(\Z^n)^{(\Z )}$ is a direct sum of $\Z$
copies of the free abelian group $\Z^n$.  This result may help in
understanding of the structure of the groups of automorphisms of
the Weyl algebra $A_n$ and the polynomial algebra $P_{2n}$.
 Explicit generators are found for
the group $\mG_1$. The stabilizers in $\mG_n$ of all the ideals of
$\mA_n$ are found, they are subgroups of {\em finite}  index in
$\mG_n$. It is shown that the group $\mG_n$ has trivial centre. An
explicit inversion formula is given for the elements of $\mG_n$.
Defining relations are found for the algebra $\mA_n$.

 {\em Key Words: the Jacobian algebras,
 the group of automorphisms, the
inner automorphisms, the Fredholm operators,  the index of an
operator, stabilizers, algebraic group, semi-direct product of
groups, the prime spectrum, the minimal primes. }

 {\em Mathematics subject classification
2000: 16W20,  14E07, 14H37, 14R10, 14R15.}

$${\bf Contents}$$
\begin{enumerate}
\item Introduction. \item Preliminaries on the algebras $\mS_n$.
\item The algebras $\mS_n$ and $\mA_n$ are generalized Weyl
algebras.
 \item Certain subgroups of $\Aut_{K-{\rm
alg}}(\mA_n )$. \item The group $\Aut_{K-{\rm alg}}(\mA_1)$. \item
The group of automorphisms of the algebra $\CA_n:= \mA_n/ \ga_n$.
 \item The group of automorphisms of the Jacobian algebra $\mA_n$.
\item Stabilizers in $\Aut_{K-{\rm alg}}(\mA_n )$ of the  ideals
of $\mA_n$.

\end{enumerate}
\end{abstract}


\section{Introduction}
Throughout, ring means an associative ring with $1$; module means
a left module;
 $\N :=\{0, 1, \ldots \}$ is the set of natural numbers; $K$ is a
field of characteristic zero and  $K^*$ is its group of units;
$P_n:= K[x_1, \ldots , x_n]$ is a polynomial algebra over $K$;
$\der_1:=\frac{\der}{\der x_1}, \ldots , \der_n:=\frac{\der}{\der
x_n}$ are the partial derivatives ($K$-linear derivations) of
$P_n$; $\End_K(P_n)$ is the algebra of all $K$-linear maps from
$P_n$ to $P_n$ and $\Aut_K(P_n)$ is its group of units (i.e. the
group of all the invertible linear maps from $P_n$ to $P_n$); the
subalgebra  $A_n:= K \langle x_1, \ldots , x_n , \der_1, \ldots ,
\der_n\rangle$ of $\End_K(P_n)$ is called the $n$'th {\em Weyl}
algebra.

$\noindent $

{\it Definition}: The {\em Jacobian algebra} $\mA_n$ is the
subalgebra of $\End_K(P_n)$ generated by the Weyl algebra $A_n$
and the elements $H_1^{-1}, \ldots , H_n^{-1}\in \End_K(P_n)$
where $$H_1:= \der_1x_1, \ldots , H_n:= \der_nx_n.$$

Clearly, $\mA_n =\bigotimes_{i=1}^n \mA_1(i) \simeq \mA_1^{\t n }$
where $\mA_1(i) := K\langle x_i, \der_i , H_i^{-1}\rangle \simeq
\mA_1$. The algebra $\mA_n$ contains all the  integrations
$\int_i: P_n\ra P_n$, $ p\mapsto \int p \, dx_i$, since  $\int_i=
x_iH_i^{-1}: x^\alpha \mapsto (\alpha_i+1)^{-1}x_ix^\alpha$. In
particular, the algebra $\mA_n$ contains the algebra $K\langle
x_1, \ldots , x_n$,
 $\frac{\der}{\der x_1}, \ldots ,\frac{\der}{\der x_n},  \int_1,
\ldots , \int_n\rangle $ of  polynomial integro-differential
operators. This fact explains why the algebras $\mA_n$ and $A_n$
have different properties and  why the group $\mA_n^*$ of units of
the algebra $\mA_n$ is huge. The algebra $\mA_n$ is neither left
nor right Noetherian. In particular, it is neither left nor right
localization of the Weyl algebra $A_n$ in the sense of Ore. The
algebra $\mA_n$ is not simple, its classical Krull dimension is
$n$, and it contains the algebra of infinite dimensional matrices,
\cite{Bav-Jacalg}. The Jacobian algebra $\mA_n$ appeared in my
study of the group of polynomial automorphisms and the Jacobian
Conjecture, which is a conjecture that makes sense {\em only} for
 the polynomial algebras in the class of all the commutative
 algebras, this was proved in
\cite{Bav-inform}. In order to solve the Jacobian Conjecture,   it
is reasonable to believe that one should create a technique which
makes sense {\em only} for polynomials;  the Jacobian algebras are
a step in this direction (they exist for polynomials but make no
sense even for Laurent polynomials). The Jacobian algebras were
studied  in \cite{Bav-Jacalg}. A closely related to the algebra
$\mA_n$ is the so-called algebra $\mS_n$ of one-sided inverses of
the polynomial algebra $P_n$. The reason for that is that the
derivation $\der_i$  is a left (but not two-sided) inverse of the
$i$th integration $\int_i$, i.e. $\der_i\int_i = \id_{P_n}$.

$\noindent $

{\it Definition}, \cite{shrekalg}. The 
{\em algebra} $\mathbb{S}_n$ {\em of one-sided inverses} of $P_n$
is an algebra generated over a field $K$  by $2n$ elements $x_1,
\ldots , x_n, y_n, \ldots , y_n$ that satisfy the defining
relations:
$$ y_1x_1=\cdots = y_nx_n=1 , \;\; [x_i, y_j]=[x_i, x_j]= [y_i,y_j]=0
\;\; {\rm for\; all}\; i\neq j,$$ where $[a,b]:= ab-ba$ is  the
algebra  commutator of elements $a$ and $b$. Let $G_n:=
\Aut_{K-{\rm alg}}(\mS_n)$.

$\noindent $

By the very definition, the algebra $\mS_n\simeq \mS_1^{\t n}$ is
obtained from the polynomial algebra $P_n$ by adding commuting,
left (but not two-sided) inverses of its canonical generators. The
algebra $\mS_1$ is a well-known primitive algebra
\cite{Jacobson-StrRing}, p. 35, Example 2. Over the field
 $\mathbb{C}$ of complex numbers, the completion of the algebra
 $\mS_1$ is the {\em Toeplitz algebra} which is the
 $C^*$-algebra generated by a unilateral shift on the
 Hilbert space $l^2(\N )$ (note that $y_1=x_1^*$). The Toeplitz
 algebra is the universal $C^*$-algebra generated by a
 proper isometry.

$\noindent $

{\it Example}, \cite{shrekalg}. Consider a vector space $V=
\bigoplus_{i\in \N}Ke_i$ and two shift operators on $V$, $X:
e_i\mapsto e_{i+1}$ and $Y:e_i\mapsto e_{i-1}$ for all $i\geq 0$
where $e_{-1}:=0$. The subalgebra of $\End_K(V)$ generated by the
operators $X$ and $Y$ is isomorphic to the algebra $\mS_1$
$(X\mapsto x$, $Y\mapsto y)$. By taking the $n$'th tensor power
$V^{\t n }=\bigoplus_{\alpha \in \N^n}Ke_\alpha$ of $V$ we see
that the algebra $\mS_n\simeq \mS_1^{\t n}$ is isomorphic to the
subalgebra of $\End_K(V^{\t n })$ generated by the $2n$ shifts
$X_1, Y_1, \ldots , X_n, Y_n$ that act in `perpendicular
directions.'

$\noindent $

The algebras $\mS_n$ are fundamental non-Noetherian algebras, they
are universal non-Noetherian algebras of their own kind in a
similar way as the polynomial algebras are universal in the class
of all the commutative algebras and the Weyl algebras are
universal in the class of algebras  of differential operators.

The algebra $\mS_n$ often appears as a subalgebra or a factor
algebra  of many non-Noetherian algebras. For example, $\mS_1$ is
a factor algebra of certain non-Noetherian down-up algebras as was
shown by Jordan  \cite{Jordan-Down-up}  (see also Benkart and Roby
\cite{Benk-Roby}; Kirkman, Musson, and Passman
\cite{Kir-Mus-Pas-99}; Kirkman and Kuzmanovich \cite{Kirk-Kuz}).
The Jacobian algebra $\mA_n$ contains the algebra $\mS_n$ where
$$y_1:=H_1^{-1}\der_1, \ldots, y_n:= H_n^{-1}\der_n.$$
Moreover, the algebra $\mA_n$ is the subalgebra of
 $\End_K(P_n)$ generated by the algebra $\mS_n$ and  the
$2n$ invertible elements $H_1^{\pm 1}, \ldots , H_n^{\pm 1}$ of
$\End_K(P_n)$.  The algebra $\mA_n$ contains the subalgebra
$K\langle  \frac{\der}{\der x_1}, \ldots ,\frac{\der}{\der x_n},
\int_1, \ldots , \int_n\rangle $ which is isomorphic to the
algebra $\mS_n$. The Gelfand-Kirillov dimensions of the algebras
$\mS_n$ and $\mA_n$
 are $2n$ and $3n$ respectively, \cite{shrekaut} and
 \cite{Bav-Jacalg}.
\begin{itemize}
\item (Corollary \ref{g21Mar9}) {\em Every automorphism of the
algebra $\mS_n$ can be uniquely extended to an automorphism of the
algebra $\mA_n$.} \item (Proposition \ref{f21Mar9}) 1. {\em The
group $G_n$ is a subgroup of}  $\mG_n$.

2. $G_n=\{ \s \in \mG_n\, | \, \s (\mS_n) = \mS_n\}$.
\end{itemize}
 The group $G_n$ of automorphisms of the
 algebra $\mS_n$ is found in \cite{shrekaut}, it is a huge group.
\begin{itemize}
\item (Theorem 5.1, \cite{shrekaut})  $G_n=S_n\ltimes \mT^n
\ltimes \Inn (\mS_n)$, {\em where $S_n$ is the symmetric group,
$\mT^n$ is the $n$-dimensional algebraic  torus, and $\Inn
(\mS_n)$ is the group of inner automorphisms of the algebra}
$\mS_n$.
\end{itemize}
 The algebras $\mS_n$ were studied
 in \cite{shrekalg} and \cite{shrekaut}. For
 the group $G_n$ explicit generators  are found in \cite{shrekaut}, \cite{K1aut}, and
 \cite{Snaut} respectively for $n=1$, $n=2$, and $n>2$.

Ignoring the non-Noetherian property, the Jacobian algebras
$\mA_n$,
 $\mS_n$, the Weyl algebras $A_n$, and the polynomial algebras
$P_{2n}$ belong to the same class of algebras (see below a reason
why this is the case) -- this is a correct approach for studying
the algebras $\mA_n$ and $\mS_n$.  It is an experimental fact
\cite{shrekalg} that the algebra $\mS_1$ has properties that are a
mixture of the properties of the polynomial algebra $P_2$ in two
variable and the {\em first Weyl } algebra $A_1$, which is not
surprising when we look at their defining relations:
\begin{eqnarray*}
 P_2:& yx-xy=0; \\
 A_1:&yx-xy=1;\\
 \mS_1:& yx=1.
\end{eqnarray*}
The same is true for their higher analogues: $P_{2n}=P_2^{\t n}$,
$A_n:=A_1^{\t n }$ (the $n$'th {\em Weyl } algebra), and $\mS_n =
\mS_1^{\t n }$. For example,
\begin{eqnarray*}
 \clKdim (\mS_n )&\stackrel{\cite{shrekalg}}{=} &2n=\clKdim (P_{2n}), \\
\gldim (\mS_n )&\stackrel{\cite{shrekalg}}{=} &n=\gldim (A_n), \\
\GK (\mS_n )&\stackrel{\cite{shrekalg}}{=} &2n=\GK (A_n) = \GK(P_{2n}), \\
\end{eqnarray*}
where $\clKdim$, $\gldim$, and $\GK$ stand for the classical Krull
dimension, the global homological dimension, and the
Gelfand-Kirillov dimension respectively. In this paper, a reason
is found why the algebras $A_n$, $\mA_n$ and $\mS_n$ are `similar'
-- they are  {\em generalized Weyl algebras} (Lemma \ref{b21Mar9}
and Lemma \ref{a21Mar9}). The big difference between the algebras
$\mA_n$ and  $\mS_n$ on the one side and the algebras $P_{2n}$ and
$A_n$ on the other side  is that the first  are neither left nor
right Noetherian and are not  domains either. This fact can also
be explained using the presentations of the algebras as
generalized Weyl algebras. In contrast to the Weyl algebra $A_n$,
 the defining algebra endomorphisms for the algebras $\mA_n$ and
 $\mS_n$ are {\em not} automorphisms, and as a result these algebras
 acquire some pathological  properties (like not being
 Noetherian).

The aim of this paper is to find the group $\mG_n:=\Aut_{K-{\rm
alg}}(\mA_n)$ of automorphisms of the algebra $\mA_n$.
\begin{itemize}
\item (Theorem \ref{10Apr9}) $ \mG_n=S_n\ltimes (\mT^n\times
\Xi_n)\ltimes \Inn (\mA_n)$. \item (Corollary \ref{d10Apr9})
$\mG_n\supseteq \mG_n':= S_n\ltimes (\mT^n\times \mU_n)\ltimes
\underbrace{\GL_\infty (K)\ltimes\cdots \ltimes \GL_\infty
(K)}_{2^n-1 \;\; {\rm times}}$ {\em and} $\mU_n\simeq (\Z^n)^{(\Z
)}$, \item (Theorem \ref{17Mar9}) $\mG_1\simeq (\mT^1\times
\Z^{(\Z )})\ltimes \GL_\infty (K)$,
\end{itemize}
where  $S_n$ is the symmetric group, $\mT^n$ is the
$n$-dimensional algebraic  torus, $\Xi_n\simeq    \Z^n$ and
$\mU_n$ are groups given explicitly,  $\Inn (\mA_n)$ is the group
of inner automorphisms of the algebra $\mA_n$, and $\GL_\infty
(K)$ is the group of all the invertible infinite dimensional
matrices of the type $1+M_\infty (K)$ where the algebra (without
1) of infinite dimensional matrices $M_\infty (K) :=\varinjlim
M_d(K)=\bigcup_{d\geq 1}M_d(K)$ is the injective limit of the
 matrix algebras. A semi-direct product $H_1\ltimes H_2\ltimes
\cdots \ltimes H_m$ of several groups means that $H_1\ltimes
(H_2\ltimes ( \cdots \ltimes (H_{m-1}\ltimes H_m)\cdots )$. In
particular, we found explicit generators for the group $\mG_1$. In
1968, for the first Weyl algebra $A_1$, explicit generators for
its group of automorphisms were found by Dixmier, \cite{Dix}. For
the higher Weyl algebras $A_n$, to find their groups of
automorphisms and generators is an old open problem.

The proof of Theorem \ref{10Apr9} is rather long (and non-trivial)
and based upon several results proved in this paper (and in
\cite{Bav-Jacalg, shrekalg, shrekaut}) which are interesting on
their own. Let me explain briefly the logical structure of the
proof. There are two cases to consider when $n=1$ and $n>1$. The
proofs of both cases are based on rather  different ideas, and the
first case serves as the base of an induction for the second one.
The case $n=1$ is a kind of a degeneration of the second case and
is much more easier. The key point in finding the group $\mG_1$ is
to use the Fredholm linear maps/operators, their {\em indices},
and the fact that each automorphism of the algebra $\mA_n$ is
determined by its action on the  set $\{ x_1, \ldots , x_n\}$ (or
$\{ y_1, \ldots , y_n\}$):

\begin{itemize}
\item (Theorem \ref{21Mar9}) (Rigidity of the group $\mG_n$)  {\em
Let $\s , \tau \in \mG_n$. Then the following statements are
equivalent.}
\begin{enumerate}
 \item $\s = \tau$. \item $\s (x_1) = \tau (x_1), \ldots , \s (x_n) = \tau
 (x_n)$.
 \item$\s (y_1) = \tau (y_1), \ldots , \s (y_n) = \tau (y_n)$.
\end{enumerate}
\end{itemize}

For $n>1$, one of the key ideas in finding the group $\mG_n$ is to
use an  induction on $n$ and  the fact that the algebra $\mA_n$
has the unique maximal ideal $\ga_n$. The ideal $\ga_n$ is the sum
$\gp_1+\cdots + \gp_n$ of all the height one prime ideals of the
algebra $\mA_n$ (there are exactly $n$ of them, and they are found
in \cite{Bav-Jacalg}). By the uniqueness of $\ga_n$, there is the
natural  group homomorphism
$$ \xi : \mG_n\ra \Aut_{K-{\rm alg}}(\CA_n), \;\; \s \mapsto (\xi (\s ) :a+\ga_n\mapsto \s (a) +\ga_n),  $$
where $\CA_n:= \mA_n/ \ga_n$ is a simple algebra which is
isomorphic to a certain  localization of the Weyl algebra $A_n$
(and has Gelfand-Kirillov dimension $3n$). Note that the
Gelfand-Kirillov dimension of the Weyl algebra $A_n$ is $2n$.
Briefly, the problem of finding the group $\mG_n$ is equivalent to
the problem of finding the kernel and the image of the group
homomorphism $\xi$. As the first step in finding the the image of
$\xi$, the  group $\Aut_{K-{\rm alg}}(\CA_n) $ is  found (Theorem
\ref{5Apr9}), and then using a delicate argument on multiplicities
of the growth of the number of the common eigenvalues of a certain
commuting family of differential operators we find the image of
the homomorphism $\xi$ explicitly (Theorem \ref{6Apr9}.(2), note
that $\im (\xi ) \neq \Aut_{K-{\rm alg}}(\CA_n) $). As a
consequence, we prove that
\begin{itemize}
\item (Theorem \ref{6Apr9}.(1)) $\mG_n = S_n\ltimes(\mT^n\times
\mU_n) \ltimes \ker (\xi)$.
\end{itemize}
Then using two characterizations of the elements of the group
$\mG_n$ via the invertible linear maps on the polynomial algebra
$P_n$ (Corollary \ref{d21Mar9}.(1) and Lemma \ref{a18Apr9}) we
prove that
\begin{itemize}
\item (Theorem \ref{6Apr9}.(3)) $\ker (\xi)\subseteq \Inn
(\mA_n)$.
\end{itemize}
This is the most difficult part  of  the  proof of finding the
group $\mG_n$. From this fact, it is relatively easy to show that
\begin{itemize}
\item (Theorem \ref{B10Apr9}) $\Inn (\mA_n)=\mU^0_n\ltimes \ker
(\xi )$,
\end{itemize}
and then that $\mG_n = S_n\ltimes (\mT^n\times \Xi_n) \ltimes \Inn
(\mA_n)$ (Theorem \ref{10Apr9}).

The structure of the group $\mG_1\simeq (\mT^1\times \Z^{(\Z
)})\ltimes \GL_\infty (K)$ (Theorem \ref{17Mar9})  is yet another
confirmation of `similarity' of the algebras $P_2$, $A_1$,
$\mA_1$,  and $\mS_1$. The groups of automorphisms of the
polynomial algebra $P_2$ and the Weyl algebra $A_1$ were found by
 Jung \cite{jung}, Van der Kulk \cite{kulk}, and Dixmier
\cite{Dix} respectively. These two groups have almost identical
structure, they are `infinite $GL$-groups' in the sense that they
are generated by the  algebraic  torus $\mT^1$ and by the obvious
automorphisms: $x\mapsto x+\l y^i$, $y\mapsto y$; $x\mapsto x$,
$y\mapsto y+ \l x^i$, where $i\in \N$ and $\l \in K$; which are
sort of `elementary infinite dimensional matrices' (i.e. `infinite
dimensional transvections'). The same picture as for the groups
$\mG_1\simeq (\mT^1\times \Z^{(\Z )})\ltimes \GL_\infty (K)$ and
$\Aut_{K-{\rm alg}}(\mS_1)=\mT^1\ltimes \GL_\infty (K)$,
\cite{shrekaut}. In prime characteristic, the group of
automorphism of the Weyl algebra $A_1$ was found by Makar-Limanov
\cite{Mak-LimBSMF84} (see also \cite{A1rescen}  for a different
approach and for further developments).

$\noindent $

The paper is organized as follows. In Section \ref{PREL}, some
useful results from \cite{shrekalg} are collected which are used
later. In Section \ref{ASGWA}, it is proved that the algebras
$\mS_n$ and $\mA_n$ are generalized Weyl algebras. As a result,
defining relations are found for the algebra $\mA_n$.

In Section \ref{ACSOA}, several subgroups of the group $\mG_n$ are
introduced, a useful description (Corollary \ref{d21Mar9}.(1)) of
the group $\mG_n$  is given, and a criterion of equality of two
elements of the group $\mG_n$ is proved (Theorem \ref{21Mar9}).

In Section \ref{ATG1S}, the group $\mG_1$ is found (Theorem
\ref{17Mar9}), and explicit generators for the group $\mG_1$ are
given (Theorem \ref{17Mar9}.(4)).

In Section \ref{ATGRP}, the groups $\Aut_{K-{\rm alg}}(\CA_n)$,
$\Inn (\CA_n)$, and $\Out (\CA_n)$ are found (Theorem
\ref{5Apr9}).

In Section \ref{GGNAN}, the groups $\mG_n$, $\Inn (\mA_n )$, and
$\Out (\mA_n )$ are found (Theorem \ref{10Apr9}, Theorem
\ref{6Apr9}, Theorem \ref{B10Apr9}, Corollary \ref{c10Apr9}).
Several corollaries are obtained. It is proved that the groups
$\mG_n$ and $\ker (\xi) $ have trivial centre (Theorem
\ref{16Apr9}, Theorem \ref{A16Apr9}). Explicit inversion formulae
for the elements of the groups $\mG_n$ and $G_n$ are found,
(\ref{mGinv}) and (\ref{1mGinv}). The groups $\mS_n^*$ and $\Inn
(\mS_n)$ are described (Corollary \ref{a19Apr9}). These
descriptions are instrumental in finding explicit generators for
the groups $G_n$, $\mS_n^*$, and $\Inn (\mS_n)$ in \cite{Snaut}.

The Problem/Conjecture of Dixmier \cite{Dix} states that {\em
every algebra endomorphism of the Weyl algebra $A_n$ is an
automorphism}. The Weyl algebra $A_n$ is a simple algebra  but the
Jacobian algebra $\mA_n$ is not. Moreover, its classical Krull
dimension is $n$ \cite{Bav-Jacalg}. Theorem \ref{A10Apr9} has
flavour of the Problem of Dixmier, it states that {\em no proper
prime factor algebra of $\mA_n$ is embeddable into $\mA_n$}.

In Section \ref{STABAN}, the stabilizers in the group $\mG_n$ of
 all the ideals of the algebra $\mA_n$ are computed (Theorem \ref{15Apr9}). In particular, the
stabilizers of all the prime ideals of $\mA_n$ are found
(Corollary \ref{b15Apr9}.(2)).

The maximal ideal $\ga_n:= \gp_1+\cdots + \gp_n$ of the algebra
$\mA_n$ is of height $n$, \cite{Bav-Jacalg}.

\begin{itemize}
\item (Corollary \ref{b15Apr9}.(3)) {\em The ideal $\ga_n$ is the
only nonzero, prime, $\mG_n$-invariant ideal of the algebra
$\mA_n$.} \item  (Corollary  \ref{b15Apr9}) {\em Let $\gp$ be a
prime ideal of $\mA_n$. Then its stabilizer $\St_{\mG_n}(\gp )$ is
a maximal subgroup of the group $\mG_n$ iff $n>1$ and $\gp$ is of
height 1, and, in this case, $[\mG_n:\St_{\mG_n}(\gp )]=n$. }\item
(Corollary \ref{d15Apr9}) {\em Let $\ga$ be a proper ideal of
$\mA_n$. Then its stabilizer $\St_{\mG_n}(\ga )$ has finite index
in the group $\mG_n$.} \item (Corollary \ref{c15Apr9}) {\em If
$\ga$ is a generic  ideal of $\mA_n$ then its stabilizer is
written via the wreath products of symmetric groups:}
$$  \St_{\mG_n}(\ga )= (S_m\times \prod_{i=1}^t(S_{h_i}\wr
S_{n_i}))\ltimes (\mT^n\times \Xi_n)\ltimes \Inn (\mA_n),  $$ {\em
where $\wr$ stands for the {\em wreath} product of groups.}
\end{itemize}
Corollary \ref{e15Apr9} classifies all the proper
$\mG_n$-invariant  ideals of the algebra $\mA_n$, there are
exactly $n$ of them.


\section{Preliminaries on the 
algebras $\mS_n$}\label{PREL}

In this section, we collect some results without proofs on the
algebras $\mS_n$ from \cite{shrekalg} that will be used in this
paper, their proofs can be found in \cite{shrekalg}.

  Clearly,
$\mathbb{S}_n=\mS_1(1)\t \cdots \t\mS_1(n)\simeq \mathbb{S}_1^{\t
n}$ where $\mS_1(i):=K\langle x_i,y_i \, | \, y_ix_i=1\rangle
\simeq \mS_1$ and $\mS_n=\bigoplus_{\alpha , \beta \in \N^n}
Kx^\alpha y^\beta$ where $x^\alpha := x_1^{\alpha_1} \cdots
x_n^{\alpha_n}$, $\alpha = (\alpha_1, \ldots , \alpha_n)$,
$y^\beta := y_1^{\beta_1} \cdots y_n^{\beta_n}$, $\beta =
(\beta_1,\ldots , \beta_n)$. In particular, the algebra $\mS_n$
contains two polynomial subalgebras $P_n$ and $Q_n:=K[y_1, \ldots
, y_n]$ and is equal,  as a vector space,  to their tensor product
$P_n\t Q_n$. Note that also the Weyl algebra $A_n$ is a tensor
product (as a vector space) $P_n\t K[\der_1, \ldots , \der_n]$ of
its two polynomial subalgebras.

When $n=1$, we usually drop the subscript `1' if this does not
lead to confusion.  So, $\mS_1= K\langle x,y\, | \,
yx=1\rangle=\bigoplus_{i,j\geq 0}Kx^iy^j$. For each natural number
$d\geq 1$, let $M_d(K):=\bigoplus_{i,j=0}^{d-1}KE_{ij}$ be the
algebra of $d$-dimensional matrices where $\{ E_{ij}\}$ are the
matrix units, and let
$$M_\infty (K) :=
\varinjlim M_d(K)=\bigoplus_{i,j\in \N}KE_{ij}$$ be the algebra
(without 1) of infinite dimensional matrices. The algebra  $
M_\infty (K)=\bigoplus_{k\in \Z}M_\infty (K)_k$ is $\Z$-graded
 where $M_\infty (K)_k:=\bigoplus_{i-j=k}KE_{ij}$ ($M_\infty (K)_kM_\infty
 (K)_l\subseteq M_\infty (K)_{k+l}$ for all $k,l\in \Z$). The algebra $\mS_1$ contains
 the ideal $F:=\bigoplus_{i,j\in
\N}KE_{ij}$, where 
\begin{equation}\label{Eijc}
E_{ij}:= x^iy^j-x^{i+1}y^{j+1}, \;\; i,j\geq 0.
\end{equation}
For all natural numbers $i$, $j$, $k$, and $l$,
$E_{ij}E_{kl}=\d_{jk}E_{il}$ where $\d_{jk}$ is the Kronecker
delta function.  The ideal $F$ is an algebra (without 1)
isomorphic to the algebra $M_\infty (K)$ via $E_{ij}\mapsto
E_{ij}$. In particular, the algebra $F=\bigoplus_{k\in \Z}F_{1,k}$
is $\Z$-graded where $F_{1,k}:=\bigoplus_{i-j=k}KE_{ij}$
($F_{1,k}F_{1,l}\subseteq F_{1,k+l}$ for all $k,l\in \Z $). For
all $i,j\geq 0$, 
\begin{equation}\label{xyEij}
xE_{ij}=E_{i+1, j}, \;\; yE_{ij} = E_{i-1, j}\;\;\; (E_{-1,j}:=0),
\end{equation}
\begin{equation}\label{xyEij1}
E_{ij}x=E_{i, j-1}, \;\; E_{ij}y = E_{i, j+1} \;\;\;
(E_{i,-1}:=0),
\end{equation}
\begin{equation}\label{xyEij2}
E_{i+1, j+1}x=xE_{i, j}, \;\; E_{ij}y = yE_{i+1, j+1}.
\end{equation}
\begin{equation}\label{mS1d}
\mS_1= K\oplus xK[x]\oplus yK[y]\oplus F,
\end{equation}
the direct sum of vector spaces. Then 
\begin{equation}\label{mS1d1}
\mS_1/F\simeq K[x,x^{-1}]=:L_1, \;\; x\mapsto x, \;\; y \mapsto
x^{-1},
\end{equation}
since $yx=1$, $xy=1-E_{00}$ and $E_{00}\in F$.

$\noindent $

The algebra $\mS_n = \bigotimes_{i=1}^n \mS_1(i)$ contains the
ideal
$$F_n:= F^{\t n }=\bigoplus_{\alpha , \beta \in
\N^n}KE_{\alpha \beta}, \;\; {\rm where}\;\; E_{\alpha
\beta}:=\prod_{i=1}^n E_{\alpha_i \beta_i}(i), \; E_{\alpha_i
\beta_i}(i):= x_i^{\alpha_i}y_i^{\beta_i}-
x_i^{\alpha_i+1}y_i^{\beta_i+1}.$$ Note that $E_{\alpha
\beta}E_{\g \rho}=\d_{\beta \g }E_{\alpha  \rho}$ for all elements
$\alpha, \beta , \g , \rho \in \N^n$ where $\d_{\beta
 \g }$ is the Kronecker delta function.

\begin{proposition}\label{a19Dec8}
\cite{shrekalg} The polynomial algebra $P_n$
 is the only faithful, simple $\mS_n$-module.
\end{proposition}

In more detail, ${}_{\mS_n}P_n\simeq \mS_n / (\sum_{i=0}^n \mS_n
y_i) =\bigoplus_{\alpha \in \N^n} Kx^\alpha \overline{1}$,
$\overline{1}:= 1+\sum_{i=1}^n \mS_ny_i$; and the action of the
elements $E_{\beta \g}$ and the  canonical generators of the
algebra $\mS_n$ on the polynomial algebra $P_n$ is given by the
rule:
$$  x_i*x^\alpha = x^{\alpha + e_i},
\;\; y_i*x^\alpha =
\begin{cases}
x^{\alpha - e_i}& \text{if } \; \alpha_i>0,\\
0& \text{if }\; \alpha_i=0,\\
\end{cases} \;\; {\rm and }\;\; E_{\beta \g }*x^\alpha = \d_{\g \alpha
} x^\beta,
$$
where $e_1:= (1,0,\ldots , 0),  \ldots , e_n:=(0, \ldots , 0,1)$
is the canonical basis for the free $\Z$-module
$\Z^n=\bigoplus_{i=1}^n \Z e_i$.  We identify the algebra $\mS_n$
with its image in the algebra $\End_K(P_n)$ of all the $K$-linear
maps from the vector space $P_n$ to itself, i.e. $\mS_n \subset
\End_K(P_n)$.

$\noindent $

{\bf The involution $\eta$ on $\mS_n$}. The algebra $\mS_n$ admits
the {\em involution}
$$ \eta : \mS_n\ra \mS_n, \;\; x_i\mapsto y_i, \;\; y_i\mapsto
x_i, \;\; i=1, \ldots , n,$$ i.e. it is a $K$-algebra
anti-isomorphism ($\eta (ab) = \eta (b) \eta (a)$ for all $a,b\in
\mS_n$) such that $\eta^2 = \id_{\mS_n}$, the identity map on
$\mS_n$. So, the algebra $\mS_n$ is {\em self-dual} (i.e. it is
isomorphic to its opposite algebra, $\eta : \mS_n\simeq
\mS_n^{op}$). The involution $\eta$ acts on the `matrix' ring
$F_n$ as the transposition,  
\begin{equation}\label{eEij1}
\eta (E_{\alpha \beta} )=E_{\beta \alpha}.
\end{equation}
 The canonical generators $x_i$,
$y_j$ $(1\leq i,j\leq n)$ determine the ascending filtration $\{
\mS_{n, \leq i}\}_{i\in \N}$ on the algebra $\mS_n$ in the obvious
way (i.e. by the total degree of the generators): $\mS_{n, \leq
i}:= \bigoplus_{|\alpha |+|\beta |\leq i} Kx^\alpha y^\beta$ where
$|\alpha | \; = \alpha_1+\cdots + \alpha_n$ ($\mS_{n, \leq
i}\mS_{n, \leq j}\subseteq \mS_{n, \leq i+j}$ for all $i,j\geq
0$). Then $\dim (\mS_{n,\leq i})={i+2n \choose 2n}$ for $i\geq 0$,
and so the Gelfand-Kirillov dimension $\GK (\mS_n )$ of the
algebra $\mS_n$ is equal to $2n$. It is not difficult to show
 that the algebra $\mS_n$ is neither left nor
right Noetherian. Moreover, it contains infinite direct sums of
left and right ideals (see \cite{shrekalg}).

\begin{itemize}
 \item   {\em The
 algebra $\mS_n$ is central, prime, and catenary. Every nonzero
 ideal  of $\mS_n$ is an essential left and right submodule of}
 $\mS_n$.
 \item  {\em The ideals of
 $\mS_n$ commute ($IJ=JI$);  and the set of ideals of $\mS_n$ satisfy the a.c.c..}
 \item  {\em The classical Krull dimension $\clKdim (\mS_n)$ of $\mS_n$ is $2n$.}
  \item  {\em Let $I$
be an ideal of $\mS_n$. Then the factor algebra $\mS_n / I$ is
left (or right) Noetherian iff the ideal $I$ contains all the
height one primes of $\mS_n$.}
\end{itemize}

{\bf The set of height 1 primes of $\mS_n$}.   Consider the ideals
of the algebra $\mS_n$:
$$\gp_1:=F\t \mS_{n-1},\; \gp_2:= \mS_1\t F\t \mS_{n-2}, \ldots ,
 \gp_n:= \mS_{n-1} \t F.$$ Then $\mS_n/\gp_i\simeq
\mS_{n-1}\t (\mS_1/F) \simeq  \mS_{n-1}\t K[x_i, x_i^{-1}]$ and
$\bigcap_{i=1}^n \gp_i = \prod_{i=1}^n \gp_i =F_n$. Clearly,
$\gp_i\not\subseteq \gp_j$ for all $i\neq j$.

\begin{itemize}
 \item
{\em The set $\CH_1$ of height 1 prime ideals of the algebra
$\mS_n$ is} $\{ \gp_1, \ldots , \gp_n\}$.
\end{itemize}

 Let
$\ga_n:= \gp_1+\cdots +\gp_n$. Then the factor algebra
\begin{equation}\label{SnSn}
\mS_n/ \ga_n\simeq (\mS_1/F)^{\t n } \simeq \bigotimes_{i=1}^n
K[x_i, x_i^{-1}]= K[x_1, x_1^{-1}, \ldots , x_n, x_n^{-1}]=:L_n
\end{equation}
is a  Laurent polynomial algebra in $n$ variables,  and so $\ga_n$
is a prime ideal of height and co-height $n$ of the algebra
$\mS_n$. The algebra $L_n$ is commutative, and so
\begin{equation}\label{comab}
 [a,b]\in \ga_n\;\; {\rm for \; all}\;\; a,b\in \mS_n .
\end{equation}
That is $[\mS_n , \mS_n]\subseteq \ga_n$. In particular, $[\mS_1,
\mS_1]\subseteq F$. Since $\eta (\ga_n) = \ga_n$, the involution
of the algebra $\mS_n$ induces the {\em automorphism}
$\overline{\eta}$ of the factor algebra $\mS_n / \ga_n$ by the
rule: 
\begin{equation}\label{1SnSn}
\overline{\eta}: L_n\ra L_n, \;\; x_i\mapsto x_i^{-1}, \;\; i=1,
\ldots , n.
\end{equation}
It follows that $\eta (ab) - \eta (a) \eta (b)\in \ga_n$ for all
elements $a, b\in \mS_n$.


\section{The algebras $\mS_n$ and $\mA_n$ are generalized Weyl algebras}\label{ASGWA}
In this section, we show that the algebras $\mS_n$ and $\mA_n$ are
generalized Weyl algebras (Lemma \ref{a21Mar9} and Lemma
\ref{b21Mar9}). As a result, we have explicit defining relations
for the algebra $\mA_n$, they are used in checking of existence of
various automorphisms and anti-automorphisms.  This fact explains
why the algebras $\mS_n$, $\mA_n$, and $A_n$ have `similar'
properties. We start by recalling some properties of generalized
Weyl algebras.

$\noindent $

{\bf Generalized Weyl Algebras}. Let $D$ be a ring, $\s
=(\s_1,...,\s_n)$ be
 an $n$-tuple  of  commuting ring endomorphisms of $D$,  and
$a=(a_1,...,a_n)$ be  an $n$-tuple  of  elements
 of $D$. 
  The {\it generalized Weyl algebra}
$A=D(\s,a)$ (briefly,  GWA) of degree $n$ is  a ring generated by
$D$  and    $2n$ elements $x_1,...,x_n,$ $y_1,...,$ $y_n$ subject
to the defining relations \cite{Bav-FA91}, \cite{Bav-AlgAnaliz92}:
\begin{align*}
y_ix_i&=a_i ,& x_iy_i&=\s_i(a_i), \\
\s_i(d)x_i &=x_i d, & d y_i&=y_i\s_i(d), \;\;\;  d \in D,
\end{align*}
$$[x_i,x_j]=[y_i,y_j]=[x_i,y_j]=0, \;\;\;  i\neq j,$$
where $[x, y]=xy-yx$.
  We say that  $a$  and $\s $ are the sets
of {\it defining } elements and endomorphisms of $A$ respectively.
 When all $\s_i$ are automorphisms, the GWAs are also known as {\em
hyperbolic rings}, see the book of Rosenberg \cite{Rosen-book}.
 For a
vector $k=(k_1,...,k_n)\in \mathbb{Z}^n$, let
$v_k=v_{k_1}(1)\cdots v_{k_n}(n)$ where, for $1\leq i\leq n$ and
$m\geq 0$: $\,v_{m}(i)=x_i^m$, $\,v_{-m}(i)=y_i^m$, $v_0(i)=1$. It
follows from  the  definition that $A=\bigoplus_{k\in
\mathbb{Z}^n} A_k$ is a $\mathbb{Z}^n$-graded ring
($A_kA_e\subseteq A_{k+e},$ for all  $k,e \in \mathbb{Z}^n$),
 where $A_k=\bigoplus_{l,k+l\in \N^n} v_{-l}Dv_{k +l}$.
 The tensor product (over the ground  field) $A\t A'$ of generalized Weyl algebras
 of degree $n$ and $n'$ respectively is a GWA of degree $n+n'$:
$$A\otimes A'=D\otimes D'((\tau ,\tau'), (a, a')).$$

Let $\CP_n$ be a polynomial algebra  $K[H_1, \ldots , H_n]$ in $n$
indeterminates  and let $\s =(\s_1,...,\s_n)$ be the  $n$-tuple of
commuting automorphisms of $\CP_n$ such that
 $\s_i(H_i)=H_i-1$ and $\s_i(H_j)=H_j$, for $i\neq j$.  The algebra homomorphism
\begin{equation}\label{AnGWA}
A_n\ra \CP_n ((\s_1,...,\s_n), (H_1, \ldots , H_n)), \;\;
x_i\mapsto  x_i, \; \; \der_i\mapsto y_i, \; \; i=1, \ldots , n,
\end{equation}
  is an isomorphism. We identify the Weyl algebra $A_n$
with the GWA above via this isomorphism. Note that $H_i= \der_i
x_i = x_i\der_i+1$. Denote by $S_n$ the multiplicative submonoid
of $\CP_n$ generated by the elements $H_i+j$, $i=1, \ldots , n$,
and $j\in \mathbb{Z}$. It follows  from the above  presentation of
the Weyl algebra $A_n$ as a GWA that $S_n$ is an Ore set in $A_n$,
and, using the $\mathbb{Z}^n$-grading, that the (two-sided)
localization $\CA_n:=S_n^{-1}A_n$ of the Weyl algebra $A_n$ at
$S_n$ is the {\em skew Laurent polynomial ring}
\begin{equation}\label{Anskewlaurent}
\CA_n=S_n^{-1}\CP_n [x_1^{\pm 1}, \ldots ,x_n^{\pm 1};
\s_1,...,\s_n]
\end{equation}
with coefficients from the algebra
$$ \CL_n:=S_n^{-1}\CP_n=K[H_1^{\pm 1}, (H_1 \pm 1)^{-1}, (H_1 \pm 2)^{-1}, \ldots ,
H_n^{\pm 1}, (H_n \pm 1)^{-1}, (H_n \pm 2)^{-1}, \ldots ],$$ which
is the localization of $\CP_n$ at $S_n$. We identify the Weyl
algebra $A_n$ with the subalgebra of $\CA_n$ via the monomorphism,
$$A_n\ra \CA_n, \;\; x_i\mapsto x_i,\;\; \der_i\mapsto  H_ix_i^{-1}, \;\;
i=1, \ldots , n.$$ Let $k_n$ be the $n$'th {\em Weyl skew field},
that is the full ring of quotients of the $n$'th Weyl algebra
$A_n$ (it exists by  Goldie's Theorem  since $A_n$ is a Noetherian
domain). Then the algebra $\CA_n$ is a $K$-subalgebra of $k_n$
generated by the elements $x_i$, $x_i^{-1}$, $H_i$ and $H_i^{-1}$,
$i=1, \ldots , n$ since, for all  $j\in \N$, 
\begin{equation}\label{Hijxi}
 (H_i\mp j)^{-1}=x_i^{\pm j}H_i^{-1}x_i^{\mp j}, \;\; i=1, \ldots
, n.
\end{equation}
Clearly, $\CA_n \simeq \CA_1 \t \cdots \t \CA_1$ ($n$ times).
 A $K$-algebra $R$ has the {\em endomorphism
property over} $K$ if, for each simple $R$-module $M$, $\End_R(M)$
is algebraic over $K$.

\begin{theorem}\label{bigAn}
\cite{Bav-Ann2001} Let $K$ be a field of characteristic zero.
\begin{enumerate}
\item The algebra $\CA_n$ is a simple, affine, Noetherian domain.
\item The Gelfand-Kirillov dimension $\GK (\CA_n)=3n$ $(\neq
2n=\GK (A_n))$. \item The (left and right) global dimension ${\rm
gl.dim} (\CA_n)=n$. \item The (left and right) Krull dimension
${\rm K.dim} (\CA_n)=n$. \item Let ${\rm d}={\rm gl.dim}$ or ${\rm
d}= {\rm K.dim}$. Let $R$  be a Noetherian $K$-algebra with ${\rm
d}(R)<\infty $ such that $R[t]$, the polynomial ring in a central
 indeterminate, has the endomorphism property over $K$. Then
${\rm d}(\CA_1\t R)= {\rm d}(R)+1$. If, in addition, the field $K$
is algebraically closed and uncountable, and the algebra $R$ is
affine, then
 ${\rm d}(\CA_n\t R)= {\rm d}(R)+n$.
\end{enumerate}
\end{theorem}

$\GK ({\cal A}_1)=3$ is due to A. Joseph \cite{Jos1}, p. 336;
 see also  \cite{KL}, Example 4.11, p. 45.

It is an experimental fact that many small quantum groups are
GWAs. For more about GWAs and their generalizations the interested
reader can refer to \cite{Alev-Far-La-Sol-00, Artam-Cohn-99,
Bav-inform, Bav-Jor-01, Bav-Len-JA-2001-KdimGWA,
Bav-Len-JA-2001-TKMin, Bav-Oyst-1998KdimGWA, Bav-Oyst-2000Adv,
Bav-Oyst-2001TrAMS, Bek-Ben-Fut-04, Far-Sol-SuAl-03, Hart-06,
Kir-Mus-Pas-99, Kir-Kuz-05,  Maz-Tur-02, Prest-Pu-02, Rich-Sol-06,
Staf-sl2-82}.

$\noindent $

Suppose that $A$ is a $K$-algebra that admits two elements $x$ and
$y$ with $yx=1$. The element $xy\in A$ is an idempotent, $(xy)^2=
xy$, and so the set $xyA xy$ is  a $K$-algebra where $xy$ is its
identity element. Consider the linear maps $\s = \s_{x,y}, \tau =
\tau_{x,y}:A \ra A$ which are defined as follows 
\begin{equation}\label{sxyt}
\s (a) =  xay, \;\; \tau (a) = yax.
\end{equation}
Then $\tau \s = \id_A$ and $\s \tau  (a) = xy \cdot a \cdot xy$,
and so 
\begin{equation}\label{Asrt}
A= \s (A) \bigoplus \ker (\tau ).
\end{equation}
Each element $a\in A$ is a unique sum $a_1+a_2$ where $a_1\in \s
(A)$ and $a_2\in \ker (\tau )$. Below, formulae for the elements
$a_1$ and $a_2$ are given. The map $\s : A\ra A$ is an algebra
{\em monomorphism} ($\s (ab) = xaby = xay \cdot xby= \s (a) \s
(b)$) with $\s (1) = xy$. The element $\s \tau \in \End_K(A)$ is
an idempotent, then so is $\id_A-\s \tau$. Since $\s \tau (\ker
(\tau ))=0$ and $\s \tau \s (a) = \s (a)$, we have $\s (A) = \im
(\s \tau )$ and $\ker (\tau ) = \im (\id_A-\s \tau )$. For each
element $a\in A$,
$$ a= \s \tau (a) +(\id_A - \s \tau ) (a) \;\;{\rm where} \;\; \s
\tau (a) \in \s (A), \;\; (\id_A - \s \tau ) (a)\in \ker (\tau
).$$ Note that $xy A xy \subseteq \s (A) = \s (1\cdot A \cdot 1) =
xy \s (A) xy \subseteq xy A xy$, and so $\s (A) = xy A xy$. The
map $\s : A\ra \s (A)= xy A xy$ is an isomorphism of algebras
where $\s(1) = xy$ and its inverse map $\tau : \s (A) \ra A$ is an
algebra isomorphism with $\tau (xy ) = 1$.

Suppose that the algebra $A$ contains a subalgebra $D$ such that
$\s (D) \subseteq D$ and $\tau (D)\subseteq D$ (and so $xy = \s
(1) \in D$), and that the algebra $A$ is generated by $D$, $x$,
and $y$. Since $yx=1$, we have $x^iD_D\simeq D$ and
${}_DDy^i\simeq D$. It follows from the relations:
\begin{align*}
yx&=1 ,& xy&=\s (1), \\
xd  &=\s (d) x, & d y&=y \s (d), \;\;\;  d \in D,
\end{align*}
that $A= \sum_{i\geq 1} y^iD+\sum_{i\geq 0}Dx^i$. Suppose, in
addition, that the sum is a {\em direct} one. Then the algebra $A$
is the GWA $D(\s , 1)$.

\begin{lemma}\label{a13Apr9}
Keep the assumptions as above, i.e. $A=D\langle x,y\rangle =
\bigoplus_{i\geq 1} y^iD\bigoplus\bigoplus_{i\geq 0}Dx^i$, $\tau
(D)\subseteq D$  and $\s (D)\subseteq D$. Then $A= D(\s , 1)$. If,
in addition, the element $xy$ is central in $D$. Then $Dx^i =
x^iD$ and $Dy^i = y^iD$ for all $i\geq 1$.
\end{lemma}

{\it Proof}. It suffices to prove the equalities for $i=1$. $xD
=\s (D)x \subseteq Dx= Dxy x =  xyDx = x \tau (D) \subseteq x D$
and so $xD = Dx$. Similarly, $Dy = y \s(D) \subseteq y D = y xyD=
 yDxy = \tau (D) y \subseteq Dy$, and so $Dy = yD$.
$\Box $

$\noindent $

{\bf The algebras $\mS_n$ are generalized Weyl algebras}. Consider
the case $n=1$. The subalgebra $\mF_1:= K+F=\bigoplus_{i\in
\Z}\mF_{1,i}$ of $\End_K(P_1)$ is a $\Z$-graded algebra where
$\mF_{1,0}:= K+F_{1,0}= K\bigoplus \bigoplus_{i\in \N}E_{ii}$ and
$\mF_{1,k}:= F_{1,k}=\bigoplus_{i-j=k}KE_{ij}$ for $k\in \Z
\backslash \{ 0\}$. Moreover, $\mF_{1,k}\mF_{1,l}\subseteq
\mF_{1,k+l}$ for all $k,l\in \Z$. By (\ref{xyEij}) and
(\ref{xyEij1}),  $\s (\mF_{1,0}) \subseteq \mF_{1,0}$ where $\s $
is as in (\ref{sxyt}).  The $K$-algebra monomorphism
\begin{equation}\label{tF10}
\s : \mF_{1,0}\ra \mF_{1,0}, \;\; E_{ii}\mapsto E_{i+1, i+1},  \;
\; i\geq 0, \;\; 1\mapsto xy=1-E_{00},
\end{equation}
is not an automorphism since $\mF_{1,0}=\im (\s ) \oplus
KE_{00}$. The algebra $\mF_{1,0}$ is a commutative algebra which
is generated by the elements $E_{ii}$, $i\geq 0$, that satisfy the
defining relations: $E_{ii}^2=E_{ii}$ and $E_{ii}E_{jj}=0$ for all
$i\neq j$. For each $i\in \N$, $KE_{ii}$ is an ideal of the
algebra $\mF_{1,0}$. The algebra $\mF_{1,0}$ is not Noetherian
since it contains the infinite direct sum of ideals
$\bigoplus_{i\in \N}KE_{ii}$.

 Note that  $yx=1$, $xy=\s (1)= 1-E_{00}$,
 $xE_{ii}=E_{i+1, i+1}x= \s(E_{ii})x$, and $E_{ii}y= yE_{i+1,
 i+1}= y\s (E_{ii})$ for all $i\in \N$. Now, it follows from
 (\ref{xyEij}), (\ref{xyEij1}), and (\ref{mS1d}) that the assumptions of Lemma \ref{a13Apr9} hold, and so the algebra
 $\mS_1$ is the GWA $\mS_1= \mF_{1,0}(\s , 1)$. In
 particular, the algebra $\mS_1=\bigoplus_{i\in \Z}\mS_{1,i}$ is
 $\Z$-graded where
$$
\mS_{1,i}=\begin{cases}
x^i\mF_{1,0}=\mF_{1,0}x^i& \text{if $i\geq 1$},\\
\mF_{1,0}& \text{if $i=0$},\\
 y^{-i}\mF_{1,0}=\mF_{1,0} y^{-i}& \text{if $i\leq -1$}.\\
\end{cases}
$$
The kernel of the map $\tau : \mS_1\ra \mS_1$, $a\mapsto yax$, is
\begin{equation}\label{kertau}
\ker (\tau ) = KE_{00}+\sum_{i\geq 1}(KE_{0i}+KE_{i0}).
\end{equation}
Indeed, since $\tau (F) \subseteq F$ and the induced map
$\overline{\tau}:\mS_1/F\ra \mS_1/F$ is the identity map, the
kernel of the map $\tau$ coincides with the kernel of its
restriction to $F$, and so (\ref{kertau}) is obvious since $\tau
(E_{ij} ) = E_{i-1, j-1}$. The kernel of the map $\tau :\mS_1\ra
\mS_1$ is not an ideal of the algebra $\mS_1$. This means that the
map $\tau $ is not an algebra endomorphism but its restriction to
the subalgebra $\mF_{1,0}$ of $\mS_1$ {\em is} a $K$-algebra
epimorphism: 
\begin{equation}\label{tF1}
\tau :\mF_{1,0}\ra \mF_{1,0}, \;\; E_{ii}\mapsto E_{i-1, i-1}\;\;
(i\geq 0),\;\; 1\mapsto 1,
\end{equation}
with $\ker (\tau |_{\mF_{1,0}})=KE_{00}$. For all $i\in \N$ and
$d\in \mF_{1,0}$, $dx^i=x^i\tau^i(d)$ and $ y^id=\tau^i(d)y^i$.

 For $n\geq 1$ and $i=1, \ldots
, n$, $\mS_1(i) = \mF_{1,0}(i)(\s_i ,
 1)$, and
 $ \mS_n=\bigotimes_{i=1}^n \mS_1(i) = \bigotimes_{i=1}^n \mF_{1,0}(i)(\s_i ,
 1)$ is the GWA
 $$ \mS_n = \mF_{n,0} ((\s_1, \ldots , \s_n), (1,
 \ldots , 1)), $$
 where $\mF_{n, 0}:= \bigotimes_{i=1}^n \mF_{1,0}(i)$ and
 $\mF_{1,0}(i):= K+F_{1,0}(i) = K\bigoplus\bigoplus_{k\geq
 0}E_{kk}(i)$.

\begin{lemma}\label{a21Mar9}
The algebra $\mS_n =\mF_{n,0} ((\s_1, \ldots , \s_n), (1,
 \ldots , 1))$ is a generalized Weyl algebra. In particular, the algebra
 $\mS_n =\bigoplus_{\alpha \in \Z^n}\mS_{n,\alpha}$ is $\Z^n$-graded  where
 $\mS_{n, \alpha} = \mF_{n,0}v_\alpha =
 v_\alpha \mF_{n, 0}$  for all  $\alpha \in \Z^n$.
\end{lemma}

For each algebra $\mS_1(i)$, let $\tau_i$ be the corresponding map
$\tau$ which is extended to the map $\tau_i:\mS_n\ra \mS_n$,
$a\mapsto y_iax_i$. It is not an algebra endomorphism but its
restriction to the subalgebra $\mF_{n,0}$ of $\mS_n$ {\em is} a
$K$-algebra epimorphism, $\tau_i(\mF_{n,0})=\mF_{n,0}$, with $\ker
(\tau_i |_{\mF_{n,0}})=KE_{00}(i)\bigotimes \bigotimes_{j\neq
i}\mF_{1,0}(j)$. For all $j\in \N$ and $d\in \mF_{n,0}$,
$dx_i^j=x_i^j\tau_i^j(d)$ and $ y_i^jd=\tau_i^j(d)y_i^j$.

$\noindent $

{\bf The Jacobian algebras $\mA_n$ are generalized Weyl algebras}.
Consider the case $n=1$ and the commutative subalgebra $C_1:=
\CP_1\bigoplus F_{1,0}= K[H]\bigoplus \bigoplus_{i\geq 0}KE_{ii}$
of the algebra $\End_K(P_1)$. Then $\mF_{1,0}=K+F_{1,0}\subset
C_1$. Then $C_1/F_{1,0}\simeq \CP_1=K[H]$. As an abstract algebra,
the commutative  algebra $C_1$ is generated by the elements $H$,
$E_{ii}$, $i\in \N$, that satisfy the defining relations:
$$ HE_{ii}= (i+1)E_{ii}, \;\; E_{ii}^2=E_{ii}, \;\;
E_{ii}E_{jj}=0\;\; {\rm for\; all}\;\; i\neq j.$$
For each natural number $i$, the element $H+i$ is a unit of the
algebra $\End_K(P_1)$ but the element $H-i$, $i\geq 1$, is not
since $\ker_{P_1}(H-i) = Kx^{i-1}$. For each scalar $\l \in K^*$
and for each integer $i\geq 1$, the element 
\begin{equation}\label{Hil}
 (H-i)_\l := H-i+\l E_{i-1, i-1}\in \mA_1
\end{equation}
 is a unit of the
algebra $\End_K(P_1)$ since $(H-i)_\l *x^{i-1}= \l x^{i-1}\neq 0$.
Let $S_1'$ be the multiplicative submonoid of $C_1$ generated by
the elements $\{ (H-i)_1, H+j\, | \, i\geq 1, j\geq 0\}$. The
localization $\mD_1:= S_1'^{-1}C_1$ of the algebra $C_1$ at $S_1'$
contains the algebra $C_1$.  In \cite{Bav-Jacalg}, it is shown
that the algebra $\mA_1$ is generated by $\mD_1$, $x$, and $y$;
$F_{1,0}$ is an ideal of the algebra $\mD_1$ such that
$\mD_1/F_{1,0}\simeq S_1^{-1}\CP_1$; and $F$ is an ideal of the
algebra $\mA_1$ such that $\mA_1/F\simeq \CA_1$.
There is the commutative diagram of natural algebra
homomorphisms with exact rows (of algebras not necessarily with
1):
$$
\xymatrix{0\ar[r] & F_{1,0}\ar[r]\ar[d]^{=}  & C_1 \ar[r]\ar[d]& C_1/F_{1,0} \simeq \CP_1\ar[r]\ar[d] & 0 \\
0\ar[r] & F_{1,0} \ar[r]  & \mD_1 \ar[r] & \mD_1/F_{1,0} \simeq
S_1^{-1}\CP_1\ar[r] & 0 \, . }
$$
 Since $\CL_1:=S_1^{-1}\CP_1= K[H^{\pm 1}, (H\pm
1)^{-1}, (H\pm 2)^{-1}, \ldots ]$ and $(H-i)_1\equiv H-i\mod
F_{1,0}$, it follows from the exact sequence at the bottom of the
commutative diagram above  that 
\begin{equation}\label{CD1d}
\mD_1=\CL_1^-\bigoplus \CL_1^+\bigoplus F_{1,0}, \;\;
\CL_1^-:=\bigoplus_{i,j\geq 1}K\frac{1}{(H-i)_1^j},\;\;
\CL_1^+:=K[H^{\pm 1}, (H+ 1)^{-1},  (H+ 2)^{-1}, \ldots ].
\end{equation}
Consider the maps $\s (a) = xay$ and $\tau (a) = yax$ of the
algebra $\mA_1$ as defined in (\ref{sxyt}). Then $\s (\mD_1) =
\mD_1 \s (1)\subseteq \mD_1$ as follows from the equalities: for
all natural numbers $i\geq 0$ and $j\geq 1$,
$$\s (E_{ii})= E_{i+1, i+1}, \;\; \s (H)=H-1=(H-1)\s (1)
= (H-1)_1\s (1), \;\; \s ((H-j)_1)= (H-j-1)_1\s (1).$$ Moreover,
$\mD_1 = \s (\mD_1) \bigoplus KE_{00}$ where $KE_{00}$ is an ideal
of the commutative algebra $\mD_1$ such that $\tau (KE_{00})=0$.
Then $\tau (\mD_1) = \mD_1$. Moreover, the map $\tau : \mD_1\ra
\mD_1$, $\s (d)+\l E_{00} \mapsto d$, is an {\em algebra
epimorphism} ($\tau (1) = \tau (1-E_{00}+E_{00})= \tau (\s
(1)+E_{00})= 1$) with kernel $KE_{00}$. In particular, $\tau (H)=
H+1$ and $\tau (E_{ii})= E_{i-1, i-1}$. For all $i\in \N$ and
$d\in \mD_1$, $dx^i=x^i\tau^i(d)$ and $y^id=\tau^i(d)y^i$. The
kernel of the
linear map $\tau : \mA_1\ra \mA_1$ is equal to 
\begin{equation}\label{1kertau}
\ker (\tau)= KE_{00}+\sum_{i\geq 1} (KE_{0i}+KE_{i0}).
\end{equation}
In more detail, since $\tau (F)\subseteq  F$ and the induced map
$\overline{\tau}: \mA_1/ F\ra \mA_1/ F$ is the inner automorphism
$a\mapsto x^{-1}ax$, the kernel of the map $\tau$ coincides with
the kernel of its restriction to $F$, and the equality
(\ref{1kertau})
 follows from (\ref{kertau}).
 The kernel of $\tau$  is not an ideal of the algebra $\mA_1$, and so the map
$\tau $ is not an algebra homomorphism, but $\tau$ {\em is} an
algebra homomorphism modulo $F$.

For all integers $i\geq 1$,
$$ \tau ((H-i)_1) = \begin{cases}
(H-(i-1))_1& \text{if } i\geq 2,\\
H& \text{if } i=1,
\end{cases}\;\;  {\rm and}\;\;\frac{1}{(H-i)_1}E_{jj}=\begin{cases}
\frac{1}{j-i+1}E_{jj}& \text{if } j\neq i-1,\\
E_{i-1,i-1}& \text{if } j=i-1.\\
\end{cases} $$
The $K[\tau ]$-module $\mD_1$  is the direct sum of the following
  {\em indecomposable}, infinite dimensional submodules:
$K[H]$; $\bigoplus_{i\geq 1}K\frac{1}{(H-i)_1^j}\bigoplus
\bigoplus_{i\geq 1} K\frac{1}{(H+i)^j}$ for each integer $j\geq
1$; and $F_{1,0}:=\bigoplus_{s\geq 0}KE_{ss}$ for each $k\in \Z$.
The map $\tau : \CL_1^-\bigoplus \CL_1^+\ra \CL_1^-\bigoplus
\CL_1^+$ is a  bijection,  the map $\tau : \mF_{1,0}\ra \mF_{1,0}$
is a surjection with $\ker (\tau ) = KE_{00}$, and the map $\s :
\mF_{1,0}\ra \mF_{1,0}$ is an injection with $\mF_{1,0}= \im (\s
)\oplus  KE_{00}$.

The $K[\s ]$-module structure of the algebra $\mD_1$ is
complicated as it follows from the equalities: for all $\l \in K$
and $i\geq 1$,
$$ \s^i (\l ) = \l ( 1-\sum_{j=0}^{i-1}E_{jj})\;\; {\rm and}\;\;
\s^i (H) = H-i+\sum_{j=0}^{i-1}(i-1-j)E_{jj}.$$ Note that $\s (x)
= xxy = x(1-E_{00}) = x-E_{10}$ and $\s (y) = xyy = (1-E_{00})y =
y-E_{01}$. Then $\s (x^i) = x^i - E_{i0}$ and $\s (y^i) = y^i -
E_{0i}$ for all $i\geq 1$.

The algebra $\mA_1$ is generated by the elements $x$,
$y:=H^{-1}\der$, $H$,  and $H^{-1}$, or, equivalently, the algebra
$\mA_1$ is obtained from the subalgebra $\mS_1=K\langle
x,y\rangle$  of $\End_K(P_1)$ by adding the elements $H$ and
$H^{-1}$ of $\End_K(P_1)$. For each integer $m\geq 1$, $\der^m =
(Hy)^m= H(H+1) \cdots (H+m-1)y^m$, and so $\mD_1\der^m =
\mD_1y^m$. In (\cite{Bav-Jacalg}, Theorem 2.3) it is proved that
the algebra $\mA_n =\bigoplus_{\alpha \in \Z^n}\mA_{n,\alpha}$ is
a $\Z^n$-graded algebra where $\mA_{n,\alpha}:=\bigotimes_{k=1}^n
\mA_{1,\alpha_k}(k)$ and, for $n=1$,

$$
\mA_{1,i}=\begin{cases}
x^i\mD_1& \text{if $i\geq 1$},\\
\mD_1& \text{if $i=0$},\\
\mD_1 \der^{-i}=\mD_1y^{-i}& \text{if $i\leq -1$}.\\
\end{cases}
$$
Since $yx=1$, $xy=\s (1) = 1-E_{00}$,  $xd=\s (d)x$ and $dy=y\s
(d) $ for all $d\in \mD_1$, $\s (\mD_1)\subseteq \mD_1$, $\tau
(\mD_1) \subseteq \mD_1$, and $xy$ is a central element of
$\mD_1$, by Lemma \ref{a13Apr9}, the Jacobian algebra $\mA_1$ is
the GWA $\mA_1= \mD_1(\s , 1)$ such that $x^i \mD_1= \mD_1x^i$ and
$y^i \mD_1= \mD_1y^i$ for all $i\geq 1$. Note that $\s (1) x=x$
and $ y\s (1)=y$.

For an arbitrary $n$ and for $i=1, \ldots , n$, $\mA_1(i) =
\mD_1(i) (\s_i, 1)$, and so the Jacobian algebra
$$ \mA_n=\bigotimes_{i=1}^n \mA_1(i) = \mD_n ((\s_1, \ldots ,
\s_n), (1, \ldots , 1))$$ is the generalized Weyl algebra where
$\mD_n:= \bigotimes_{i=1}^n \mD_1(i)$.

\begin{lemma}\label{b21Mar9}
\begin{enumerate}
\item The Jacobian algebra $ \mA_n=\mD_n ((\s_1, \ldots , \s_n),
(1, \ldots , 1))$ is the generalized Weyl algebra. In particular,
the algebra $\mA_n = \bigoplus_{\alpha \in \Z^n} \mA_{n, \alpha}$
is $\Z^n$-graded where   $\mA_{n, \alpha} = \mD_n v_\alpha =
v_\alpha \mD_n$ for all $\alpha \in \Z^n$. \item The Jacobian
algebra $\mA_n$ is obtained from the subalgebra $\mS_n$ of
$\End_K(P_n)$ by adding $2n$ invertible elements $H_1^{\pm 1},
\ldots , H_n^{\pm 1}$ of $\End_K(P_n)$.
\end{enumerate}
\end{lemma}

Lemma \ref{b21Mar9}.(1) gives defining relations for the Jacobian
algebra $\mA_n$. The algebras $\CL_n^+ :=\bigotimes_{i=1}^n
\CL_1(i)^+$ and $C_n:=\bigotimes_{i=0}^nC_1(i)$ are subalgebras of
$\mD_n$ where $\CL_1(i)^+:=K[H_i^{\pm 1}, (H_i+1)^{-1},
(H_i+2)^{-1},\ldots ]$ and $C_1(i):= K[H_i]+F_{1,0}(i) =
K[H_i]\bigoplus \bigoplus_{j\geq 0}KE_{jj}(i)$.

For each algebra $\mA_1(i)$, let $\tau_i$ be the corresponding map
$\tau$ which is extended to the map $\tau_i:\mA_n\ra \mA_n$,
$a\mapsto y_iax_i$. It is not an algebra endomorphism but its
restriction to the subalgebra $\mD_n$ of $\mA_n$ {\em is} a
$K$-algebra epimorphism, $\tau_i(\mD_n)=\mD_n$, with $\ker (\tau_i
|_{\mD_n})=KE_{00}(i)\bigotimes\bigotimes_{j\neq i}\mD_1(j)$. In
more detail, for $a,b\in \mD_1$,
$$\tau (a) \tau (b) = ya(1-E_{00})bx= \tau (ab)-yaE_{00}bx= \tau
(ab ),$$ since $aE_{00}b \in KE_{00}$ and $yE_{00}=0$. For all
$j\in \N$ and $d\in \mD_n$, $dx_i^j=x_i^j\tau_i^j(d)$ and $
y_i^jd=\tau_i^j(d)y_i^j$. Indeed, when $n=1$, $x\tau (d) =
(1-E_{00})dx=dx-E_{00}dx=dx$ since $E_{00}d\in KE_{00}$ and
$E_{00}x=0$.

$\noindent $

Next, we consider two involutions of the algebra $\mA_n$. Using
them we construct  the subgroup $\Xi_n$ of $\mG_n$ which is one of
the building blocks of the group $\mG_n$. We use the fact that a
product of two involutions is an {\em automorphism}.

 {\bf The
involution $\th$ on $\mA_n$}. The {\em involution}
$$\th : A_n\ra A_n, \;\; x_i\mapsto \der_i, \;\; \der_i\mapsto
x_i, \;\; i=1, \ldots , n,$$ of the Weyl algebra $A_n$ can be
uniquely extended to the involution (see \cite{Bav-Jacalg})
\begin{equation}\label{thinv}
\th : \mA_n\ra \mA_n, \;\; x_i\mapsto \der_i, \;\; \der_i\mapsto
x_i,\;\;  H_i^{\pm 1}\mapsto  H_i^{\pm 1},  \;\; i=1, \ldots , n,
\end{equation}
of the Jacobian algebra $\mA_n$.  Moreover, 
\begin{equation}\label{nthEij}
\th (E_{\alpha \beta})= \frac{\alpha !}{\beta !}\, E_{\beta
\alpha},
\end{equation}
where $\alpha !:= \alpha_1!\cdots \alpha_n!$ and $0!:=1$; and so
\begin{equation}\label{thFn}
\th (F_n ) =  F_n.
\end{equation}
Since $ \th (x_i) = H_iy_i$  and $\th (y_i) = x_iH_i^{-1}\not\in
\mS_n$ for $ i=1, \ldots , n$, we see that $\th (\mS_n)
\not\subseteq \mS_n$, i.e. the involution $\th$ does not preserve
the subalgebra $\mS_n$ of $\mA_n$.

$\noindent $

{\bf The involution $\eta$ of the algebra $\mA_n$}. The involution
$\eta$ of the algebra $\mS_n$ can be uniquely extended to the
involution 
\begin{equation}\label{exteta}
\eta : \mA_n\ra \mA_n , \;\; x_i\mapsto y_i, \;\; y_i\mapsto x_i,
\;\;  H_i^{\pm 1}\mapsto H_i^{\pm 1}, \;\; i=1, \ldots , n.
\end{equation}
This can be easily verified using the defining relations of the
algebra $\mA_n$ given by the presentation of the algebra $\mA_n$
as a GWA.  In particular, $\eta (\der_i) = x_iH_i$. Since $\eta
(x_i) = y_i = H_i^{-1}\der_i\not\in A_n$, wee see that $\eta (A_n)
\not\subseteq A_n$, i.e. the involution $\eta$ does not preserve
the Weyl algebra $A_n$. The compositions of the involutions  $\eta
\th$ and $\th \eta$ are {\em automorphisms} of the algebra
$\mA_n$. Moreover, they  belong to the stabilizer
$$ \St_{\mG_n}(H_1, \ldots , H_n):=\{ \s \in \mG_n \, | \, \s
(H_1)=H_1, \ldots , \s (H_n) = H_n\},$$ and \begin{eqnarray*}
 \eta \th : & x_i\mapsto x_iH_i, \;\;\;\; y_i\mapsto H_i^{-1}y_i, \;\; H_i^{\pm 1}\mapsto H_i^{\pm 1}, \\
  \th \eta :  & x_i\mapsto x_iH_i^{-1}, \;\; y_i\mapsto H_iy_i, \;\;\;\; H_i^{\pm 1}\mapsto H_i^{\pm 1}.  \\
\end{eqnarray*}
Note that $(\eta \th )^{-1}= \th \eta $ and $(\eta \th )^m:
x_i\mapsto x_iH_i^m$, $ y_i\mapsto H_i^{-m}y_i$ for all $m\in \Z$,
and so the automorphism $\eta \th$ generates a free cyclic
subgroup of $\mG_n$. These automorphisms yield an idea of
introducing the subgroup $\mU_n$  of $\mG_n$ (see Section
\ref{ACSOA}).

Since $\mA_n=\bigotimes_{i=1}^n \mA_1(i)$, there is a natural
inclusion of groups
$$ \prod_{i=1}^n \Aut_{K-{\rm alg}}(\mA_1(i))\subseteq \mG_n, \;\;
(\s_1, \ldots , \s_n) \mapsto (\bigotimes_{i=1}^n \s_i :
\bigotimes_{i=1}^n a_i\mapsto \bigotimes_{i=1}^n \s_i(a_i)), $$
where $a_i\in \mA_1(i)$.

$\noindent $

{\bf The group $\Xi_n$}.  For each number $i=1, \ldots , n$, let
$\eta_i$ and $\th_i$ be the involutions $\eta$ and $\th$ for the
algebra $\mA_1(i)$. Consider the subgroup of $\mG_n$,
$$ \Xi_n :=\prod_{i=1}^n \Xi_1(i)=  \langle \eta_1\th_1 \rangle \times \cdots \times \langle \eta_n\th_n\rangle \simeq \Z^n,
 \;\; {\rm where}\;\;  \Xi_1(i) := \langle \eta_i\th_i \rangle .$$
Note that, for all $\alpha = (\alpha_1, \ldots , \alpha_n)\in
\Z^n$,
$$\prod_{j=1}^n (\eta_j\th_j)^{\alpha_j}: x_i\mapsto
x_iH_i^{\alpha_i}, \;\; y_i\mapsto H_i^{-\alpha_i}y_i, \;\;
H_i^{\pm 1}\mapsto H_i^{\pm 1}, $$ for all $i=1, \ldots , n$. We
will see that  the group  $\Xi_n $ consists of  outer
automorphisms of the algebra $\mA_n$ (Corollary \ref{c10Apr9}).


\section{Certain subgroups of $\Aut_{K-{\rm alg}}(\mA_n
)$}\label{ACSOA}

Recall that  $\mG_n:=\Aut_{K-{\rm alg}}(\mA_n)$ is the group of
automorphisms of the algebra $\mA_n$. In this section, a useful
description of the group $\mG_n$ is given (Corollary
\ref{d21Mar9}.(1)), an important (rather peculiar) criterion of
equality of two elements of $\mG_n$ (Theorem \ref{21Mar9}) is
found, and several subgroups of $\mG_n$ are introduced that are
building blocks of the group $\mG_n$. These results are important
in finding the group $\mG_n$. Our goal is to prove that
$\mG_n=S_n\ltimes (\mT^n\times \mU_n)\ltimes \ker (\xi )$ (Theorem
\ref{6Apr9}) and $\mG_n = S_n \ltimes (\mT^n \times \Xi_n) \ltimes
\Inn (\mA_n)$ (Theorem \ref{10Apr9}). In this section, the groups
$S_n$, $\mT^n$, $\mU_n$, $\ker (\xi )$, and $\Inn (\mA_n)$ are
introduced and it is proved that the inclusion $S_n\ltimes
(\mT^n\times \mU_n)\ltimes \Inn (\mS_n) \subseteq \mG_n$ holds
(Lemma \ref{x23Mar9}). For $n=1$, the inclusion is the equality
(Theorem \ref{17Mar9}).

{\bf Ideals and the prime ideals of the algebra $\mA_n$}. Recall
some facts about ideals of the algebra $\mA_n$ that are proved in
\cite{Bav-Jacalg} and will be used later in the paper.

$0$ is a prime ideal of $\mA_n$.
$$ \gp_1:=F\t\mA_{n-1}, \gp_2:= \mA_1\t F\t \mA_{n-2} , \ldots ,
\gp_n :=\mA_{n-1}\t F,$$  are precisely the prime ideals  of
height one of $\mA_n$. Let $\Sub_n$ be the set of all subsets of
$\{ 1, \ldots , n\}$.
\begin{itemize}
\item (Corollary 3.5, \cite{Bav-Jacalg}) {\em The map $\Sub_n\ra
\Spec (\mA_n)$, $ I\mapsto \gp_I:= \sum_{i\in I}\gp_i$, $\emptyset
\mapsto 0$, is a bijection, i.e. any nonzero prime ideal of
$\mA_n$ is a unique sum of primes of height 1; $|\Spec
(\mA_n)|=2^n$; the height of $\gp_I$ is $| I|$;  and}
 \item (Lemma 3.6, \cite{Bav-Jacalg}) $\gp_I\subseteq \gp_J$ {\em iff} $I\subseteq
 J$.
\item (Corollary 3.15, \cite{Bav-Jacalg}) $\ga_n:= \gp_1+\cdots
+\gp_n$ {\em is the only prime ideal of $\mA_n$ which is
completely prime; $\ga_n$ is the only ideal $\ga$ of $\mA_n$ such
that $\ga \neq \mA_n$ and $\mA_n/\ga$ is a Noetherian (resp. left
Noetherian, resp. right Noetherian) ring.} \ \item (Theorem 3.1,
\cite{Bav-Jacalg}) {\em Each ideal $I$ of $\mA_n$ is an idempotent
ideal, i.e. $I^2= I$; and ideals of $\mA_n$ commute} ($IJ= JI$).
\item (Theorem 3.11, \cite{Bav-Jacalg}) {\em The lattice of ideals
of $\mA_n$ is distributive}. \item (Corollary 2.7.(4,7),
\cite{Bav-Jacalg}) {\em The ideal $\ga_n$ is the largest (hence,
the only maximal) ideal of $\mA_n$ distinct from $\mA_n$,  and
$F_n$ is the smallest nonzero ideal of $\mA_n$. } \item (Corollary
2.7.(11), \cite{Bav-Jacalg}) $\GK (\mA_n/ \ga )= 3n$ {\em for all
ideals $\ga$ of $\mA_n$ such that $\ga \neq \mA_n$. }
\end{itemize}

$\noindent $

{\bf The automorphism $\heta\in \Aut (\mG_n)$}. The involution
$\eta$ of the algebra $\mA_n$ yields the automorphism $\heta\in
\Aut (\mG_n)$ of the group $\mG_n$: 
\begin{equation}\label{Aetah}
\heta : \mG_n\ra \mG_n, \;\; \s \mapsto \eta \s \eta^{-1}.
\end{equation}
Clearly, $\heta^2= e$ and $\heta (\s ) = \eta \s \eta$ since
$\eta^2= e$.

$\noindent $

{\bf The group homomorphism $\xi$}.

\begin{lemma}\label{e21Mar9}
$\s (\ga_n) = \ga_n$ for all $\s \in \mG_n$.
\end{lemma}
{\it Proof}. The statement is obvious since the ideal $\ga_n$ is
the only maximal ideal of the algebra $\mA_n$. $\Box $

$\noindent $

 By Lemma \ref{e21Mar9}, we have the group
homomorphism (recall that $\CA_n= \mA_n/ \ga_n$):
\begin{equation}\label{Axiaut}
\xi : \mG_n \ra \Aut_{K-{\rm alg}}(\CA_n), \;\; \s\mapsto
(\overline{\s} : a+\ga_n \mapsto \s (a) +\ga_n).
\end{equation}
The homomorphisms $\heta $ and $\xi$ will be used often in the
study of the group $\mG_n$. We will find the group $\Aut_{K-{\rm
alg}}(\CA_n)$ of algebra automorphisms of the skew Laurent
polynomial algebra $\CA_n$ (Theorem \ref{5Apr9}). We are
interested in finding the image and the kernel of the homomorphism
$\xi$ (Theorem \ref{6Apr9}.(2) and Corollary \ref{b18Apr9}). This
will help us to find the group $\mG_n$.  We will see that the
image of $\xi$ is relatively small and the kernel of $\xi$  is
large.

$\noindent $

{\bf The $\mA_n$-module $P_n$}. By the very definition of the
algebra $\mA_n$ as a subalgebra of $\End_K(P_n)$, the
$\mA_n$-module $P_n$ is faithful.

\begin{proposition}\label{J15Ma7}
 {\rm (Corollary 2.7, \cite{Bav-Jacalg})} The polynomial algebra $P_n$ is
the only (up to isomorphism)  faithful, simple $\mA_n$-module.
\end{proposition}
 The $\mA_n$-module $P_n$ is a very special module for the algebra
 $\mA_n$. Its properties, especially the uniqueness, are used
 often in this paper. The polynomial algebra $P_n=
 \bigoplus_{\alpha \in \N^n}Kx^\alpha$ is a naturally $\N^n$-graded
 algebra. This grading is compatible with the $\Z^n$-grading of
 the algebra $\mA_n$, i.e. the polynomial algebra $P_n$ is a $\Z^n$-graded $\mA_n$-module.
  Each element $y_i\in \mA_n\subseteq
 \End_K(P_n)$ is a {\em locally nilpotent} map, that is
 $P_n=\bigcup_{j\geq 1}\ker_{P_n}(y_i^j)$. Moreover,
 $$ \bigcap_{i=1}^n \ker_{P_n}(y_i) = K.$$
Each element $x_i\in \mA_n\subseteq
 \End_K(P_n)$ is an injective (but not a surjective) map. Each
 element $H_i\in \mA_n\subseteq
 \End_K(P_n)$ is a {\em semi-simple} map (that is $P_n=
 \bigoplus_{\l \in K} \ker_{P_n} (H_i-\l )$) with the set of
 eigenvalues $\Z_+:=\{ 1,2, \ldots \}$ since $H_i*x^\alpha =
 (\alpha_i+1) x^\alpha$ for all $\alpha \in \N^n$. Moreover,
\begin{equation}\label{kerHia}
\bigcap_{i=1}^n \ker_{P_n}(H_i-(\alpha_i+1))=Kx^\alpha,
 \;\;
 \alpha \in \N^n.
\end{equation}
In particular, the $K[H_1, \ldots , H_n]$-module
$P_n=\bigoplus_{\alpha \in \N^n}Kx^\alpha$ is the sum of simple,
non-isomorphic,  one-dimensional submodules $Kx^\alpha$, and so
$P_n$ is a semi-simple $K[H_1, \ldots , H_n]$-module.

For the Weyl algebra $A_n$, the $A_n$-module $A_n/\sum_{i=1}^n
A_n\der_i$ is isomorphic to $P_n$ via $1+\sum_{i=1}^n
A_n\der_i\mapsto 1$. The same statement is true for the algebra
$\mA_n$.

\begin{proposition}\label{x16Apr9}
The $\mA_n$-module $\mA_n/\sum_{i=1}^n \mA_n\der_i=
\mA_n/\sum_{i=1}^n \mA_ny_i$ is isomorphic to $P_n$ via
$1+\sum_{i=1}^n \mA_ny_i\mapsto 1$.
\end{proposition}

{\it Proof}. Since $ \mA_n/\sum_{i=1}^n \mA_ny_i\simeq
\bigotimes_{i=1}^n \mA_1(i)/ \mA_1(i) y_i$, it suffices to prove
the statement for $n=1$. In this case, there is the $\mA_1$-module
epimorphism $f: \mA_1/ \mA_1y\ra P_1$, $1+\mA_1y\mapsto 1$. Now,
the statement follows from Lemma \ref{a14Apr9}.(2) since
$f(E_{i0})=x^i$ for all $i\geq 0$, and so $\ker (f) = \mA_1y$.
$\Box $

\begin{lemma}\label{a14Apr9}
\begin{enumerate}
\item $\mA_1= x\mA_1\bigoplus \bigoplus_{i\geq 0}KE_{0i}$. \item
$\mA_1= \mA_1y\bigoplus \bigoplus_{i\geq 0}KE_{i0}$.
\end{enumerate}
\end{lemma}

{\it Proof}. 1. Recall that $\mA_1=\bigoplus_{i\geq 1}
\mD_1y^i\bigoplus \bigoplus_{i\geq 0} x^i\mD_1$ and $\mD_1= \s
(\mD_1)  \bigoplus KE_{00}$. For each $i\geq 1$, $x\mD_1y^i= \s
(\mD_1) xy^i= \s (\mD_1) xy y^{i-1} = \s (\mD_1) \s (1) y^{i-1}
=\s (\mD_1 1) y^{i-1} =\s (\mD_1 ) y^{i-1}$, and so $\mD_1y^{i-1}=
\s (\mD_1)y^{i-1}\bigoplus KE_{00}y^{i-1} = x\mD_1y^i\bigoplus
KE_{0,i-1}$. Therefore, $\mA_1=x\mA_1\bigoplus \bigoplus_{i\geq 0}
KE_{0i}$.

2. Statement 2 is obtained from statement 1 by applying the
involution $\eta$. $\Box $

$\noindent $

Let $A$ be an algebra and $\s$ be its automorphism. For an
$A$-module $M$, the {\em twisted} $A$-module  ${}^{\s}M$, as a
vector space,  coincides with the module $M$ but the action of the
algebra $A$ is given by the rule: $a\cdot m :=\s (a)m$ where $a\in
A$ and $m\in M$.

\begin{corollary}\label{y16Apr9}
\begin{enumerate}
\item Let $M$ be an $\mA_n$-module. Then $\Hom_{\mA_n}(P_n,
M)\simeq \bigcap_{i=1}^n \ker (y_{i,M})$, $f\mapsto f(1)$,  where
$y_{i,M}:M\ra M$, $m\mapsto y_i m$. In particular,
$\End_{\mA_n}(P_n) \simeq K$.  \item By Proposition \ref{J15Ma7},
for each automorphism $\s \in \mG_n$, the $\mA_n$-modules $P_n$
and ${}^\s P_n$ are isomorphic, and each isomorphism $f:P_n\ra
{}^\s P_n$ is
 given by the rule: $f(p) = \s (p) *v$, where  $v=f(1)$ is any
 nonzero element
 of the 1-dimensional vector space $\bigcap_{i=1}^n \ker (\s
 (y_i)_{P_n})$.
\end{enumerate}
\end{corollary}

As an application of these results to the $\mA_n$-module $P_n$, we
have a useful criterion of equality of two elements in the group
$\mG_n$. This criterion is used in the proof of the fact that the
kernel $\ker (\xi )$ has trivial centre (Theorem \ref{A16Apr9}).
Another (even more unexpected) criterion is given in Theorem
\ref{21Mar9} which is used in many proofs in this paper.

\begin{corollary}\label{yz16Apr9}
Let $\s , \tau \in \mG_n$. Then $\s = \tau$ iff $\s (E_{\alpha
0})=\tau (E_{\alpha 0})$ for all $\alpha \in \N^n$ iff $\s
(E_{0\alpha })=\tau (E_{0\alpha })$ for all $\alpha \in \N^n$ iff
$\s (E_{\alpha \beta })=\tau (E_{\alpha \beta})$ for all $\alpha ,
\beta \in \N^n$.
\end{corollary}

{\it Proof}.  The last `iff' follows from the previous two. The
second `iff' follows from the first one by using the automorphism
$\heta$ of the group $\mG_n$: $\s =\tau$ iff $\heta (\s ) = \heta
(\tau )$ iff $\heta (\s ) (E_{\alpha 0})= \heta (\tau ) (E_{\alpha
0})$ for all $\alpha \in \N^n$ (by the first `iff') iff  $\eta \s
 (E_{0\alpha })= \eta \tau  (E_{0\alpha })$ for all $\alpha \in
\N^n$ (since $\eta (E_{\alpha 0}) = E_{0\alpha}$) iff $\s
(E_{0\alpha })= \eta (E_{0\alpha })$ for all $\alpha \in \N^n$ (by
applying $\eta^{-1}$ to the previous equality).

So, it remains to prove that if $\s (E_{\alpha 0}) = \tau
(E_{\alpha 0})$ for all $\alpha \in \N^n$  then $\s = \tau$.
Without loss of generality we may assume that $\tau =e$, the
identity of the group $\mG_n$. So, we have to prove that if $\s
(E_{\alpha 0})= E_{\alpha 0}$ for all $\alpha \in \N^n$ then $\s
=e$. For each number $i=1, \ldots , n$,
$$ 0=(1-E_{00}(i))*1=\s (1-E_{00}(i))*1=\s (x_iy_i) *1= \s (x_i)
\s(y_i)*1, $$ and so $0=\s (y_i) \s (x_i) \s (y_i) *1= \s (y_ix_i)
\s (y_i) *1=\s (y_i) *1$, i.e. $\bigcap_{i=1}^n \ker (\s
(y_i)_{P_n})=K$. By Corollary \ref{y16Apr9}.(2), the map $f:P_n\ra
{}^\s P_n$, $p\mapsto \s (p) *1$, is an $\mA_n$-module
isomorphism. Now, $f(x^\alpha ) = f(E_{\alpha 0}*1)= \s (E_{\alpha
0})*1=E_{\alpha 0}*1=x^\alpha$ for all $\alpha \in \N^n$. This
means that $f$ is the identity map. For all $a\in \mA_n$ and $p\in
P_n$, $a*p= f(a*p)= \s (a)*f(p)= \s (a) *p$, and so $\s (a) = a$
since the $\mA_n$-module $P_n$ is faithful. That is $\s = e$, as
required.  $\Box $

$\noindent $

{\bf A description of the group $\mG_n$}.  The next lemma is
extremely useful in finding the group of automorphisms of algebras
that have a {\em unique faithful} module that satisfies an
isomorphism-invariant property.

\begin{lemma}\label{c21Mar9}
Suppose that an algebra $A$ has a unique (up to isomorphism)
faithful $A$-module $M$ that satisfies an isomorphism-invariant
property, say $\CP$. Then
$$ \Aut_{K-{\rm alg}}(A) = \{ \s_\v \, | \, \v \in \Aut_K(M), \;
\v A\v^{-1} = A\}$$ where  $\s_\v (a) := \v a \v^{-1}$ for $a\in
A$,  and the algebra $A$ is identified with its isomorphic copy in
$\End_K(M)$ via the algebra monomorphism $a\mapsto (m\mapsto am)$.
\end{lemma}

{\it Proof}.  Let $\s \in  \Aut_{K-{\rm alg}}(A)$. The twisted
$\mA_n$-module  ${}^{\s}M$ is faithful and satisfies the property
$\CP$. By the uniqueness of $M$, the $A$-modules $M$ and
${}^{\s}M$ are isomorphic. So, there exists an element $\v \in
\Aut_K(M)$ such that $\v a = \s (a)\v$ for all $a\in A$, and so
$\s (a) = \v a \v^{-1}$, as required. $\Box $

$\noindent $

{\it Example}. The matrix algebra $M_d(K)$ has a unique (up to
isomorphism)  simple module  which is $K^n$. Then, by Lemma
\ref{c21Mar9}, every automorphism of $M_d(K)$ is inner.

$\noindent $

Recall that the polynomial algebra $P_n$ is a unique (up to
isomorphism) faithful,  simple module for the algebra $\mA_n$ and
the algebra $\mS_n$ (see Corollary 2.7,  \cite{Bav-Jacalg}; and
Corollary 3.3, \cite{shrekalg},  respectively).
\begin{corollary}\label{d21Mar9}
\begin{enumerate}
\item  $\mG_n = \{ \s_\v \, | \, \v \in \Aut_K(P_n), \; \v \mA_n
\v^{-1} = \mA_n\}$ where $\s_\v (a) := \v a \v^{-1}$, $a\in
\mA_n$. \item (Theorem 3.2, \cite{shrekaut}) $\mS_n = \{ \s_\v \,
| \, \v \in \Aut_K(P_n), \; \v \mS_n \v^{-1} = \mS_n\}$ where
$\s_\v (a) := \v a \v^{-1}$, $a\in \mS_n$.
\end{enumerate}
\end{corollary}

In \cite{shrekaut}, Corollary \ref{d21Mar9}.(2) was used in
finding the group $G_n$. Next, several important subgroups of
$\mG_n$ are introduced, they are building blocks of the group
$\mG_n$ (Theorem \ref{6Apr9}, Theorem \ref{10Apr9}).

$\noindent $

{\bf The group $\Inn (\mA_n)$ of inner automorphism of $\mA_n$}.
Let $\mA_n^*$ be the group of units of the algebra $\mA_n$. The
centre $Z(\mA_n)$ of the algebra $\mA_n$ is $K$,
\cite{Bav-Jacalg}. For each element $u\in \mA_n^*$, let
$\o_u:\mA_n\ra \mA_n$, $ a\mapsto uau^{-1}$, be the inner
automorphism associated with the element $u$. Then the group of
inner automorphisms of the algebra $\mA_n$, $\Inn (\mA_n)= \{
\o_u\, | \, u\in \mA_n^*\}\simeq \mA_n^*/ K^*$,  is a normal
subgroup of $\mG_n$. It follows from the equality
\begin{equation}\label{AetaH1}
\heta (\o_u) = \o_{\eta (u)^{-1}}, \;\; u\in \mA_n^*,
\end{equation}
that $\heta (\Inn (\mA_n)) = \Inn (\mA_n)$. Clearly, $\xi (\Inn
(\mA_n))\subseteq \Inn (\CA_n)$. The group $\mA_1^*$ was found in
\cite{Bav-Jacalg}, Theorem 4.2.

$\noindent $

{\bf The  algebraic  torus $\mT^n$}. The $n$-dimensional
algebraic   torus $\mT^n:= \{ t_\l \, | \, \l = (\l_1, \ldots ,
\l_n) \in K^{*n}\}$ is a subgroup of the groups $\mG_n$ and  $G_n$
where
$$t_\l (x_i) = \l_ix_i, \;\; t_\l (y_i)= \l_i^{-1}y_i, \;\; t_\l(H_i^{\pm 1})= H_i^{\pm 1}, \;\; i=1,
\ldots , n.$$ Moreover, $\mT^n$ is a subgroup of the group
$\Aut_{K-{\rm alg}}(A_n)$ since $t_\l (\der_i) = t_\l (H_iy_i) =
\l_i^{-1} H_iy_i= \l_i^{-1}\der_i$.  The  algebraic  torus $\mT^n$
is also a subgroup of the groups $\Aut_{K-{\rm alg}}(\CA_n)$ and
$\Aut_{K-{\rm alg}}(L_n)$ where
$$t_\l
(x_i) = \l_ix_i, \;\;  t_\l(H_i)= H_i, \;\; i=1, \ldots , n.$$
Then $\heta (\mT^n) = \mT^n$ and $\heta (t_\l ) = t_\l^{-1} =
t_{\l^{-1}}$ where $\l^{-1} := (\l_1^{-1}, \ldots , \l_n^{-1})$;
$\xi (\mT^n) = \mT^n$ and $\xi (t_\l ) = t_\l$. So, the maps
$\heta : \mT^n\ra \mT^n$ and $\xi : \mT^n\ra \mT^n$ are group
isomorphisms. Note that 
\begin{equation}\label{tlEab}
t_\l (E_{\alpha \beta})=\l^{\alpha - \beta}E_{\alpha, \beta }
\end{equation}
where $\l^{\alpha - \beta}:= \prod_{i=1}^n\l_i^{\alpha_i -
\beta_i}$.

$\noindent $

{\bf The symmetric group $S_n$}. The groups $\mG_n$ and  $G_n$
contain the symmetric group $S_n$ where each elements $\tau $ of
$S_n$ is identified with the automorphism of the algebras $\mA_n$
and  $\mS_n$ given by the rule:
$$ \tau (x_i) = x_{\tau (i)}, \;\; \tau (y_i) = y_{\tau (i)}, \;\;
\tau (H^{\pm 1}_i) = H^{\pm 1}_{\tau (i)},\;\; i=1, \ldots , n.$$
The group $S_n$ is also a subgroup of the groups $\Aut_{K-{\rm
alg}}(\CA_n)$ and $\Aut_{K-{\rm alg}}(L_n)$ where
$$ \tau (x_i) = x_{\tau (i)},\;\; \tau (H_i) = H_{\tau (i)},\;\;
i=1, \ldots , n.$$ Clearly, $\heta (S_n) = S_n$ and $\heta (\tau )
= \tau$ for all $\tau \in S_n$;  $\xi (S_n) = S_n$ and $\xi (\tau
) = \tau$ for all $\tau \in S_n$. Note that 
\begin{equation}\label{tlEab1}
\tau (E_{\alpha \beta })= E_{\tau (\alpha ) \tau (\beta )}
\end{equation}
where $\tau (\alpha ):= (\alpha_{\tau^{-1}(1)}, \ldots ,
\alpha_{\tau^{-1}(n)})$.

$\noindent $

The groups $\mG_n$, $G_n$, $\Aut_{K-{\rm alg}}(\CA_n)$,  and
$\Aut_{K-{\rm alg}}(L_n)$ contain the semi-direct product
$S_n\ltimes \mT^n$ since $\mT^n\cap S_n = \{ e\}$ and
\begin{equation}\label{tt1}
 \tau t_\l \tau^{-1}= t_{\tau (\l )}\;\; {\rm where}\;\; \tau (\l
) := (\l_{\tau^{-1}(1)},  \ldots , \l_{\tau^{-1}(n)}),
\end{equation}
 for
all $\tau \in S_n$ and $t_\l \in \mT^n$. Clearly, the maps
\begin{eqnarray*}
\heta : S_n\ltimes \mT^n\ra S_n\ltimes \mT^n, &  \tau t_\l\mapsto \tau t_\l^{-1},  \\
 \xi : S_n\ltimes \mT^n\ra S_n\ltimes \mT^n, &  \tau t_\l\mapsto \tau t_\l,
\end{eqnarray*}
are group isomorphisms.

By (Theorem 5.1, \cite{shrekaut}), $G_n=S_n\ltimes \mT^n\ltimes
\Inn (\mS_n)$. The algebras $\mS_n$ and  $\mA_n$ are central, and
so $\Inn (\mS_n) \simeq \mS_n^*/K^*\subseteq \mA_n^*/K^*\simeq
\Inn (\mA_n)$. Now, statement 1 of the next proposition follows.

\begin{proposition}\label{f21Mar9}
\begin{enumerate}
\item $G_n = S_n\ltimes \mT^n \ltimes \Inn (\mS_n)\subseteq \mG_n
$.\item $G_n= \{ \s \in \mG_n\, | \, \s (\mS_n ) = \mS_n\}$.
\end{enumerate}

\end{proposition}

{\it Proof}. 2. Statement 2 follows from statement 1 and Corollary
\ref{d21Mar9}. $\Box$.

$\noindent $

Proposition \ref{f21Mar9} means that every automorphism of the
algebra $\mS_n$ can be extended to an automorphism of the algebra
$\mA_n$. Moreover, it can be extended uniquely (Corollary
\ref{g21Mar9}).

$\noindent $

{\bf The subgroup $\mU_n$ of $\mG_n$}. The group $\mA_1^*$ of
units of the algebra $\mA_1$ contains the following infinite
discrete subgroup Theorem 4.2, \cite{Bav-Jacalg}:
\begin{equation}\label{defH1}
\CH := \{ \prod_{i\geq 0} (H+i)^{n_i}\cdot \prod_{i\geq
1}(H-i)^{n_{-i}}_1\, | \, (n_i)\in \Z^{(\Z )}\}\simeq \Z^{(\Z )}.
\end{equation}
 For each
tensor multiple $\mA_1(i)$ of the algebra
 $\mA_n=\bigotimes_{i=1}^n\mA_1(i)$, let $\CH_1 (i)$ be the corresponding
group $\CH$. Their (direct) product 
\begin{equation}\label{defH2}
 \CH_n:= \CH_1 (1) \cdots \CH_1
(n)= \prod_{i=1}^n\CH_1 (i)
\end{equation}
 is a (discrete) subgroup of the group
$\mA_n^*$  of units of the algebra $\mA_n$, and $\CH_n\simeq
\CH^n\simeq (\Z^n)^{(\Z )}$.

For each element $u=u_1\ldots u_n=(u_1, \ldots , u_n)\in \CH_n =
\prod_{i=1}^n\CH_1 (i)$, using the presentation of the algebra
$\mA_n$ as a GWA and the fact that $\mD_n$ is a {\em commutative}
algebra we can easily check that the map $\mu_u : \mA_n\ra \mA_n$
is an algebra automorphism where 
\begin{equation}\label{muu}
\mu_u: x_i\mapsto x_iu_i, \;\; y_i\mapsto u_i^{-1}y_i, \;\;
H_i^{\pm 1}\mapsto H_i^{\pm 1}, \;\; i=1, \ldots , n.
\end{equation}
Let us give two typical verifications that the map $\mu_u$ is an
automorphism of the algebra $\mA_n$ which use the fact that the
algebra $\mA_n$ is a GWA:
$$ \mu_u(x_k) \mu_u(y_k) = x_ku_ku_k^{-1}y_k= x_ky_k = 1-E_{00}(k)
= \mu_u (1-E_{00}(k)), $$ and, for each element $d\in \mD_n$,
$\mu_u(d) \mu_u(x_k) = dx_ku_k= x_k\tau_k(d) u_k = x_ku_k\tau_k(d)
= \mu_u(x_k) \mu_u(\tau_k(d))$, by (\ref{muDn=Dn}). Note that
$\mu_u\mu_v= \mu_{uv}$ and $\mu_u^{-1} = \mu_{u^{-1}}$, and so we
have the subgroup of
$\mG_n$, 
\begin{equation}\label{UnZZ}
\mU_n:= \{ \mu_u\, | \, u=(u_1, \ldots , u_n) \in \CH_n =
\prod_{i=1}^n \CH_1 (i)\} \simeq (\Z^n)^{(\Z)},\;\; \mu_u\lra
(n_i(k)),
\end{equation}
where $u_k=\prod_{i\geq 0} (H_k+i)^{n_i(k)}\cdot \prod_{i\geq 1}
(H_k-i)_1^{n_{-i}(k)}$. A direct calculation shows that
\begin{equation}\label{muEij}
\mu_u(E_{ij}(k))= \begin{cases}
E_{ij}(k)u_k\tau_k(u_k) \cdots \tau_k^{i-j-1}(u_k)& \text{if }i>j,\\
E_{ii}(k)& \text{if }i=j,\\
u_k^{-1}\tau_k(u_k^{-1}) \cdots \tau_k^{j-i-1}(u_k^{-1})E_{ij}(k)& \text{if }i<j.\\
\end{cases}
\end{equation}
It follows that 
\begin{equation}\label{muDn=Dn}
\mu_u (\mD_n) = \mD_n, \;\; \mu_u|_{\mD_n} = \id_{\mD_n}.
\end{equation}
By the very definition of the group $\mU_n$, it is  the direct
product of its subgroups $\mU_1(k)$:
$$ \mU_n = \mU_1(1) \times \cdots \times \mU_1(n)$$
where $\mU_1(k) := \{\mu_{u_k}\, | \, u_k\in \CH_1 (k)\}\simeq
\Z^{(\Z )}$.

Note that $\mT^n\cap \mU_n = \{ e\}$ and $t_\l \mu_u = \mu_ut_\l$
for all $t_\l \in \mT^n$ and $\mu_u\in \mU_n$. Therefore,
$\mG_n\supset \mT^n\mU_n= \mT^n\times \mU_n$. Note that $S_n\cap
\mT^n\mU_n = \{ e\}$ and, 
\begin{equation}\label{smus}
s\mu_{(u_1, \ldots , u_n)}s^{-1} = \mu_{(s(u_{s^{-1}(1)}),\ldots ,
s(u_{s^{-1}(n )}))},
\end{equation}
for all $s\in S_n$ and $\mu \in \mU_n$. Therefore, $ S_n\ltimes
(\mT^n\times \mU_n)\subseteq \mG_n$.  The restriction of the
homomorphism $\xi$ (see (\ref{Axiaut})) to the subgroup $\mU_n$ of
$\mG_n$ yields the group isomorphism $\xi : \mU_n\simeq \xi
(\mU_n)=\{ \mu_{\xi (u)}:x_k\mapsto x_k\xi (u_k), H_k\mapsto H_k\,
| \, u\in \CH_n \}$. Moreover, $\xi :S_n\ltimes (\mT^n\times
\mU_n)\simeq S_n \ltimes (\mT^n \times \xi (\mU_n))$.

\begin{lemma}\label{x23Mar9}
$S_n\ltimes (\mT^n\times \mU_n)\ltimes \Inn (\mS_n) \subseteq
\mG_n$.
\end{lemma}

{\it Remark}. In fact, the equality holds for $n=1$ (Theorem
\ref{17Mar9}.(1)), and the strict inclusion holds when $n>1$. The
last statement can be easily deduced from Theorem \ref{10Apr9} and
the fact that the set $\mA_n^*/K^*$ is much more massive than
$\mS_n^*/K^*$. Note that $\Inn (\mA_n) \simeq \mA_n/K^*\supset
\mS_n^*/K^* \simeq \Inn (\mS_n)$.

{\it Proof}. The groups $S_n\ltimes (\mT^n\times \mU_n)$ and $
\Inn (\mS_n)$ are subgroups of the group $\mG_n$ (Proposition
\ref{f21Mar9}). Since $\Inn (\mS_n) \subseteq \ker (\xi )$ and
$\xi :S_n\ltimes (\mT^n\times \mU_n)\simeq S_n \ltimes (\mT^n
\times \xi (\mU_n))$, we see that $S_n\ltimes (\mT^n\times
\mU_n)\cap \Inn (\mS_n) = \{ e\}$, and the statement follows.
$\Box $

$\noindent $

{\bf The subgroup $\mU^0_n$ of $\mG_n$}. For each number $i=1,
\ldots , n$, let $\CH_1^0(i)$ be the kernel of the group
epimorphism 
\begin{equation}\label{degHi}
\deg_{H_i} : \CH_1(i) \ra \Z, \;\; \prod_{j\geq
0}(H_i+j)^{n_j}\cdot \prod_{j\geq 1}(H_i-j)_1^{n_{-j}}\mapsto
\sum_{k\in \Z}n_k.
\end{equation}
Then $\CH_1(i) = \{ H_i^j\, | \, j\in \Z\} \times \CH_1^0(i)$, and
so $\CH_n = \{ H^\alpha \, | \, \alpha \in \Z^n\} \times \CH_n^0$
where $\CH_n^0:= \prod_{i=1}^n \CH_1^0(i)$. Correspondingly,
\begin{equation}\label{degHi1}
\mU_n= \{ \mu_{H^\alpha}\, | \, \alpha \in \Z^n\} \times
\mU_n^0=\Xi_n  \times \mU^0_n \;\; {\rm where}\;\; \mU_n^0:=\{
\mu_u\, | \, u\in \CH_n^0\} \simeq \CH_n^0, \;\; \mu_u\lra u.
\end{equation}
Let $r$ be an element of a ring $R$. The
 element $r$ is called {\em regular} if $ \lann_R(r)=0$ and
 $\rann_r(r)=0$ where $ \lann_R(r):= \{ s\in R\, | \, sr=0\}$ is
 the {\em left annihilator} of $r$ and $ \rann_R(r):= \{ s\in R\, | \, rs=0\}$ is
 the {\em right annihilator} of $r$ in $R$.

The next lemma shows that the elements $x$ and $y$ of the algebra
$\mS_1$ are not regular. Note that the element $x\in A_1$ is  a
regular element of the Weyl algebra $A_1$ which is not regular as
an element of the algebra $\mA_1$.
\begin{lemma}\cite{shrekalg}\label{a7Dec8}
\begin{enumerate}
\item $\lann_{\mS_1}(x) = \mS_1E_{00} = \bigoplus_{i\geq 0}
KE_{i,0}=\bigoplus_{i\geq 0} Kx^i(1-xy)$ and $\rann_{\mS_1}(x)=0$.
\item $\rann_{\mS_1}(y) = E_{00}\mS_1 =  \bigoplus_{i\geq 0}
KE_{0, i}=\bigoplus_{i\geq 0} K (1-xy)y^i$ and
$\lann_{\mS_1}(y)=0$.
\end{enumerate}
\end{lemma}

It follows from Lemma \ref{a7Dec8} and the presentation of the
algebra $\mA_n$ as a GWA (see also (\ref{CD1d})) that, for each
$i=1, \ldots , n$, 
\begin{equation}\label{Alanxin}
\lann_{\mA_n}(x_i) =
 \bigoplus_{j\geq 0} \mA_{n-1, i}E_{j,0}(i)=\bigoplus_{j\geq 0}
 \mA_{n-1, i}x_i^jE_{00}(i),
\end{equation}
\begin{equation}\label{Alanyin}
\rann_{\mA_n}(y_i) =
 \bigoplus_{j\geq 0} E_{0,j}(i)\mA_{n-1, i}=\bigoplus_{j\geq
 0}E_{00}(i)y_i^j\mA_{n-1, i},
\end{equation}
where $\mA_{n-1, i}$ stands for $\bigotimes_{k\neq i}\mA_1(k)$.

For an algebra $A$ and a subset $S\subseteq A$, $\Cen_A(S):=\{
a\in A\, | \, as=sa$ for all $s\in S\}$ is the {\em centralizer}
of the set $S$ in $A$. It is a subalgebra of $A$. It follows from
the presentation of the algebra $\mA_n$ as a GWA that
\begin{equation}\label{ACenxy}
\Cen_{\mA_n}(x_1, \ldots , x_n)=K[x_1, \ldots , x_n], \;\;
\Cen_{\mA_n}(y_1, \ldots , y_n)=K[y_1, \ldots , y_n].
\end{equation}

Let $\mE_n:= \End_{K-{\rm alg}}(\mA_n)$ be the monoid of all the
$K$-algebra endomorphisms of $\mA_n$. The group of units of this
monoid is $\mG_n$. The automorphism $\heta \in \Aut (\mG_n)$ can
be extended to an automorphism $\heta \in \Aut (\mE_n)$ of the
monoid $\mE_n$: 
\begin{equation}\label{A1etah}
\heta : \mE_n\ra \mE_n, \;\; \s \mapsto \eta \s \eta^{-1}.
\end{equation}
The next  result is instrumental in finding the group of
automorphisms of the algebra $\mA_n$.
\begin{theorem}\label{21Mar9}
Let $\s , \tau \in \mG_n$. Then the following statements are
equivalent.
\begin{enumerate}
 \item $\s = \tau$. \item $\s (x_1) = \tau (x_1), \ldots , \s (x_n) = \tau
 (x_n)$.
 \item $\s (y_1) = \tau (y_1), \ldots , \s (y_n) = \tau (y_n)$.
\end{enumerate}
\end{theorem}

{\it Proof}. Without loss of generality we may assume that $\tau
=e$, the identity automorphism. Consider the following two
subgroups of $\mG_n$, the stabilizers of the sets $\{ x_1, \ldots
, x_n\}$ and $\{ y_1, \ldots , y_n\}$:
\begin{eqnarray*}
\St_{\mG_n}  (x_1, \ldots , x_n):= & \{ g\in \mG_n\, | \, g(x_1) =
x_1,
\ldots , g(x_n)=x_n\}, \\
 \St_{\mG_n}    (y_1, \ldots , y_n):= & \{ g\in \mG_n\,
| \, g(y_1) = y_1, \ldots , g(y_n)=y_n\}.
\end{eqnarray*}
Then $$\heta (\St_{\mG_n}    (x_1, \ldots , x_n))=\St_{\mG_n}
(y_1, \ldots , y_n), \;\; \heta (\St_{\mG_n}    (y_1, \ldots ,
y_n))= \St_{\mG_n} (x_1, \ldots , x_n).$$ Therefore, the theorem
(where $\tau = e$) is equivalent to the single statement that
$\St_{\mG_n} (x_1, \ldots , x_n)=\{ e\}$. So, let $\s \in
\St_{\mG_n} (x_1, \ldots , x_n)$. We have to show that $\s = e$,
i.e. $\s (y_i) = y_i$ and $\s (H_i) = H_i$ for all $i$. For each
$i=1, \ldots , n$, $1=\s (y_ix_i) = \s (y_i) x_i$ and $1=y_ix_i$.
By taking the difference of these equalities we see that $a_i:= \s
(y_i) -y_i\in \lann_{\mA_n}(x_i)$. By (\ref{Alanxin}),
$a_i=\sum_{j\geq 0} \l_{ij}E_{j0}(i)$ for some elements
$\l_{ij}\in \bigotimes_{k\neq i}\mA_1(k)$, and so $$ \s (y_i) =
y_i+\sum_{j\geq 0}\l_{ij} E_{j0}(i).$$ The element $\s (y_i)$
commutes with the elements $\s (x_k) = x_k$, $k\neq i$, hence all
$\l_{ij}\in K[x_1, \ldots , \hx_i, \ldots , x_n]$, by
(\ref{ACenxy}). Since $E_{j0}(i) = x_i^j E_{00}(i)$, we can write
$$ \s (y_i) =y_i+ p_iE_{00}(i) \;\; {\rm for \; some}\;\; p_i\in
P_n.$$ We have to show that all $p_i=0$. Suppose that this is not
the case. Then $p_i\neq 0$ for some $i$. We seek a contradiction.
Note that $\s^{-1}\in \St_{\mG_n}   (x_1, \ldots , x_n)$, and so
$\s^{-1} (y_i) =y_i+ q_iE_{00}(i)$ for some  $q_i\in P_n$. Recall
that $E_{00}(i) = 1-x_iy_i$. Then $\s^{-1} (E_{00}(i))=
1-x_i(y_i+q_iE_{00}(i))= (1-x_iq_i)E_{00}(i)$, and
$$ y_i=\s^{-1}\s (y_i) = \s^{-1} (y_i+p_iE_{00}(i))=
y_i+(q_i+p_i(1-x_iq_i))E_{00}(i), $$ and so $q_i+p_i= x_ip_iq_i$
since the map $P_n\ra P_nE_{00}(i)$, $p\mapsto pE_{00}(i)$, is an
isomorphism of $P_n$-modules as it follows from (\ref{xyEij}).
This is impossible by comparing the degrees of the polynomials on
both sides of the equality. Therefore, $\s (y_i) = y_i$ for all
$i$.

By Proposition \ref{J15Ma7}, there is an $\mA_n$-module
isomorphism $\v : P_n\ra {}^\s P_n$, $p\mapsto \s (p)*v$, where
$v:=\v (1) \in \bigcap_{i=1}^n \ker_{{}^\s P_n}(\s (y_i))=
\bigcap_{i=1}^n \ker_{P_n}(y_i) = K1$.  Without loss of generality
we may assume that $v=1$. Then  $1=\v (1) = \v (H_i*1) = \s (H_i)
*1$ for all $i$. For each  $\alpha \in \N^n$ and $i=1, \ldots ,
n$,
\begin{eqnarray*}
 \s (H_i) *x^\alpha &=&\s (H_i) x^\alpha  *1= \s (H_i) \s (x^\alpha ) *1= \s (H_ix^\alpha)*1= \s(x^\alpha (H_i+\alpha_i))*1 \\
 &=& \s (x^\alpha ) ( \s (H_i)+\alpha_i)*1= x^\alpha (\alpha_i+1)*1=
 (\alpha_i+1) x^\alpha.
\end{eqnarray*}
This means that the linear maps $\s (H_i), H_i\in \End_K(P_n)$
coincide.  Therefore, $\s (H_i) = H_i$ for all $i$ since the
$\mA_n$-module $P_n$ is faithful. This proves that $\s =e$. $\Box
$

\begin{corollary}\label{g21Mar9}
$\St_{\mG_n}(\mS_n ):= \{ \s \in \mG_n\, | \, \s (\mS_n) = \mS_n ,
\s|_{\mS_n}= {\rm id}_{\mS_n}\} = \{ e\}$, and so each
automorphism of the algebra $\mS_n$ can be uniquely extended to an
automorphism of the algebra $\mA_n$ (see Proposition
\ref{f21Mar9}.(1)).
\end{corollary}

{\it Proof}. Exactly the same result as Theorem \ref{21Mar9} is
true for the algebra $\mS_n$, i.e. when we replace the group
$\mG_n$ by the group $G_n$ (Theorem 3.7, \cite{shrekaut}). Now,
the statement is obvious.  $\Box $

$\noindent $

Theorem \ref{21Mar9} states that each automorphism of the
non-commutative, finitely generated, non-Noetherian algebra
$\mA_n$ is uniquely determined by its action on its commutative,
finitely generated subalgebra $P_n$. A similar result is true for
the algebra $\mS_n$ (Theorem 3.7, \cite{shrekaut}) and for  the
ring $\CD (P_n)$ of differential operators on the polynomial
algebra $P_n$ over a field of {\em prime} characteristic. The
algebra $\CD (P_n)$ is a non-commutative, {\em not finitely
generated}, non-Noetherian algebra.

\begin{theorem}
{\rm \cite{autlaur}   (Rigidity of the group $\Aut_{K-{\rm
alg}}(\CD (P_n)$))} Let $K$ be a field of prime characteristic,
and   $\s , \tau \in \Aut_{K-{\rm alg}}(\CD (P_n)$. Then $\s =
\tau$ iff  $\s (x_1) = \tau (x_1), \ldots , \s (x_n) = \tau
 (x_n)$.
\end{theorem}
The above theorem doest not hold in characteristic zero and does
not hold in prime characteristic for the ring of differential
operators on a Laurent polynomial algebra \cite{autlaur}.


\section{The group $\Aut_{K-{\rm alg}}(\mA_1)$ }\label{ATG1S}

In this section, the group $\mG_1$ and its explicit generators are
 found (Theorem \ref{17Mar9}) and it is proved that any algebra
endomorphism of the algebra $\mA_1$ is a {\em monomorphism}
(Theorem \ref{24Mar9}) (note that the algebra $\mA_1$ is not a
simple algebra). The case $n=1$ is rather special and much more
simpler than the general case. It is a sort of  degeneration of
the general case. The key idea in finding the group $\mG_1$ of
automorphisms of the algebra $\mA_1$ is to use Theorem
\ref{21Mar9}, some of the properties of the index of linear maps
in the vector space $P_1=K[x]$, and the
 explicit structure  of the group $\Aut_{K-{\rm alg}}(\mA_1/ F)$ (Theorem
\ref{A18Mar9}). We start this section with recalling some of the
results on the index from \cite{shrekaut} and  on the algebra
$\mA_1$ from \cite{Bav-Jacalg}.

$\noindent $

{\bf The index $\ind$ of linear maps and its properties}. Let $\CC
= \CC (K)$ be the family of all $K$-linear maps with finite
dimensional kernel and cokernel (such maps are called the {\em
Fredholm linear maps/operatorts}).  So, $\CC$ is the family of
{\em Fredholm} linear maps/operators. For vector spaces $V$ and
$U$, let $\CC (V,U)$ be the set of all the linear maps from $V$ to
$U$ with finite dimensional kernel and cokernel. So, $\CC
=\bigcup_{V,U}\CC (V,U)$ is the disjoint union.

$\noindent $

{\it Definition}. For a linear map $\v \in \CC$, the integer
$$ \ind (\v ) := \dim \, \ker (\v ) - \dim \, \coker (\v )$$
is called the {\em index} of the map $\v$.

$\noindent $

{\it Example}. Note that $\mA_1\subset \End_K(P_1)$, and, in
particular, $x^i, y^i\in \End_K(P_1)$. One can easily prove that
\begin{equation}\label{indxy}
\ind (x^i)= -i\;\; {\rm and }\;\; \ind (y^i)= i, \;\; i\geq 1.
\end{equation}

Lemma \ref{b8Feb9} shows that $\CC$ is a multiplicative semigroup
with zero element (if the composition of two elements of $C$ is
undefined we set their product to be zero). The next two lemmas
are well known.

\begin{lemma}\label{b8Feb9}
Let $\psi : M\ra N $ and $\v : N\ra L$ be $K$-linear maps. If two
of the following three maps: $\psi$, $\v$,  and $\v \psi$, belong
to the set $\CC$ then so does the third; and in this case, $$ \ind
(\v\psi ) = \ind (\v ) + \ind (\psi ).$$
\end{lemma}
\begin{lemma}\label{a8Feb9}
Let
$$
\xymatrix{0\ar[r] & V_1\ar[r]\ar[d]^{\v_1}  & V_2 \ar[r]\ar[d]^{\v_2} & V_3 \ar[r]\ar[d]^{\v_3 } & 0 \\
0\ar[r] & U_1\ar[r]  & U_2\ar[r] & U_3 \ar[r] & 0 }
$$
be a commutative diagram of $K$-linear maps with exact rows.
Suppose that $\v_1, \v_2, \v_3\in \CC$. Then
$$ \ind (\v_2) = \ind (\v_1)+\ind (\v_3).$$
\end{lemma}

Each nonzero element $u$ of the skew Laurent polynomial algebra
$\CA_1= \CL_1[x,x^{-1}; \s_1 ]$  (where $\s_1(H)=H-1$) is a unique
sum $u=\l_sx^s+\l_{s+1}x^{s+1}+\cdots +\l_dx^d$ where all $\l_i\in
\CL_1$, $\l_d\neq 0$, and $\l_dx^d$ is the {\em leading term} of
the element $u$. The integer $\deg_x(u)=d$ is called the {\em
degree} of the element $u$, $\deg_x(0):=-\infty $. For all $u,v\in
\CA_1$, $\deg_x(uv)= \deg_x(u) + \deg_x(v)$ and $\deg_x(u+v) \leq
\max \{ \deg_x(u), \deg_x(v)\}$. The next lemma explains how to
compute the index of the elements $\mA_1\backslash F$ via the
degree function $\deg_x$ and proves that the index is a
$\mG_1$-invariant concept. Note that $F\cap \CC = \emptyset$.
\begin{lemma}\label{c24Mar9}
\begin{enumerate}
\item $\CC \cap \mA_1= \mA_1\backslash F $ (recall that
$\mA_1\subset \End_K(P_1)$) and, for each element $a\in
\mA_1\backslash F$,
$$ \ind (a) = -\deg_x(\oa )$$
where $\oa = a+F\in \mA_1/F= \CA_1$. \item $\ind (\s (a)) = \ind
(a)$ for all $\s \in \mG_1$ and $a\in \mA_1\backslash F$.
\end{enumerate}
\end{lemma}

{\it Proof}. 1. Since $\CC \cap F = \emptyset$, to prove that the
equality $\CC \cap \mA_1= \mA_1\backslash F $ holds it suffices to
show that $ \ind (a) = -\deg_x( \oa )$  for all $a\in
\mA_1\backslash F$. Let $a\in \mA_1\backslash F$ and $d:=
\deg_x(\oa )$. The element of the algebra $\mA_1$,
$$b:= \begin{cases}
y^da& \text{if }\; d\geq 0,\\
ax^{-d}& \text{if }\; d<0, \\
\end{cases}
$$
does not belong to the ideal $F$ (since $\ob = x^{-d}\oa \neq 0$),
and $\deg_x(\ob )=0$. By Lemma \ref{b8Feb9} and (\ref{indxy}), it
suffices to prove that $\ind (b)=0$ since then
$$ 0=\ind (b) =
d+\ind (a),$$   that is $\ind (a) = -\deg_x(\oa )$. The element
$b$ can be written as a sum $b=\l +\sum_{i\geq 1}\l_iy^i+f$ for
some elements $\l , \l_i\in \CL_1$ with $\l \neq 0$, and $f\in F$.
Fix a natural number $m$ such that the vector space $V:=
\bigoplus_{i=0}^m Kx^i$ is $b$-invariant ($bV\subseteq V$),  the
element $\l$ acts as an isomorphism on the vector space $U:=
P_1/V\simeq \bigoplus_{i>d}Kx^i$, and $fP_1\subseteq V$. Let $b_1$
and $b_2$ be the restrictions of the linear map $b$ to the vector
spaces $V$ and $U$ respectively. Then $\ind (b_1)=0$ since $\dim
(V)<\infty$; and $\ind (b_2)=0$ since $b_2=\l +\sum_{i\geq
1}\l_iy^i$ is a bijection. Applying Lemma \ref{a8Feb9} to the
commutative diagram
$$
\xymatrix{0\ar[r] & V\ar[r]\ar[d]^{b_1}  & P_1 \ar[r]\ar[d]^{b} & U \ar[r]\ar[d]^{b_2 } & 0 \\
0\ar[r] & V\ar[r]  & P_1\ar[r] & U \ar[r] & 0 }
$$
we have the result: $\ind (b) = \ind (b_1)+\ind (b_2)=0$.

 2. By Corollary  \ref{d21Mar9}.(1), $\ind (\s (a))= \ind (\v a \v^{-1}) =
 \ind (a)$ where $\s = \s_\v$. $\Box $

$\noindent $

{\bf The group of automorphisms of the algebra $\CA_1$}. Let
$\CG_1:= \Aut_{K-{\rm alg}}(\CA_1)$. We are going to find the
group of automorphisms of the algebra $\CA_1:= \CL_1[x^{\pm 1};
\s_1]$ (Theorem \ref{A18Mar9}) where $\CL_1:= K[H,
(H+i)^{-1}]_{i\in \Z}$ and $\s_1 (H) = H-1$. Let $\Z^{(\Z )}$ be
the direct sum of $\Z$ copies of the abelian group $\Z$. The group
of
 units  of the algebra $\CL_1$ is
 $$\CL_1^* = \{ \l \prod_{i\in \Z}(H+i)^{n_i}\, | \, \l \in K^*,
 (n_i) \in \Z^{(\Z )}\} \simeq K^*\times \CL_1'\simeq K^*\times
 \Z^{(\Z )}$$
 where $\CL_1':= \{ \prod_{i\in \Z}(H+i)^{n_i}\, | \,
 (n_i) \in \Z^{(\Z )}\} \simeq
 \Z^{(\Z )}$, $\prod_{i\in \Z}(H+i)^{n_i}\lra (n_i)$. The kernel
 $\CL_1^0$ of the group epimorphism $\CL_1'\ra \Z$,
 $ \prod_{i\in \Z}(H+i)^{n_i}\mapsto \sum_{i\in \Z}n_i$, is equal
 to
 $$ \CL_1^0=\{ \prod_{i\in \Z} (\frac{H+i}{H+i+1})^{n_i}\, | \,
 (n_i) \in \Z^{(\Z )}\} \simeq \Z^{(\Z )}, \;\;
  \prod_{i\in \Z} (\frac{H+i}{H+i+1})^{n_i} \lra (n_i).$$
Similarly, since $(\frac{H+i}{H+i+1})/(\frac{H+i+1}{H+i+2})=
\frac{(H+i)(H+i+2)}{(H+i+1)^2}$, the kernel $\CL_1^{00}$ of the
group epimorphism $\CL_1^0\ra \Z$, $\prod_{i\in
\Z}(\frac{H+i}{H+i+1})^{n_i}\mapsto \sum_{i\in \Z}n_i$, is equal
to
$$ \CL_1^{00}=\{ \prod_{i\in \Z} (\frac{(H+i)(H+i+2)}{(H+i+1)^2})^{n_i}
\, | \, (n_i) \in \Z^{(\Z )}\} \simeq \Z^{(\Z )}.$$
 Consider the following automorphisms of the algebra
 $\CA_1$:
\begin{eqnarray*}
 \s_u:& x\mapsto ux, &\;\; H\mapsto H, \;\;\;\;\; (u\in \CL_1^*) \\
 s_i:=\o_{x^i}:&   x\mapsto x, &\;\; H\mapsto H-i, \;\; (i\in \Z)  \\
 \zeta :&  x\mapsto x^{-1},& \;\; H\mapsto - H,
\end{eqnarray*}
and the subgroups they generate in the group $\Aut_{K-{\rm
alg}}(\CA_1)$: $\mL_1^*:= \{ \s_u\, | \, u\in \CL_1^*\}\simeq
\CL_1^*$, $ \s_u \lra u$; $\{ s_i\, | \, i\in \Z\}\simeq \Z$,
$s_i\lra i $; and $\langle \zeta \rangle \simeq \Z_2:= \Z / 2\Z$
since $\zeta^2= e$. We can easily check that 
\begin{equation}\label{Hpm0}
s_i\s_u s_i^{-1} = \s_{s_i(u)}, \;\; \zeta \s_u\zeta^{-1} =
\s_{\zeta \s_1^{-1} (u^{-1})}, \;\; \zeta s_i\zeta^{-1} =s_i^{-1}.
\end{equation}
Note that $\zeta \s_1 \zeta^{-1} = \s_1^{-1}$ in $\Aut_{K-{\rm
alg}}(\CL_1)$ where $\s_1(H)=H-1$. It follows that the subgroup of
the group of automorphisms of $\CA_1$ generated by all the
automorphisms above is, in fact, the semidirect product
$$ \langle \zeta \rangle \ltimes \{ s_i\, | \, i\in \Z\}
\ltimes \{ \s_u\, | \, u\in \CL_1^*\} \simeq \Z_2\ltimes \Z
\ltimes (K^*\times \Z^{(\Z )})$$ since, for $\varepsilon =0,1$;
 $i\in \Z$; and $u\in \CL_1^*$,
$$ \zeta^\varepsilon s_i\s_u : x\mapsto \zeta^\varepsilon s_i(u)
\cdot \zeta^\varepsilon(x), \;\; H\mapsto (-1)^\varepsilon H-i,$$
and so $\zeta^\varepsilon s_i\s_u= e$ iff $\varepsilon =0$, $i=0$,
 and $u=1$. Theorem \ref{A18Mar9} shows that this is the whole
 group of automorphisms of the algebra $\CA_1$. Let us collect
 results that are used in the proof of Theorem \ref{A18Mar9}. It
 follows from the fact that the fixed ring $\CL_1^{\s_1}:=\{ a\in
 \CL_1\, | \, \s_1(a) = a)\}$ is equal to $K$ that
 $$ \Cen_{\CA_1}(x) = K[x,x^{-1}]\;\; {\rm and}\;\; Z(\CA_1) =
 K.$$
 The algebra $\CA_1$ is a Noetherian domain. Let $\Frac (\CA_1)$
 be its skew field of fractions. The algebra $\CA_1$ contains two
 skew polynomial rings, $K[H][x; \s_1]$ and $K[H][x^{-1};
 \s_1^{-1}]$ which are Noetherian domains. Moreover, $\Frac (\CA_1)$ is their common skew field
 of fractions. The usual concepts of degree $\deg_x(\cdot )$ in
 $x$ and $\deg_{x^{-1}}(\cdot )$ in $x^{-1}$ for the algebras $K[H][x; \s_1]$ and $K[H][x^{-1};
 \s_1^{-1}]$ respectively can be extended to  valuations of their
 skew field of fractions,
 $$\deg_x(\cdot ), \deg_{x^{-1}}(\cdot ): \Frac (\CA_1) \ra \Z
 \cup \{ \infty \} ,$$ by the rule
 $ \deg_{x^{\pm 1}}(a^{-1}b) =\deg_{x^{\pm 1}}(b)-\deg_{x^{\pm
 1}}(a)$ where $a, b\in \CA_1$ with $a\neq 0$. Note that
 $$ \St (x) := \{ \s \in \Aut_{K-{\rm alg}}(\CA_1)\, | \, \s (x) =
 x\} = \{ s_i\, | \, i\in \Z\}$$
 since, for each element $\s \in \St (x)$, $[\s (H)-H, x]= \s
 ([H,x]) - [H,x]= x-x=0$, and so $i:=\s (H) -H\in \Cen_{\CA_1}(x)
 = K[x,x^{-1}]$. Then $i\in \Z$ since the element $\s (H)=H+i$ must be invertible
 in $\CA_1$. We extend the degree function $\deg_H(\cdot)$ from
 $K[H]$ to its field of fractions by the rule $\deg_H(ab^{-1}) =
 \deg_H(a) - \deg_H(b)$. For a group $G$, $[G,G]$ denotes its {\em
 commutant}, i.e. the subgroup of $G$ generated by all the {\em
 commutators} $[a,b]:=aba^{-1}b^{-1}$ of the elements $a,b\in G$.
 The centre of a group $G$ is denoted by $Z(G)$.

\begin{lemma}\label{a6Apr9}
\begin{enumerate}
\item The commutant $[A\ltimes B, A\ltimes B]$  of a skew product
$A\ltimes B$ of two groups is equal to $[A,A]\ltimes ([A,B]\cdot
[B,B])$ where $[A,B]$ is the subgroup of $B$ generated by all the
commutators $[a,b]:= aba^{-1}b$ for $a\in A$ and $b\in B$. Hence,
$B\cap [A\ltimes B, A\ltimes B]=[A,B]\cdot [B,B]$ and
$\frac{A\ltimes B}{[A\ltimes B, A\ltimes B]}\simeq
\frac{A}{[A,A]}\times \frac{B}{[A,B]\cdot [B,B]}$. \item
 If, in addition, the group $B$ is a direct product of groups
 $\prod_{i=1}^mB_i$ such that $aB_ia^{-1} \subseteq B_i$ for all
 elements $a\in A$ and $i=1, \ldots , m$ then
 $[A\ltimes B, A\ltimes B]= [A,A]\ltimes \prod_{i=1}^m([A,B_i][B_i, B_i]).$
\end{enumerate}

\end{lemma}

{\it Proof}. 1.  Note that $[a,b]= \o_a(b) b^{-1}$ where $\o_a (b)
:= aba^{-1}$. For $a\in A$ and $b, c\in B$,
\begin{eqnarray*}
 c[a,b]&=&c\o_a(b)b^{-1} = \o_a( \o_{a^{-1}}(c) b) (\o_{a^{-1}}(c)b)^{-1} \o_{a^{-1}}(c) bb^{-1} \\
 &=&\o_a( \o_{a^{-1}}(c) b) (\o_{a^{-1}}(c)b)^{-1}\cdot
 \o_{a^{-1}}(c)\\
 &=& [ a,\o_{a^{-1}}(c) b]\cdot \o_{a^{-1}}(c).
\end{eqnarray*}
It follows from these equalities (when, in addition, we chose
$c\in [B,B]$) that the subgroup of $B$ which is generated by its
two subgroups, $[A,B]$ and $[B,B]$, is equal their set theoretic
product $[A,B][B,B]:=\{ ef\, | \, e\in [A,B], f\in [B,B]\}$. Then
the subgroup of $C:=[A\ltimes B, A\ltimes B]$ which is  generated
by its three subgroups $[A,A]$, $[A,B]$, and $[B,B]$ is equal to
the RHS, say $R$, of the equality of the lemma. It remains to
prove that $C\subseteq R$. This inclusion follows from the fact
that,  for all $a_1,a_2\in A$ and $b_1, b_2\in B$,
\begin{equation}\label{a1b1}
 [a_1b_1, a_2b_2]= \o_{a_1}([b_1, a_2]) \o_{a_1a_2}([b_1, b_2]) [
a_1, a_2] \o_{a_2}([a_1, b_2])
\end{equation}
 which follows from the equalities
$[ab, c] = \o_a([b,c])[a,c]$ and $[a,b]^{-1} = [ b,a]$:
\begin{eqnarray*}
 [a_1b_1, a_2b_2]&=& \o_{a_1}([b_1, a_2b_2]) [ a_1, a_2b_2]=
 ([a_2b_2, a_1]\o_{a_1}([a_2b_2, b_1]))^{-1}\\
 & =& (\o_{a_2}(b_2,
 a_1]) [ a_2, a_1] \o_{a_1} ( \o_{a_2}([b_2, b_1]) [a_2,
 b_1])^{-1}\\
 &=&  \o_{a_1}([b_1, a_2]) \o_{a_1a_2}([b_1, b_2])
[ a_1, a_2] \o_{a_2}([a_1, b_2]).
\end{eqnarray*}
2. By statement 1, it suffices to show that $[A, \prod_{i=1}^m
B_i]=\prod_{i=1}^m [ A,B_i]$. The general case follows easily from
the case when $m=2$ (by induction). The case $m=2$ follows from
(\ref{a1b1}) where we put $b_1=1$, $a_1\in A$, $a_2\in B_1$, and
$b_2\in B_2$. $\Box$

\begin{theorem}\label{A18Mar9}
\begin{enumerate}
\item $\CG_1=\{ \s_{u,i,\pm 1}: x\mapsto ux^{\pm 1}, H\mapsto \pm
H-i\, | \, u\in \CL_1^*, i\in \Z\}$. \item $\CG_1 =\langle \zeta
\rangle \ltimes \{ s_i\, | \, i\in \Z\} \ltimes \{ \s_u\, | \,
u\in \CL_1^*\} \simeq \Z_2\ltimes \Z \ltimes (K^*\times \Z^{(\Z
)})$ and $Z(\CG_1)= \{ \s_\l \, | \, \l =\pm 1\}$. \item $\Inn
(\CA_1) =\{ s_i\, | \, i\in \Z\} \ltimes \mL_1^0$ where $\mL_1^0:=
\{ \s_u\, | \, u\in \CL_1^0\} $. \item $\Out (\CA_1) = \{
\zeta^\varepsilon \s_{\l H^i}\cdot \Inn (\CA_1)\, | \, \varepsilon
=0,1; \l \in K^*;  i\in \Z\} \simeq \Z_2\ltimes (K^*\times  \Z)$.
\item $[\Inn (\CA_1),\Inn (\CA_1)]=\mL_1^{00}$ where
$\mL_1^{00}:=\{ \s_u\, | \, u\in \CL_1^{00}\}\simeq \Z^{(\Z )}$.
\item $[\CG_1, \CG_1]=\{ s_{2i}\, | \, i\in \Z\} \ltimes \{ \s_u\,
| \, u\in K^{*2}\CL_1^0\langle -H^2\rangle \} $ and $\CG_1/[\CG_1,
\CG_1]\simeq \langle \zeta\rangle \times \frac{\langle s_1\rangle
}{\langle s_1^2\rangle }\times \frac{(K/K^{*2}\times \langle
H\rangle )}{\langle -H^2\rangle}.$
\end{enumerate}
\end{theorem}

{\it Proof}. 2. Note that $\s_{u,i, \pm 1}=
\zeta^{\varepsilon_\pm}s_i\s_{s_{-i}\zeta^{-\varepsilon_\pm}(u)}$
where $\varepsilon_+=0$ and $\varepsilon_-=1$, and so the first
equality follows from statement 1. The centre of the group $\CG_1$
is $\{ \s_\l \, | \, \l =\pm 1\}$ follows from (\ref{Hpm0}).

1. Let $\s \in \CG_1$ and $H'$ be the skew direct product of
groups in statement 2. We have to show that $\s \in H'$. The
automorphism $\s$ can be extended to an automorphism of the skew
field $\Frac (\CA_1)$. Since $\deg_x(\s (F^*))= \Z$ and
$\deg_{x^{-1}}(\s (F^*))= \Z$ where $F^*:= \Frac (\CA_1)
\backslash \{ 0\}$, we must have $\s (x) = ax^{-1} +b+cx$ for some
elements $a,b,c\in \CL_1$. Since $\s (x)$ is a unit, we must have
either $\s (x) = ux$ or $\s (x) = ux^{-1}$ for some element $u\in
\CL_1^*$. In both cases, we can find an element, say $\tau$, of
the group $H'$ such that $\tau \s (x) = x$, i.e. $\tau \s \in \St
(x) = \{ s_i\, | \, i\in \Z\}$, and so $\s \in H'$.

3. The algebra $\CA_1$ is central, i.e. $Z(\CA_1) = K$, and so
 $\Inn (\CA_1) \simeq  \CA_1^*/ K^*\simeq  \{ x^i\,
 | \, i\in \Z\}\ltimes\CL_1'$. For each $u\in \CL_1'$ and $i\in \Z$, the inner
 automorphism of the algebra $\CA_1$ generated by the element
 $ux^i$ is equal to $\o_{ux^i} = \s_{\frac{u}{\s_1(u)}, i, +1}$.
 Note that for each element $m\in \Z\backslash \{ 0 \}$,
\begin{equation}\label{Hpm}
H+m= \s_1^{-m}(H)=H\cdot \prod_{i=0}^{m-1}
\frac{\s_1^{-m+i}(H)}{\s_1^{-m+i+1}(H)}=H\cdot
\frac{\prod_{i=0}^{m-1} \s_1^{-m+i}(H)}{\s_1(\prod_{i=0}^{m-1}
\s_1^{-m+i}(H))}.
\end{equation}
So, an element $v\in \CL_1'$ is of type $\frac{u}{\s_1(u)}$ for
some element $u\in \CL_1'$ iff $\deg_H(v)=0$ iff $v\in \CL_1^0$.
Now, it is obvious that $\Inn (\CA_1) = \{ \o_{ux^i}\, | \, u\in
\CL_1', i\in \Z\}= \{ s_i\, | \, i\in \Z \} \ltimes \{ \s_u\, | \,
u\in \CL_1^0\}= \{ s_i\, | \, i\in \Z\}\ltimes \mL_1^0$.

4. Statement 4 follows from statements 2 and 3.

5. By statement 3, the group $\Inn (\CA_1)$ is the semi-direct
product of two abelian groups. By Lemma \ref{a6Apr9}.(1), the
commutant $C$ of the group $\Inn (\CA_1)$ is generated by the
commutators $[s_i, \s_u]=\s_{\frac{\s_1^i(u)}{u}}$ where $i\in \Z$
and $u\in \CL_1^0$. Repeating the arguments of the proof of
statement 3 where the element $H$ is replaced by $\frac{H}{H+1}$,
we see that $C=\mL_1^{00}$.

6. Let us prove the equality. The commutant of the group $H_1:=
\langle \zeta \rangle \ltimes \{ s_i\, | \, i\in \Z\}$ is $\{
s_{2i}\, | \, i\in \Z\}$, by (\ref{Hpm0}) and Lemma
\ref{a6Apr9}.(1). By statement 2, the group $\CG_1$ is the
semi-direct product $H_1\ltimes \mL_1^*$. By Lemma
\ref{a6Apr9}.(1), $[\CG_1, \CG_1]=\{ s_{2i}\, | \, i\in
\Z\}\ltimes [ H_1, \mL_1^*]$ since the group $\mL_1^*$ is abelian.
It remains to show that $[H_1, \mL_1^*]=\{ \s_u \, | \, u\in
K^{*2}\CL_1^0\langle -H^2\rangle \}$. The group $H_1$ is the
disjoint union $\{ s_i\, | \, i\in \Z\} \cup \{ \zeta s_i\, | \,
i\in \Z\}$. The commutators  $[s_i, \s_u]=
\s_{\frac{s_i(u)}{u}}=\s_{\frac{\s_1^i(u)}{u}}$ where $i\in \Z$
and $u\in \CL_1^*$ generate the group $\mL_1^0$. The group
$\mL_1^*$ is abelian and $\mL_1^0\subseteq \mL_1^*$. The
commutator
$$ [\zeta s_i, \s_u]=\s_{\zeta \s_1^{i-1}(u^{-1})u^{-1}}=
\s_{\s_1^{-i+1}\zeta (u^{-1})u^{-1}}$$ is an element of the group
$\mL_1^*$. Now,
$$[\zeta s_i, \s_u] \equiv
\s_{\frac{\s_1^{-i+1}\zeta (u^{-1})}{\zeta (u^{-1})}\cdot \zeta
(u^{-1})u^{-1} }\equiv \s_{\zeta (u^{-1})u^{-1}}\mod \mL_1^0.$$
The element $u\in \CL_1^* = K^*\times \CL_1'$ is a unique product
$u=\l v$ where $\l \in K^*$ and $v\in \CL_1'$. The groups
$\mL_1^*$ and $\CL_1^*$ are isomorphic via $\s_u\mapsto u$. The
map $\CL_1^*\ra \CL_1^*$, $ u\mapsto \zeta (u^{-1})u^{-1}$, is a
group homomorphism since the group $\CL_1^*$ is abelian. When
$u=\l$, $\zeta (\l^{-1})\l^{-1}= \l^{-2}$, and so $\{ \s_\mu , \,
| \, \mu \in K^{*2}\}\subseteq [ H_1, \mL_1^*]$. For an arbitrary
element $u= \l v$,
$$\zeta (u^{-1})u^{-1}\equiv \zeta
(v^{-1})v^{-1}\equiv (-H^2)^{-d}\mod K^{*2}\times \CL_1^0,$$ where
$d=\deg_H(v)$. Therefore, $[H_1, \mL_1^*]= \{ \s_u\, | \, u\in
K^{*2}\CL_1^0\langle -H^2\rangle \}$, as required. It follows from
the equality for the commutant of the group $\CG_1$ that in
statement 6 the isomorphism holds.   $\Box $

{\bf The group of units $(1+F)^*$ and $\mS_1^*$}. Recall that the
 algebra (without 1) $F=\bigoplus_{i,j\in \N} KE_{ij}$ is the union
$M_\infty (K) := \bigcup_{d\geq 1}M_d(K)= \varinjlim M_d(K)$ of
the matrix algebras $M_d(K):= \bigoplus_{1\leq i,j\leq
d-1}KE_{ij}$, i.e. $F= M_\infty (K)$. For each $d\geq 1$, consider
the (usual) determinant $\det_d=\det : 1+M_d(K)\ra K$, $u\mapsto
\det (u)$. These determinants determine the (global) {\em
determinant}, 
\begin{equation}\label{gldet}
\det : 1+M_\infty (K)= 1+F\ra K, \;\; u\mapsto \det (u),
\end{equation}
where $\det (u)$ is the common value of all determinants
$\det_d(u)$, $d\gg 1$. The (global) determinant has usual
properties of the determinant. In particular, for all $u,v\in
1+M_\infty (K)$, $\det (uv) = \det (u) \cdot \det (v).$ It follows
from this equality and the  Cramer's formula for the inverse of a
matrix that the group $\GL_\infty (K):= (1+M_\infty (K))^*$ of
units of the monoid $1+M_\infty (K)$ is equal to 
\begin{equation}\label{GLiK}
\GL_\infty (K) = \{ u\in 1+M_\infty (K) \, | \, \det (u) \neq 0\}.
\end{equation}
Therefore, 
\begin{equation}\label{1GLiK}
(1+F)^* = \{ u\in 1+F \, | \, \det (u) \neq 0\}=\GL_\infty (K).
\end{equation}
 The kernel
$$\SL_\infty (K):= \{ u\in
\GL_\infty (K)\, | \, \det (u) =1\}$$
 of the group epimorphism
$\det : \GL_\infty (K)\ra K^*$ is a {\em normal} subgroup of
$\GL_\infty (K)$. The group $\GL_\infty (K)= U(K)\ltimes E_\infty
(K)$ is the semi-direct product of its two subgroups where
$U(K):=\{ \mu(\l ) :=\l E_{00}+1-E_{00}\, | \, \l \in K^*\}\simeq
K^*$, $\mu (\l ) \lra \l $, and $E_\infty (K)$ is the subgroup of
$\GL_\infty (K)$ generated by all the elementary matrices $1+\l
E_{ij}$ where $\l \in K$ and $i\neq j$. It is obvious that
$E_\infty (K) = \SL_\infty (K) $, $E_\infty (K) = [E_\infty (K) ,
E_\infty (K) ] = [ \GL_\infty (K) , \GL_\infty (K) ] $ and $\det
(\mu (\l )) = \l$.

For all elements $\s \in \mT^1\times \mU_1$ and $\mu (\l ) \in
U(K)$, $\s (\mu (\l )) = \mu (\l )$. Therefore, 
\begin{equation}\label{soml}
\s \o_{\mu (\l )}= \o_{\mu (\l )}\s
\end{equation}
since $\s \o_{\mu (\l )}= \s \o_{\mu (\l )}\s^{-1}\cdot \s =
\o_{\s (\mu (\l ))}\s = \o_{\mu (\l )}\s$.

\begin{theorem}\label{a13Dec8}
{\rm (Theorem 4.6, \cite{shrekaut})}
\begin{enumerate}
\item $\mS_1^*=K^*(1+F)^*\simeq K^*\times \GL_\infty (K)$.\item
$Z(\mS_1^*) = K^*$ and $Z((1+F)^*)=\{ 1 \}$.  \item $\Inn
(\mS_1)\simeq \GL_\infty (K)$, $\o_u\lra u$.
\end{enumerate}
\end{theorem}

\begin{theorem}\label{17Mar9}
\begin{enumerate}
\item $\mG_1= (\mT^1\times \mU_1) \ltimes \Inn (\mS_1)$. \item
$\mG_1\simeq (\mT^1\times \Z^{(\Z )})\ltimes \GL_\infty (K)$.
\item $[\mG_1, \mG_1] \simeq \{ \o_u \, | \, u\in E_\infty (K)\}$
and $\mG_1/[\mG_1, \mG_1]\simeq \mT^1\times \mU_1\times U(K)\simeq
K^*\times \Z^{(\Z )}\times K^*$. \item The group $\mG_1$ is
generated by the following elements: $t_\l$ where $\l \in K^*$,
$\mu_u$ where $u\in \{ (H+i), (H-j)_1\, | \, i\in \N, j\in \N
\backslash \{ 0\} \}$, and $\o_v$ where $v\in (1+F)^*\simeq
\GL_\infty (K)$.
\end{enumerate}
\end{theorem}

{\it Proof}. 1. Let $\s \in \mG_1$. By Lemma \ref{x23Mar9},
 $(\mT^1\times \mU_1) \ltimes \Inn (\mS_1)\subseteq \mG_1$. It
 remains to show that the reverse inclusion holds, that is $\s \in (\mT^1\times \mU_1) \ltimes \Inn
 (\mS_1)$.  By Lemma \ref{e21Mar9}, $\s (F)
= F$, and so the map
$$ \overline{\s}: \CA_1=\mA_1/F\ra  \CA_1=\mA_1/F, \;\; \oa = a+F\mapsto \s (a)+F,$$
is an isomorphism of the skew Laurent polynomial algebra $\CA_1=
\CL_1[x,x^{-1}; \s_1]$ where $\s_1 (H)= H-1$. By Theorem
\ref{A18Mar9}.(1), either  either $ \overline{\s}(y) =
\overline{\l } x^{-1}$ or, otherwise, $ \overline{\s} (y) =
\overline{\l} x$ for some element  $\overline{\l} \in \CL_1^*$.
Equivalently, either $\s (y) = \l y +f$ or $\s (y) = \l x +f$ for
some element $\l \in K^*\times \CH $ (see (\ref{defH1})) and
$f\in F$. By Lemma \ref{c24Mar9}, the second case is impossible
since, by (\ref{indxy}),
$$ 1=\ind (y) = \ind (\s (y))= \ind (\l x+f) = -\deg_x(\overline{\l} x) =
-1.$$ Therefore, $\s (y) = \l y +f$. The element $\l \in K^*\times
\CH$ is a unique product $\l = \nu u$ where $\nu \in K^*$ and $
u\in \CH$. Then, $t_\nu \mu_u\s (y) = y+g$ where $g:= t_\nu
\mu_u(f) \in F$ since $t_\nu \mu_u(F) = F$ (Lemma \ref{e21Mar9}).
Fix a natural number $m$ such that $g\in \sum_{i,j=0}^mKE_{ij}$.
Then the finite dimensional vector spaces
$$ V:= \bigoplus_{i=0}^m Kx^i\subset V':=
\bigoplus_{i=0}^{m+1}Kx^i$$ are $y'$-invariant where
$y':=t_\nu\mu_u \s (y)= y+g$. Note that $y'*x^{m+1} = y*x^{m+1}=
x^m$ since $g*x^{m+1}=0$.
 Note that $P_1=\bigcup_{i\geq 1}\ker (y^i)$ and $\dim \, \ker_{P_1} (y)
 =1$. Since the $\mA_1$-modules $P_1$ and ${}^{t_\nu \mu_u \s}P_1$ are
 isomorphic (Proposition \ref{J15Ma7}), $P_1=\bigcup_{i\geq 1}\ker (y'^i)$ and $\dim \, \ker_{P_1} (y')
 =1$. This implies that the elements $x_0', x_1', \ldots , x_m',
 x^{m+1}$ are a $K$-basis for the vector space $V'$ where
 $$ x_i':= y'^{m+1-i}*x^{m+1}, \;\; i=0,1, \ldots , m; $$
 and the elements  $x_0', x_1', \ldots , x_m'$ are a $K$-basis for the vector space
 $V$. Then the elements
 $$ x_0', x_1', \ldots , x_m', x^{m+1}, x^{m+2},  \ldots  $$
are a $K$-basis for the vector space $P_1$. The $K$-linear map
\begin{equation}\label{chm}
\v : P_1\ra P_1, \;\; x^i\mapsto x_i'\; (i=0,1,\ldots , m), \;
x^j\mapsto x^j\; (j>m),
\end{equation}
belongs to the group $(1+F)^*= \GL_\infty (K)\simeq \Inn (\mS_1)$
(by Theorem \ref{a13Dec8}) and satisfies the property that $ y'\v
= \v y, $ the equality is in $\End_K(P_1)$. This equality can be
rewritten as follows:
$$ \o_{\v^{-1}}t_\nu\mu_u \s (y) = y \;\; {\rm where}\;\;
\o_{\v^{-1}}\in \Inn (\mS_1).$$ By Theorem \ref{21Mar9}, $\s =
t_{\nu^{-1}}\mu_{u^{-1}}\o_\v \in (\mT^1\times \mU_1)\ltimes \Inn
(\mS_1)$, as required.

2. Statement 2 follows from statement 1 since $\mU_1\simeq \Z^{(\Z
)}$ (by (\ref{UnZZ})) and $\Inn (\mS_1) \simeq \GL_\infty (K)$ (by
 Theorem \ref{a13Dec8}).

3. Let us identify the groups $\Inn (\mS_1)$ and $\GL_\infty (K)$
via $\o_u\lra u$ (Theorem \ref{a13Dec8}.(3)). By statement 1,
$\mG_1= (\mT^1\times \mU_1)\ltimes \GL_\infty (K)$. By Lemma
\ref{a6Apr9},
\begin{eqnarray*}
[\mG_1, \mG_1]&=& [\mT^1\times \mU_1, \GL_\infty (K)]\cdot [
\GL_\infty (K), \GL_\infty (K)]= [ \mT^1\times \mU_1, U(K)\ltimes
E_\infty (K)]\cdot E_\infty (K)\\ &=& [ \mT^1\times \mU_1,
U(K)]\cdot
E_\infty (K)= E_\infty (K)\;\;\;\; ({\rm by}\;\; (\ref{soml})),\\
 \mG_1/[\mG_1,
\mG_1]&\simeq & (\mT^1\times \mU_1\times U(K))\ltimes E_\infty (K)
/ E_\infty (K) \simeq \mT^1\times \mU_1\times U(K)\simeq K^*\times
\Z^{(\Z )}\times K^*.
\end{eqnarray*}

4. Statement 4 follows from statement 1. $\Box $

\begin{corollary}\label{d24Mar9}
Each automorphism $\s$ of the algebra $\mA_1$ is a unique product
$\s = t_{\nu^{-1}}\mu_{u^{-1}}\o_\v$ where $\s (y) \equiv \nu u y
\mod F$ and $\v\in (1+F)^*= \GL_\infty (K)$ is defined  in
(\ref{chm}).
\end{corollary}

{\it Proof}. The result was established in the proof of Theorem
\ref{17Mar9}.(1) apart from the uniqueness of $\v$ which follows
from the fact that the centre of the group $(1+F)^* = \GL_\infty
(K)$ is $\{ 1\}$, Theorem \ref{a13Dec8}.(2). $\Box $

\begin{corollary}\label{b6Apr9}
$\xi (\mG_1) = \{ \s_u\, | \, u\in \CL_1^*\}$ and $\ker (\xi ) =
\Inn (\mS_1)$.
\end{corollary}

{\it Proof}. Since $\xi : \mT^1\ltimes \mU_1\simeq \{ \s_u\, | \,
u\in \CL_1^*\}$ and $\Inn (\mS_1) \subseteq \ker (\xi )$, the
statements follow from the equality $\mG_1= (\mT^1\times \mU_1)
\ltimes \Inn (\mS_1)$ (Theorem \ref{17Mar9}.(1)). $\Box $

$\noindent $

  By Theorem 4.2, \cite{Bav-Jacalg}, the group of units of the algebra
$\mA_1$ is 
\begin{equation}\label{A1un}
\mA_1^* = K^* \times (\CH \ltimes (1+F)^*).
\end{equation}

It is a trivial observation that every algebra endomorphism of a
{\em simple} algebra is a {\em monomorphism}. The algebra $\mA_1$
is not simple but for it, as the following theorem shows,  the
same result is true as for the simple algebras.

\begin{theorem}\label{24Mar9}
Every algebra endomorphism of the algebra $\mA_1$ is a
monomorphism.
\end{theorem}

{\it Proof}. Recall  that $F$ is the  only proper  ideal of the
algebra $\mA_1$, and $\mA_1/F =\CA_1:=\CL_1\simeq K[x,x^{-1};
\s_1]$ is a simple algebra where $\s_1(H) = H-1$. Suppose that
$\g$ is an algebra endomorphism of $\mA_1$ which is not a
monomorphism, then necessarily $\g (F)=0$, and the endomorphism
$\g$ induces the algebra monomorphism $\overline{\g}: \CA_1\ra
\mA_1$, $a+F\mapsto \g (a)$. We seek a contradiction. Since $yx=1$
and $xy=1-E_{00}$, we have  the equalities $\g (y) \g (x) =1$ and
$\g (x) \g (y) =1$, i.e. the elements $\g (x) $ and $\g (y)$ are
units of the algebra $\mA_1$. The algebra $\mA_1$ is generated by
the elements $x$, $y$, and $H^{\pm 1}$, and so the {\em
noncommutative} simple algebra $\g (\mA_1)$ is generated by the
units $\g (x)$ , $\g (y)$, and $\g (H^{\pm 1})$ of the algebra
$\mA_1$. By (\ref{A1un}), the images of the elements $\g (x)$ ,
$\g (y)$, and $\g (H^{\pm 1})$ under the epimorphism $\pi
:\mA_1\ra \mA_1/ F$ belong to the commutative group $\pi (\mA_1^*)
= K^*\times \pi (\CH )$, and so they commute, this contradicts to
the fact that $\pi \overline{\g} (\CA_1)\simeq \CA_1$ is a
noncommutative algebra.
 $\Box $


\section{The group of automorphisms of the algebra $\CA_n:= \mA_n/ \ga_n$}\label{ATGRP}

 Let
$\CG_n:= \Aut_{K-{\rm alg}}(\CA_n)$. Since $\CA_n =
\bigotimes_{i=1}^n \CA_1(i) \simeq \CA_1^{\t n}$ where $\CA_1(i) =
\CL_1(i) [ x_i^{\pm 1}; \s_i]$, $\s_i (H_i) = H_i-1$, $\CL_1(i) =
K[H_i, (H_i+j)^{-1}]_{j\in \Z}$, the group $\CG_n$ contains the
symmetric group $S_n$ (elements of which permutes the tensor
multiples of the algebra $\CA_n$) and the direct product
$\prod_{i=1}^n \CG_1(i)$ where $\CG_1(i) := \Aut_{K-{\rm
alg}}(\CA_1(i))$. Moreover, $S_n\ltimes \prod_{i=1}^n
\CG_1(i)\subseteq \CG_n$. Our goal is to show that the equality
holds (Theorem \ref{5Apr9}.(1)). The group of units $\CA_n^*$ of
the algebra $\CA_n$ is equal to
$$ \CA_n^*= K^*\times (\mX_n\ltimes \CL_n')$$
where $\mX_n:= \{ x^\alpha \, | \, \alpha \in \Z^n\}$, $\CL_n':=
\prod_{i=1}^n \CL_1'(i)$, and $ \CL_1'(i) :=\{ \prod_{j\in \Z
}(H_i+j)^{n_{ij}}\, | \, (n_{ij})\in \Z^{(\Z )} \} \simeq  \Z^{(\Z
)}$.  For each number $i=1, \ldots , n$, $\CL_1^0(i)= \{ u\in
\CL_1'(i) \, | \, \deg_{H_i}(u)=0\}$ is the subgroup of
$\CL_1'(i)$. Then
$$ \CL_n^0:= \prod_{i=1}^n \CL_1^0(i)= \{ u\in \CL_n'\, | \,
\deg_{H_1}(u) = \cdots = \deg_{H_n}(u)=0\}$$ is the  subgroup of
$\CL_n'$. The commutant of the group $\CA_n^*$ is equal to
\begin{equation}\label{CAcom}
[\CA_n^*, \CA_n^*] =\CL_n^0.
\end{equation}
In more detail, by Lemma \ref{a6Apr9}, $[\CA_n^*, \CA_n^*]=
[\mX_n, \CL_n']=\prod_{i=1}^n [ \mX_1(i), \CL_1'(i)]=\prod_{i=1}^n
\CL_1^0(i) = \CL_n^0$. Let $A$ be a $K$-algebra and $Z(A)$ be its
centre. For each element $a\in A$, $\ad (a) : b\mapsto ab-ba$ is
the {\em inner derivation} of the algebra $A$ associated with the
element $a$. An element $a\in A \backslash Z(A)$ is called a {\em
strongly semi-simple} element of the algebra $A$ if $A=
\bigoplus_{\l \in E} \ker (\ad (a) - \l )$ where $E= \Ev (\ad (a)
, A)$ is the set of eigenvalues for the inner derivation $\ad (a)$
 that belong to the field $K$ . If $a$ is a strongly semi-simple element of
the algebra $A$ then so is the element $\s (a)$ for every
automorphism $\s $ of the algebra $A$.

\begin{theorem}\label{5Apr9}
\begin{enumerate}
\item $\CG_n = S_n\ltimes \prod_{i=1}^n \CG_1(i)$. \item $Z(\CG_n)
= \{ e, \s_{(-1, \ldots , -1)}\} $ where $\s_{(-1, \ldots , -1)}:
x_i\mapsto -x_i$, $H_i\mapsto H_i$, $i=1, \ldots ,n$. \item $\Inn
(\CA_n) \simeq \prod_{i=1}^n \Inn (\CA_1(i))$. \item $\Out (\CA_n)
\simeq S_n\ltimes \prod_{i=1}^n \Out (\CA_1(i))$. \item $[\Inn
(\CA_n) , \Inn (\CA_n)]= \prod_{i=1}^n[\Inn (\CA_1(i)) , \Inn
(\CA_1(i))]\simeq \prod_{i=1}^n \mL_1^{00}(i)$.\item $[\CG_n ,
\CG_n] = [S_n, S_n] \ltimes ( [S_n, \prod_{i=1}^n \CG_1(i)]\cdot
\prod_{i=1}^n [ \CG_1(i), \CG_1(i)])$ and $\frac{\CG_n}{[\CG_n ,
\CG_n]}\simeq \frac{S_n}{[S_n,S_n]}\times \frac{\CG_1}{[\CG_1 ,
\CG_1]}$, see Lemma \ref{a6Apr9} and Theorem  \ref{A18Mar9}.(6)
for details. If, in addition, the field $K$ is algebraically
closed then $\frac{\CG_n}{[\CG_n , \CG_n]}\simeq \begin{cases}
\Z_2^3& \text{if }n=1,\\
\Z_2^4& \text{if }n>1.\\
\end{cases}$

\end{enumerate}
\end{theorem}

{\it Proof}. 1. Let $\s \in \CG_n$. We have to show that $\s \in
S_n\ltimes \prod_{i=1}^n \CG_1(i)$. The restriction of the
automorphism $\s$ to the group $[\CA_n^*, \CA_n^*] = \CL_n^0$, see
(\ref{CAcom}), yields its automorphism. For each number $i=1,
\ldots , n$, $1+H_i^{-1}= \frac{H_i+1}{H_i}\in \CL_n^0$, and so
$\s ( 1+H_i^{-1})\in \CL_n^0\subset \CL_n^*$, i.e. $\s (H_i) \in
\CL_n^* = K^* \times \CL_n'$. The element $H_i$ is a {\em strongly
semi-simple} element of the algebra $\CA_n$ with $\Ev ( \ad (H_i)
, \CA_n) = \Z$. Therefore, $\s (H_i)$ is a strongly semi-simple
element of the algebra $\CA_n$ with $\Ev ( \ad (\s (H_i)), \CA_n)=
\Z$. It is not difficult to show that each strongly semi-simple
element of the set $\CL_n^*$ is of type $\l ( H_k+j)$ for some $\l
\in K^*$ and $j\in \Z$. Therefore, $\s (H_i) = \l_i
(H_{s(i)}+j_i)$, $i=1, \ldots , n$, for some $s\in S_n$, $\l_i\in
K^*$, and $j_i\in \Z$. Up to action of the group $S_n$ we may
assume that $s=e$, i.e. $\s (H_i) = \l_i (H_i+j_i)$. Then, up to
action of the group $\{ \o_v\, | \, v\in \mX_n\}$, see
(\ref{Hijxi}), we may assume that all $j_i=0$, i.e. $\s (H_i) =
\l_iH_i$. Now, $$\Z = \Ev ( \ad (\s (H_i)), \CA_n ) = \Ev (\l_i
\ad (H_i) , \CA_n) = \l_i \Ev ( \ad (H_i) , \CA_n) =\l_i \Z , $$
and so $\l_i= \pm 1$. The group $\CG_1(i)$ contains the
automorphism $\zeta_i : x_i\mapsto x_i^{-1}$, $H_i\mapsto - H_i$,
of order 2. Up to action of the group $\prod_{i=1}^n\langle
\zeta_i\rangle$, we may assume that all $\l_i =1$, i.e. $\s (H_i)
= H_i$. Then $\s^{-1}(H_i) = H_i$ for all $i$. For each element
$\alpha \in \Z^n$, $\CL_n x^\alpha = \bigcap_{i=1}^n \ker_{\CA_n}
(\ad (H_i) -\alpha_i)$. Therefore, $\s ( \CL_n x^\alpha )
\subseteq \CL_n x^\alpha $ and $\s^{-1} ( \CL_n x^\alpha )
\subseteq \CL_n x^\alpha $, i.e.  $\s ( \CL_n x^\alpha ) =\CL_n
x^\alpha $. Therefore, $\s (x_i) = u_ix_i$ for some element
$u_i\in \CL_n^*$, $i=1, \ldots , n$. By Theorem \ref{A18Mar9}.(2),
the group $\CG_1(i)$ contains the subgroup $\{ \s_{v_i}\, | \,
v_i\in \CL_1^*(i)\}$. Up to action of the group $\prod_{i=1}^n\{
\s_{v_i}\, | \, v_i\in \CL_1^*(i)\}$, we may assume that $u_i\in
K^* \prod_{j\neq i}^n \CL_1'(j)$. In particular, $\s_i(u_i) = u_i$
for all $i$. For each pair of indices $i\neq j$, the elements $\s
(x_i)$ and $\s (x_j)$ commute, and so
\begin{eqnarray*}
\s_i(u_iu_j)x_ix_j&= & u_i\s_i(u_j) x_ix_j= u_ix_iu_jx_j = \s (x_i) \s (x_j) = \s (x_j) \s (x_i) \\
 & = & u_jx_j u_ix_i =u_j\s_j (u_i) x_jx_i = \s_j (u_i u_j)
 x_ix_j,
\end{eqnarray*}
hence $\s_i ( u_iu_j) = \s_j(u_iu_j)$. By the very choice of the
elements $u_i$, this is possible iff $\s_i(u_iu_j) = u_iu_j$ and
$\s_j(u_iu_j) = u_iu_j$  iff $u_i\s_i(u_j) = u_iu_j$ and
$\s_j(u_i) u_j= u_iu_j$ iff $\s_i (u_j) = u_j$ and $\s_j(u_i) =
u_i$ iff all $u_i\in K^*$. Now, it is obvious that $\s \in
\prod_{i=1}^n \CG_1(i)$.

2. We assume that $n>1$ since for $n=1$ the statement is true
(Theorem \ref{A18Mar9}.(2)). For each element $s\in S_n$, the
restriction of the inner automorphism $\o_s: t\mapsto sts^{-1}$ of
the group $\CG_n$ to its normal subgroup $\prod_{i=1}^n \CG_1(i)$
permutes the components. Therefore, the centre $Z$ of the group
$\CG_n$  is a subgroup of $\prod_{i=1}^n \CG_1(i)$, and so
$Z\subseteq \prod_{i=1}^n Z( \CG_1(i))$. Now, statement 2 follows
from Theorem \ref{A18Mar9}.(2).

3. $\Inn (\CA_n) \simeq \CA_n^* / K^*\simeq \prod_{i=1}^n
(\CA_1^*(i)/ K^*) \simeq \prod_{i=1}^n \Inn (\CA_1^*(i))$.

4. $\Out (\CA_n) = \CG_n / \Inn (\CA_n) \simeq (S_n \ltimes
\prod_{i=1}^n \CG_1(i))/ \prod_{i=1}^n \Inn (\CA_1(i))\simeq S_n
\ltimes \prod_{i=1}^n \CG_1(i)/ \Inn (\CA_1(i))\simeq S_n \ltimes
\prod_{i=1}^n \Out (\CA_1(i))$.

5. Statement 5 follows from statement 3.

6. Statement 6 follows from statement 1, Lemma \ref{a6Apr9}, and
Theorem  \ref{A18Mar9}.(6).
 $\Box $


\section{The group  of automorphisms of the Jacobian  algebra $\mA_n$}\label{GGNAN}

In this section, the groups $\mG_n$, $\Inn (\mA_n)$ and $\Out
(\mA_n)$ are found (Theorem \ref{10Apr9}, Theorem \ref{B10Apr9},
and Corollary \ref{c10Apr9}). The image of the homomorphism $\xi$
(Theorem \ref{6Apr9}), its kernel (Corollary \ref{b18Apr9}), and
the group $\Aut_{\Z^n-{\rm gr}}(\mA_n)$ of the $\Z^n$-grading
preserving automorphisms of the algebra $\mA_n$ are found
(Corollary \ref{a11Apr9}), and the stabilizer $\St_{\mG_n}(\CH_1)$
of all the height one prime ideals of the algebra $\mA_n$  is
described (Corollary \ref{a12Apr9}). These results are used in
Section \ref{STABAN} where the stabilizers of all the ideals of
the algebra $\mA_n$ are found. It is proved that the groups
$\mG_n$ and $\ker (\xi )$ have trivial centre (Theorem
\ref{16Apr9}, Theorem \ref{A16Apr9}). Explicit inversion formulae
for automorphisms $\s \in \mG_n$ and $\s \in G_n$ are obtained,
(\ref{mGinv}) and (\ref{1mGinv}).

$\noindent $

{\bf A characterization of the elements of the group $\mG_n$}. For
each automorphism $\s $ of the algebra $\mA_n$, the next corollary
gives explicitly the map $\v \in \Aut_K(P_n)$ such that $\s =
\s_\v$ (see Corollary \ref{d21Mar9}). Corollary \ref{a18Apr9} is
used at the final stage of the proof of Theorem \ref{6Apr9} and
Theorem \ref{10Apr9}.

\begin{lemma}\label{a18Apr9}
For each automorphism $\s $ of the algebra $\mA_n$, there exists a
$K$-basis $\{ x'^\alpha\}_{\alpha \in \N^n}$ of the polynomial
algebra $P_n$ such that $\s (H_i) *x'^\alpha = (\alpha_i + 1)
x'^\alpha$ and $\s (y_i) *x'^\alpha = x'^{\alpha - e_i}$ for all
$i=1, \ldots , n$ (where $x'^\beta :=0$ if $\beta \in
\Z^n\backslash \N^n$). Then
\begin{enumerate}
\item  $\s =\s_\v$ where the map $\v \in \Aut_K(P_n)$: $x^\alpha
\mapsto x'^\alpha$ is the change-of-the-basis map, \item $\s (x_i)
*x'^\alpha = x'^{\alpha +e_i}$ for all $i=1, \ldots , n$, and
\item  the basis $\{ x'^\alpha\}_{\alpha \in \N^n}$ is unique up
to a simultaneous multiplication of each element of the basis by
the same nonzero scalar.
\end{enumerate}
\end{lemma}

{\it Proof}.  Recall that the polynomial algebra
$P_n=\bigoplus_{\alpha \in \N^n}Kx^\alpha$ is the direct sum of
non-isomorphic, one-dimensional, simple $K[H_1, \ldots ,
H_n]$-modules (see (\ref{kerHia})) such that $y_i*x^\alpha
=x^{\alpha - e_i}$ for all $\alpha \in \N^n$ and $i=1,\ldots , n$
(where $x^\beta :=0$ if $\beta \in \Z^n\backslash \N^n$). Recall
that $\s = \s_\v$ for some linear map $\v \in \Aut_K(P_n)$,  the
linear map $\v : P_n\ra {}^\s P_n$ is an $\mA_n$-module
isomorphism (Corollary \ref{d21Mar9}.(1)), and the map $\v$ is
unique up to a multiplication by a nonzero scalar since
$\End_{\mA_n} (P_n) \simeq K$ (Corollary \ref{y16Apr9}.(1)).
 Let $x'^\alpha := \v (x^\alpha )$ for $\alpha\in \N^n$. Then the
 fact that the map $\v$ is an $\mA_n$-module homomorphism is equivalent to
 the fact that the following equations hold:
 \begin{eqnarray*}
 \s (H_i)*x'^\alpha &=& \v H_i\v^{-1}\v *x^\alpha = (\alpha_i+1) \v *x^\alpha =(\alpha_i+1) x'^\alpha ,  \\
  \s (y_i)*x'^\alpha &=& \v y_i\v^{-1}\v *x^\alpha =  \v *x^{\alpha -e_i} = x'^{\alpha_i-e_i} ,  \\
   \s (x_i)*x'^\alpha &=& \v x_i\v^{-1}\v *x^\alpha =  \v *x^{\alpha +e_i} =
   x'^{\alpha_i+e_i}.
\end{eqnarray*}
Note that the last equality follows from the previous two: by the
first equality, the polynomial algebra $P_n=\bigoplus_{\alpha \in
\N^n}Kx'^\alpha$ is the direct sum of non-isomorphic,
one-dimensional, simple $K[\s (H_1),\ldots , \s (H_n)]$-modules.
Since $\s (H_i) \s (x_i) *x'^\alpha = \s (H_ix_i)*x'^\alpha = \s
(x_i(H_i+1))*x'^\alpha= \s (x_i) (\s (H_i) +1)*x'^\alpha =
(\alpha_i+2) \s (x_i)*x'^{\alpha}$ for all $i$, we have $\s (x_i)
*x'^\alpha = \l_{i,\alpha} x'^{\alpha +e_i}$ for  a scalar $
\l_{i,\alpha}$ which is necessarily equal to $1$ since
$$ x'^\alpha = \s (y_i)\s (x_i)*x'^\alpha = \l_{i,\alpha} \s (y_i)
*x'^{\alpha +e_i}= \l_{i,\alpha} x'^\alpha.$$Since the isomorphism
$\v$ is unique up to a multiplication by a nonzero scalar, the
basis $\{ x'^\alpha \}$ is unique up to a simultaneous
multiplication of each element of it by the same nonzero scalar.
The proof of the lemma is complete.  $\Box $

$\noindent $

Let $A=\bigoplus_{\alpha \in \Z^n}A_\alpha$ be a $\Z^n$-graded
algebra. Each nonzero element $a$ of $A$ is a unique sum $\sum
a_\alpha$ where $a_\alpha \in A_\alpha$. Define $|a|:= \max \{
|\alpha_1|, \ldots , |\alpha_n|\, | \, a_\alpha\neq 0\}$ and set
 $|0|:=0$. For all elements $a, b \in A$ and $\l \in K^*$,
 $|a+b|\leq \max \{|a|, |b|\}$, $|ab|\leq  |a|+|b|$, and $|\l a|=
 |a|$. In particular, $|a^m|\leq m|a|$ for all $m\geq 1$. In the
 proof of Theorem \ref{6Apr9}, we use the concept of $|\cdot |$ in
 the case of the algebra $\mA_n$.

For each natural number $d\geq 1 $, there is the decomposition
$K[x_i]= (\bigoplus_{j=0}^{d-1}Kx_i^j)\bigoplus (\bigoplus_{k\geq
d} Kx_i^k)$. The idempotents of the algebra $\mA_n$,
$p(i,d):=\sum_{j=0}^{d-1}E_{jj}(i)$ and $q(i,d):=1-p(i,d)$, are
the projections onto the first and the second summand
correspondingly. Since $P_n=\bigotimes_{i=1}^n K[x_i]$, the
identity map $1=\id_{P_n}$ on the vector space  $P_n$ is the sum
\begin{equation}\label{1pIq}
1=\bigotimes_{i=1}^n (p(i,d)+q(i,d))=\sum_{I\subseteq \{ 1, \ldots
, n\}}p(I,d)q(CI,d)
\end{equation}
of orthogonal idempotents where $p(I,d):= \prod_{i\in I}p(i,d)$,
$q(CI,d):= \prod_{i\in CI}q(i,d)$, $p(\emptyset , d):=1$, and
 $q(\emptyset , d):=1$. Each idempotent $p(I,d)q(CI,d)\in
 \mS_n\subset \mA_n \subset \End_K(P_n)$ is the projection onto
 the summand $P_n(I, d)$ in the decomposition for $P_n$,
\begin{equation}\label{PnSI}
P_n=\bigoplus_{I\subseteq \{ 1, \ldots , n\} }P_n(I,d), \;\;
P_n(I,d) :=\bigoplus \{ Kx^\alpha \, | \, \alpha_i<d, \; {\rm
if}\; i\in I; \alpha_j\geq d, \; {\rm if}\; j\in CI\}.
\end{equation}
In particular, the idempotent $q(\{ 1, \ldots , n\}, d)$ is the
projection onto the subspace $P_n(\{ 1, \ldots , n\}, d)=\bigoplus
\{ Kx^\alpha\, | \, {\rm all}\; \alpha_i\geq d\}$.

$\noindent $

{\bf A formula for the map $\v $ such that $\s = \s_\v$}. By
Theorem \ref{21Mar9}, each element $\s = \s_\v \in \mG_n$ is
uniquely determined by the elements $\s (y_1), \ldots , \s (y_n)$.
 In the proof of Theorem \ref{6Apr9}, an explicit formula for the
 map $\v$ is given, (\ref{vf18Apr9}), via the elements $\s (y_1),  \ldots , \s (y_n)$.

By the very definition, the group $\ker (\xi )$ contains precisely
all the automorphisms $\s \in \mG_n$ such that
$$ \s (x_i) \equiv x_i\mod \ga_n ,\;\;  \s (y_i) \equiv y_i\mod \ga_n , \;\; \s (H_i) \equiv H_i\mod \ga_n
,\;\; i=1, \ldots , n.$$

\begin{theorem}\label{6Apr9}
\begin{enumerate}
\item $\mG_n = S_n \ltimes (\mT^n \times \mU_n) \ltimes \ker (\xi
)$. \item $\im (\xi ) = S_n \ltimes \prod_{i=1}^n \{ \s_{u_i} \, |
\, u_i\in \CL_1^*(i)\}$ where $\s_{u_i}\in \CG_1(i): x_i\mapsto
u_ix_i$, $H_i\mapsto H_i$. \item $\ker (\xi ) \subseteq \Inn
(\mA_n)$.
\end{enumerate}
\end{theorem}

{\it Proof}. 1. Statement 1 follows from statement 2. Indeed,
suppose that statement 2 holds. Then the restriction of the
homomorphism $\xi$ to the subgroup $S_n\ltimes (\mT^n \times
\mU_n)$ of $\mG_n$ yields an isomorphism to $\im (\xi )$, and so
statement 1 follows from the short exact sequence of groups: $1\ra
\ker (\xi )  \ra \mG_n \ra \im (\xi ) \ra 1$.

2. For $n=1$, statement 2 is obvious (Corollary \ref{b6Apr9}). Let
$n\geq 2$. Let $R$ be the RHS of the equality in statement 2. For
each number $i=1, \ldots , n$, consider the subgroup $\G (i) :=
\langle \zeta_i\rangle\ltimes \{ s_{i,j}\, | \, j\in \Z\}$ of
$\CG_1(i)$ where $\zeta_i : x_i\mapsto x_i^{-1}$, $ H_i\mapsto
-H_i$; and $s_{i,j} : x_i\mapsto  x_i$, $H_i\mapsto H_i-j$. Then
$\G_n:= \prod_{i=1}^m \G (i)$ is the subgroup of $\prod_{i=1}^n
\CG_1(i)$. Using the descriptions of the groups $\CG_n$ and
$\CG_1(i)$ (Theorem \ref{5Apr9}.(1) and Theorem \ref{A18Mar9}.(2))
and (\ref{Hpm0}), we see that the group $\CG_n$ is equal to the
set theoretic product $R\cdot \G_n :=\{ r\g \, | \, r\in R, \g \in
\G_n\}$ with $R\cap \G_n = \{ e\}$. Since $\xi : S_n \ltimes
(\mT^n \times \mU_n)\simeq R$, we have
$$ \im (\xi ) = \im (\xi ) \cap \CG_n = \im (\xi ) \cap R\G_n =
R\cdot (\im (\xi ) \cap \G_n ).$$ All the elements $y_i\in
\End_K(P_n)$, $i=1, \ldots , n$, are locally nilpotent maps. For
all automorphisms $\s\in \mG_n$, the elements $\s (y_i)\in
\mA_n\subseteq \End_K(P_n)$, $i=1, \ldots , n$, are locally
nilpotent maps since ${}^\s P_n \simeq P_n$. Therefore,
$$\im (\xi ) \cap \G_n = \im (\xi ) \cap \Sh_n$$ where $\Sh_n := \prod_{i=1}^n
\{ s_{i,j}\, | \, j\in \Z\}$ is the subgroup of $\G_n$ (otherwise,
take $\s \in \im (\xi ) \cap \G_n\backslash \im (\xi ) \cap
\Sh_n$, then $\s (y_i) = x_i+a_i$ for some $i$ and $a_i\in \ga_n$,
but the element $x_i+a_i$ is not a locally nilpotent map as
$(x_i+a_i)^j*(x_1\cdots x_n)^k= x_i^j(x_1\cdots x_n)^k\neq 0$ for
all $j\geq 1$ and $k\gg 0$, a contradiction).

 It remains to show that  $\im (\xi ) \cap \Sh_n =
\{ e\}$. Let $\xi (\s ) \in \im (\xi ) \cap \Sh_n$ for an element
$\s \in \mG_n$. There exist elements $j := (j_1, \ldots , j_n) \in
\Z^n$ and $d\in \N$ such that 
\begin{equation}\label{dsHxy}
\s (H_i)-H_i - j_i, \; \s (x_i) - x_i,\;  \s(y_i) - y_i\in
\sum_{k=1}^n \mA_{n-1, k}\t (\sum_{s,t=0}^{d-1} KE_{st}(k))\;\;
{\rm for\; all}\;\; i,
\end{equation}
  where $\mA_{n-1, k}:=\bigotimes_{l\neq
k}\mA_1(l)$. We have to show that all $j_i=0$. In fact, it
suffices to show  only that all $j_i\leq 0$ (indeed, since
$\s^{-1} (H_i) \equiv  H_i-j_i\mod \ga_n$, then $-j_i\leq 0$, i.e.
$j_i=0$ for all $i$). Suppose that there exists an index, say $n$,
such that $j_n>0$. We seek a contradiction. For vectors $\alpha ,
\beta \in \N^n$, we write $\alpha \succeq \beta$ if $\alpha_1\geq
\beta_1, \ldots , \alpha_n\geq \beta_n$. By the choice of the
number $d$, for all elements $\alpha \in \N^n$ such that $\alpha
\succeq \underline{d}:=(d, \ldots , d)$ and for all $i=1, \ldots ,
n$,
\begin{eqnarray*}
 \s (H_i)*x^\alpha &=& (H_i+j_i) *x^\alpha= (\alpha_i+1+j_i) x^\alpha,  \\
 \s (x_i)*x^\alpha &=& x_i*x^\alpha= x^{\alpha+e_i}, \\
 \s (y_i)*x^\alpha &=& y_i*x^\alpha= x^{\alpha-e_i} \;\;\;\;
 (x^\beta:=0 \;\; {\rm if}\;\; \beta \in \Z^n\backslash \N^n).
\end{eqnarray*}
By Lemma \ref{a18Apr9}, for the automorphism $\s = \s_\v$ there
exists a {\em unique} $K$-basis $\{ x'^\alpha\}_{\alpha \in \N^n}$
of the polynomial algebra $P_n$ such that $x'^{\alpha +j}=
x^\alpha$ for all $\alpha \succeq \underline{d}$, and the map $\v
:P_n\ra {}^\s P_n$, $x^\alpha \mapsto x'^\alpha$ $(\alpha \in
\N^n)$, is an $\mA_n$-module isomorphism.  For each natural number
$t$, the set $C_t:=\{ \alpha \in \N^{n-1}\, \, | \, {\rm all}\;
\alpha_i<t\}$ contains $t^{n-1}$ elements. Let $c:= 2(d+j_n) (|\s
(y_n)|+1)$. Fix a natural number $s$ such that $s>c+d$, and let
$W_t:= (s, s, \ldots , s, d)+C_t$ (recall that $C_t\subseteq
\N^{n-1}\subset \N^n$) and $\Pi_t:= \bigoplus_{\alpha \in W_t}
Kx^\alpha$. For all elements $\alpha \in W_t$, $\s (H_n) *x^\alpha
= (d+j_n+1) x^\alpha$. It follows that the set $\Pi'_t:=
\sum_{k=1}^{d+j_n}\s (y_n^k)*\Pi_t$ is the direct sum
$\bigoplus_{k=1}^{d+j_n}\s (y_n^k)*\Pi_t$, and $\Pi_t'\subseteq
\Pi_t'':= \sum_{\alpha \in W_t'}Kx^\alpha$ where $W_t':=
((s-c,\ldots ,s-c)+\Pi_{t+2c})\times \{ 0,1, \ldots , d-1\}$.
Comparing the leading terms of both ends of the inequality,
$$ (d+j_n) t^{n-1}= \dim (\Pi_t') \leq \dim (\Pi_t'') =
d(t+2c)^{n-1}= dt^{n-1}+\cdots \;\; (t\geq 1), $$ we conclude that
$j_n\leq 0$, a contradiction. Therefore, all $j_i=0$, and the
proof of statement 2  is complete.

3. Let $\s \in \ker (\xi )$. Then $\xi (\s )=e \in \im (\xi ) \cap
\Sh_n= \{ e\}$, and we keep the notation of the proof of statement
2 for the automorphism $\s = \s_\v$ (where all $j_i=0$) . It
remains to show that $\v \in \mA_n$. This is obvious since

\begin{equation}\label{vf18Apr9}
\v = q(\{ 1, \ldots , n\}, d)+\sum_{\emptyset \neq I\subseteq  \{
1, \ldots , n\}}(\sum_{\alpha \in C_d(I)}\prod_{j\in I}\s
(y_j^{d-\alpha_j})\cdot \prod_{i\in I}x_i^{d-\alpha_i}\cdot
E_{\alpha \alpha}(I))p(I,d)q(CI, d)
\end{equation}
where $C_d(I):=\{ (\alpha_i)_{i\in I}\in \N^I\, | \, {\rm all }\;
\alpha_i<d\}$,  $E_{\alpha \alpha}(I):=\prod_{i\in
I}E_{\alpha_i\alpha_i}(i)$, and $d$ is from (\ref{dsHxy}). To
prove that this formula holds for the map  $\v$ we have to show
that $\v *x^\alpha = x'^\alpha$ for all $\alpha \in \N^n$. For
each monomial $x^\alpha$, let $I:= \{ i \, | \, \alpha_i<d\}$.
Then $x^\alpha = \prod_{i\in I}x_i^{\alpha_i}\cdot \prod_{k\in
CI}x_k^{\alpha_k}$ and, if $I\neq \emptyset$,
\begin{eqnarray*}
 \v *x^\alpha &=& \prod_{j\in I}\s (y_j^{d-\alpha_j})\cdot \prod_{i\in I}x_i^{d-\alpha_i}*\prod_{i\in
I}x_i^{\alpha_i}\cdot \prod_{k\in CI}x_k^{\alpha_k}= \prod_{j\in
I}\s (y_j^{d-\alpha_j})*\prod_{i\in
I}x_i^d\cdot \prod_{k\in CI}x_k^{\alpha_k}  \\
&=& \prod_{j\in I}\s (y_j^{d-\alpha_j})*x'^{(\sum_{i\in
I}de_i+\sum_{k\in CI}\alpha_ke_k)}=x'^{(\sum_{i\in
I}\alpha_ie_i+\sum_{k\in CI}\alpha_ke_k)}=x'^\alpha .
\end{eqnarray*}
If $I=\emptyset$ then $\v *x^\alpha = q(\{ 1, \ldots , n\} , d)
*x^\alpha = x^\alpha= x'^\alpha$.  The proof of the theorem is
complete. $\Box $

\begin{corollary}\label{d10Apr9}
 $\mG_n\supseteq \mG_n':= S_n\ltimes
(\mT^n\times \mU_n)\ltimes \underbrace{\GL_\infty (K)\ltimes\cdots
\ltimes \GL_\infty (K)}_{2^n-1 \;\; {\rm times}}$  and
$\mU_n\simeq (\Z^n)^{(\Z )}$.
\end{corollary}

{\it Proof}. Since (see Corollary \ref{a19Apr9}) $\ker  (\xi )
\supseteq \Inn (\mS_n)\supseteq \underbrace{\GL_\infty
(K)\ltimes\cdots \ltimes \GL_\infty (K)}_{2^n-1 \;\; {\rm
times}}$, \cite{shrekaut},  the statement follows from Theorem
\ref{6Apr9}.(1,3). $\Box $

\begin{corollary}\label{b18Apr9}
Let $\s = \s_\v\in \ker (\xi )$ where $\v \in \mA_n^*$ is from
(\ref{vf18Apr9}). Then
\begin{enumerate}
\item $\v \equiv 1 \mod \ga_n$, i. e. $\v \in (1+\ga_n)^*$, and so
the element $\v $ from (\ref{vf18Apr9}) is the {\em unique}
element $\v'$ such that $\s = \s_{\v'}$ and $\v'\equiv 1 \mod
\ga_n$; and \item $\ker (\xi )= \{ \o_\v \, | \, \v \in
(1+\ga_n)^* \} \simeq (1+\ga_n)^*$, $\o_\v \lra \v$.
\end{enumerate}
\end{corollary}

{\it Proof}. 1. By (\ref{vf18Apr9}), $\v \equiv q(\{1, \ldots ,
n\}, d)\equiv 1 \mod \ga_n$.

2. Statement 2 follows from statement 1.  $\Box $

$\noindent $

{\bf The canonical presentation of $\s \in \mG_n$}. By Theorem
\ref{6Apr9}.(1) and Corollary \ref{b18Apr9}, each automorphism
$\s$ of the algebra $\mA_n$ is a {\em unique} product $st_\l
\mu_u\o_\v$ where $s\in S_n$, $t_\l \in \mT^n$, $\mu_u\in \mU_n$,
and $\o_\v$ is an inner automorphism of the algebra $\mA_n$ with
$\v \in (1+\ga_n)^*$.

$\noindent $

{\it Definition}. The unique product $\s =st_\l \mu_u\o_\v$  is
called the {\em canonical presentation} of $\s $.

\begin{corollary}\label{c18Apr9}
Let $\s \in \mG_n$ and $ \s =st_\l \mu_u\o_\v$ be its canonical
presentation. Then the automorphisms $s$, $t_\l$, $\mu_u$, and
$\o_\v$ can be effectively (in finitely many steps) found from the
action of the automorphism $\s$ on the elements $\{ H_i, x_i,
y_i\, | \, i=1, \ldots , n\}$:
$$ \s (H_i) \equiv H_{s(i)}\mod \ga_n, \;\; \s (x_i) \equiv
x_{s_(i)}\l_is (u_i)\mod \ga_n, $$ and  the element $\v$ is given
by the formula (\ref{vf18Apr9}) for the automorphism
$(st_\l\mu_u)^{-1}\s \in \ker (\xi )$.
\end{corollary}

{\bf A formula for  the inverse map $\v^{-1}$ where  $\s= \s_\v\in
\ker (\xi )$}. For an  automorphism $\s = \s_\v \in \ker (\xi )$,
i.e. $\v \in (1+\ga_n)^*$ is as in (\ref{vf18Apr9}), we are going
to produce a formula for the inverse $\v^{-1}$ (see
(\ref{1vf18Apr9})) which is the most difficult part in finding
the inversion formula (\ref{mGinv}) for $\s \in \mG_n$. Since $\v
: P_n=\bigoplus_{\alpha \in \N^n} Kx^\alpha \simeq
P_n=\bigoplus_{\alpha \in \N^n} Kx'^\alpha$ and $x'^\alpha = \v
(x^\alpha )$, for each element $\alpha \in \N^n$, the projection
$E'_{\alpha \alpha}$ onto the summand $Kx'^\alpha$ is equal to $\v
E_{\alpha\alpha}\v^{-1} = \s (E_{\alpha\alpha})$.
 For each subset $I$ of the set $\{ 1,\ldots , n\}$, let
 $$P_n'(I,d):=\v (P_n(I,d))=\bigoplus \{ Kx'^\alpha \, | \,
 \alpha_i<d\; {\rm if}\; i\in I; \alpha_j\geq d\; {\rm if}\; j\in
 CI\}.$$ Since $\v : P_n=\bigoplus P_n(I,d)
\simeq P_n=\bigoplus P_n'(I,d)$, the projections onto the summand
$P'_n(I,d)$ is equal to
$$ \v p(I,d)q(I,d) \v^{-1} = \s (p(I,d)q(I,d)) = \s (p(I,d)) \s
(q(I,d))= p'(I,d)q'(I,d)$$ where $p'(I,d):= \s (p(I,d))$ and
$q'(I,d):=\s (q(I,d))$.

The inverse map $\v^{-1}$ of the map $\v$ from (\ref{vf18Apr9}) is
given by the rule: 
\begin{equation}\label{1vf18Apr9}
\v^{-1} = q'(\{ 1, \ldots , n\}, d)+\sum_{\emptyset \neq
I\subseteq \{ 1, \ldots , n\}}(\sum_{\alpha \in C_d(I)}\prod_{j\in
I}y_j^{d-\alpha_j}\cdot \prod_{i\in I}\s (x_i)^{d-\alpha_i}\cdot
E_{\alpha \alpha}'(I))p'(I,d)q'(CI, d)
\end{equation}
where $E_{\alpha \alpha}'(I):=\prod_{i\in
I}E_{\alpha_i\alpha_i}'(i)$ and $E_{jj}'(i):=\s (E_{jj}(i))=\s
(x_i)^j\s (y_i)^j-\s (x_i)^{j+1}\s (y_i)^{j+1}$;
$p'(I,d):=\prod_{i\in I}p'(i,d)$ and $p'(i,d):=  \s
(p(i,d))=\sum_{j=0}^{d-1}E_{jj}'(i)$;  $q'(CI, d):= \s (q(CI,
d))=\prod_{i\in CI}(1-p'(d,i))$. Let $\psi$ be the RHS of
(\ref{1vf18Apr9}). We have to show that $\psi : x'^\alpha \ra
x^\alpha$ for all $\alpha \in \N^n$. Fix $\alpha$, and let $I:= \{
i \, | \, \alpha_i<d\}$. Then $x'^\alpha = \prod_{i\in
I}x_i'^{\alpha_i}\cdot \prod_{k\in CI}x_k'^{\alpha_k}$ and, if
$I\neq \emptyset$ then,  by Lemma \ref{a18Apr9}.(2),
\begin{eqnarray*}
 \psi *x'^\alpha &=& \prod_{j\in I}y_j^{d-\alpha_j}\cdot \prod_{i\in I}\s (x_i)^{d-\alpha_i}*\prod_{i\in
I}x_i'^{\alpha_i}\cdot \prod_{k\in CI}x_k'^{\alpha_k}= \prod_{j\in
I}y_j^{d-\alpha_j}*\prod_{i\in
I}x_i'^d\cdot \prod_{k\in CI}x_k'^{\alpha_k}  \\
&=& \prod_{j\in I}y_j^{d-\alpha_j}*\prod_{i\in I}x_i^d\cdot
\prod_{k\in CI}x_k^{\alpha_k}=x^\alpha .
\end{eqnarray*}
If $I\neq \emptyset$ then $\psi *x'^\alpha = q'(\{ 1, \ldots , n\}
, d)*x'^\alpha = x'^\alpha = x^\alpha$. This finishes the proof of
(\ref{1vf18Apr9}).

$\noindent $

{\bf An inversion  formula for  $\s \in \mG_n$}. Let $\s \in
\mG_n$ and $\s = st_\l\mu_u\o_\v$ be its canonical presentation.
Then the inversion formula for $\s^{-1}$ is given by rule:
\begin{equation}\label{mGinv}
\s^{-1}= s^{-1} t_{s(\l^{-1})} \mu_{s(u^{-1})}\o_{st_\l \mu_u
(\v^{-1})},
\end{equation}
this is the canonical presentation for the automorphism $\s^{-1}$
where $\v^{-1}$ is given by (\ref{1vf18Apr9}). Indeed,
$$ \s^{-1}= s^{-1} \cdot s(t_\l \mu_u)^{-1}s^{-1}\cdot (st_\l
\mu_u)\o_{\v^{-1}}(st_\l \mu_u)^{-1} =s^{-1} t_{s(\l^{-1})}
\mu_{s(u^{-1})}\o_{st_\l \mu_u (\v^{-1})}. $$ {\bf The groups
$\mG_n$ and $\Inn (\mA_n)$}.

\begin{theorem}\label{B10Apr9}
$\Inn (\mA_n) = \mU^0_n\ltimes \ker (\xi )\simeq \mU^0_n\ltimes
(1+\ga_n)^*$.
\end{theorem}

{\it Proof}. It suffices to prove the equality since then the
isomorphism holds, by Corollary \ref{b18Apr9}.(2). Since $\mG_n =
S_n\ltimes (\mT^n \times \mU_n) \ltimes \ker (\xi ) $ and $\ker
(\xi ) \subseteq \Inn (\mA_n)$ (Theorem \ref{6Apr9}.(1,3)), we see
that $\Inn (\mA_n ) = \CI \ltimes \ker (\xi )$ where $\CI := \Inn
(\mA_n) \cap S_n\ltimes (\mT^n\times \mU_n)$. Since $\xi :
S_n\ltimes (\mT^n \times \mU_n)\simeq S_n\ltimes (\mT^n \times \xi
(\mU_n))$, we have $\xi : \CI \simeq \xi (\CI )$. By the very
definition of the group $\CI$,
\begin{eqnarray*}
\xi (\CI ) &\subseteq  & \Inn (\CA_n) \bigcap S_n\ltimes (\mT^n
\times \xi (\mU_n)) \\
&=& \prod_{i=1}^n \Inn (\CA_1(i))\bigcap
S_n\ltimes (\mT^n \times \xi (\mU_n))\;\; ({\rm Theorem}\; \ref{5Apr9}.(3))\\
 &=& \prod_{i=1}^n \Inn (\CA_1(i))\bigcap
\mT^n \times \xi (\mU_n)\;\; ({\rm Theorem}\; \ref{5Apr9}.(1, 3))\\
&=& \prod_{i=1}^n \Inn (\CA_1(i))\bigcap
\prod_{i=1}^n (\mT^1(i)\times \xi (\mU_1(i)))\\
&=& \prod_{i=1}^n (\Inn (\CA_1(i))\bigcap (\mT^1(i)\times \xi
(\mU_1(i))))\\
&=& \prod_{i=1}^n \xi (\mU_1^0(i))\;\; ({\rm Theorem } \;
\ref{A18Mar9}.(3))\\
&=& \xi (\prod_{i=1}^n \mU_1^0(i))= \xi (\mU_n^0),
\end{eqnarray*}
hence $\CI = \mU_n^0$ since $\CI , \mU_n^0\subseteq S_n \ltimes
(\mT^n \times \mU_n)$ and $\xi :S_n\ltimes (\mT^n\times
\mU_n)\simeq S_n\ltimes (\mT^n\times \xi (\mU_n))$. Therefore,
$\Inn (\mA_n) = \mU_n^0 \ltimes \ker (\xi )$.  $\Box $

\begin{theorem}\label{10Apr9}
$\mG_n = S_n\ltimes (\mT^n \times \Xi_n)\ltimes \Inn (\mA_n)$.
\end{theorem}

{\it Proof}. The statement follows from two facts: $\mG_n =
S_n\ltimes (\mT^n \times \Xi_n \times \mU_n^0)\ltimes \ker (\xi )
$ (Theorem \ref{6Apr9}.(1)) and $\Inn (\mA_n) = \mU_n^0\ltimes
\ker (\xi )$ (Theorem \ref{B10Apr9}). $\Box $

\begin{corollary}\label{c10Apr9}
$\Out (\mA_n) \simeq S_n\ltimes (\mT^n\times \Xi_n)$.
\end{corollary}


Recall that the set $\CH_1$ of height one  prime ideals  of the
algebra $\mA_n$ is $\{ \gp_1, \ldots , \gp_n\}$. The next
corollary describes its stabilizer,
$$ \St_{\mG_n} (\CH_1):=\{ \s \in \mG_n\, | \, \s (\gp_1) = \gp_1,
\ldots, \s (\gp_n) = \gp_n\}.$$

\begin{corollary}\label{a12Apr9}
$ \St_{\mG_n} (\CH_1)=(\mT^n\times \Xi_n)\ltimes \Inn (\mA_n)=
(\mT^n\times \mU_n) \ltimes \ker (\xi )$.
\end{corollary}

{\it Proof}. The second equality is obvious since
$\mU_n=\Xi_n\times \mU_n^0$ and $\Inn (\mA_n) = \mU_n^0\ltimes
\ker (\xi )$ (Theorem \ref{B10Apr9}). It remains to show that the
first equality holds, that is $S=R$ where $S:=\St_{\mG_n} (\CH_1)$
and $R:=(\mT^n\times \Xi_n)\ltimes \Inn (\mA_n)$. The inclusion
$R\subseteq S$ follows from (\ref{tlEab}) and (\ref{muEij}).
Clearly, $S\cap S_n=\{ e\}$. By Theorem \ref{10Apr9}, $S=S\cap
\mG_n = S\cap (S_n\ltimes R)= (S\cap S_n)\cdot R= \{ e\} \cdot R =
R$. $\Box$

The algebra $\mA_n =  \bigoplus_{\alpha \in \Z^n} \mA_{n , \alpha
}$ is a $\Z^n$-graded algebra. Let $\Aut_{\Z^n-{\rm gr}}(\mA_n):=
\{ \s \in \mG_n\, | \, \s (\mA_{n , \alpha }) = \mA_{n , \alpha }$
for all $\alpha \in \Z^n\}$ be the group of automorphisms of the
algebra $\mA_n$ that respect the $\Z^n$-grading, and
$\St_{\mG_n}(H_1, \ldots, H_n) :=\{ \s \in \mG_n\, | \, \s (H_1) =
H_1, \ldots , \s (H_n) = H_n\}$.

\begin{corollary}\label{a11Apr9}
\begin{enumerate}
\item  $\St_{\mG_n}(H_1, \ldots, H_n) = \mT^n \times \mU_n$. \item
$\Aut_{\Z^n-{\rm gr}}(\mA_n)=\St_{\mG_n}(H_1, \ldots, H_n)$.
\end{enumerate}
\end{corollary}

{\it Proof}. Note that $\mT^n \times \mU_n\subseteq
S:=\St_{\mG_n}(H_1, \ldots, H_n)$ and $S\subseteq \Aut_{\Z^n-{\rm
gr}}(\mA_n)$ since $\mA_{n, \alpha}=\bigcap_{i=1}^n \ker_{\mA_n}
(\ad (H_i) -\alpha_i)$ for all $\alpha \in \Z^n$. To finish the
proof of both statements it remains to show that $\Aut_{\Z^n-{\rm
gr}}(\mA_n)\subseteq \mT^n \times \mU_n$. Let $\s \in
\Aut_{\Z^n-{\rm gr}}(\mA_n)$. For each number $i=1, \ldots , n$,
$\mA_{n, e_i}= x_i\mD_n$ and $\mA_{n, -e_i}= \mD_ny_i$, and so $\s
(x_i) = x_iu_i$ and $\s (y_i) = v_iy_i$ for some elements $u_i,
v_i\in \mD_n$ such that $1= \s (1) = \s (y_ix_i) = v_iy_ix_iu_i =
v_iu_i$, and so $v_i= u_i^{-1}$ (as $\mD_n$ is a commutative
algebra). Now, $\mu_u^{-1}\s (x_i) = x_i$ for all $i$ where
$\mu_u\in \mT^n \times \mU_n$ and $u=(u_1, \ldots, u_n)$. By
Theorem \ref{21Mar9}, $\mu_u^{-1}\s = e$, and so $\s = \mu_u\in
\mT^n \times \mU_n$, as required.  $\Box $

$\noindent $

{\bf The inner canonical presentation of $\s \in \mG_n$}. By
Theorem \ref{10Apr9}, each automorphism $\s $ of the algebra
$\mA_n$ is a {\em unique} product $st_\l\mu_{H^\alpha}\o_w$ where
$s\in S_n$, $t_\l \in \mT^n$, $\mu_{H^\alpha}\in \Xi_n$ ($\alpha
\in \N^n$), and $\o_w\in\Inn (\mA_n)$ where $w\in \CH_n^0\ltimes
(1+\ga_n)^*$, by Theorem \ref{B10Apr9}, Corollary
\ref{b18Apr9}.(2), and (\ref{degHi1}).

$\noindent $

{\it Definition}. The unique product $\s =st_\l\mu_{H^\alpha}\o_w$
is called the {\em inner canonical presentation} of $\s $.

$\noindent $

Corollary \ref{a29Jun9} shows  how to find the inner canonical
presentation of $\s$ effectively (in finitely many steps) via the
elements $\{ \s (H_i), \s (x_i), \s (y_i)\, | \, i=1, \ldots ,
n\}$. First, we need to introduce more notation. For each number
$i=1, \ldots , n$, the map
$$ \psi_i : \CH_1(i)\ra \CH_1^0(i), \;\; u\mapsto
\frac{\tau_i(u)}{u}, $$is a group  isomorphism where
$\tau_i((H_i+j)_*)=(H_i+j+1)_*$ for all $j\in \Z$ where $
(H_i+j)_*:=\begin{cases}
(H_i+j)_1& \text{if } j<0,\\
H_i+j& \text{if }j\geq 0.\\
\end{cases}$. The inverse map $\psi_i^{-1}$ is given by the
formula 
\begin{equation}\label{pinsi}
\psi_i^{-1}(\prod_{j\in \Z}(H_i+j)_*^{n_{ij}})=\prod_{j\in
\Z}(H_i+j)_*^{\sum_{k>j}n_{ik}}.
\end{equation}
Therefore, the map $\psi:=\prod_{i=1}^n\psi_i :
\CH_n=\prod_{i=1}^n\CH_1(i)\ra \CH_n^0=\prod_{i=1}^n\CH_1^0(i)$ is
a group isomorphism  and the inverse map
$\psi^{-1}=\prod_{i=1}^n\psi_i^{-1}$ is determined by
(\ref{pinsi}). Let us show that 
\begin{equation}\label{pinsi1}
\mu_{\psi (u)}=\o_u, \;\; u\in \mU_n.
\end{equation}
For each $i=1, \ldots , n$, $\o_u(x_i)=ux_iu^{-1}=
x_i\frac{\tau_i(u)}{u}=x_i\frac{\tau_i(u_i)}{u_i}=x_i\psi_i(u_i)=\mu_{\psi
(u)}(x_i)$ where $u=\prod_{i=1}^nu_i$, and  $u_i\in \CH_1(i)$.
Therefore, $\o_u = \mu_{\psi (u)}$, by Theorem \ref{21Mar9}.
\begin{corollary}\label{a29Jun9}
Let $\s \in \mG_n$ and $\s = st_\l \mu_u\o_\v$ be its canonical
presentation (Lemma \ref{c18Apr9}),  $u=H^\alpha v$ where $\alpha
\in \N^n$,  $v\in \CH_n^0$. Then  $\s = st_\l
\mu_{H^\alpha}\o_{\psi^{-1}(v)\v }$ is the inner canonical
presentation of $\s$.
\end{corollary}

{\it Proof}. By (\ref{pinsi1}), $\mu_v=\o_{\psi^{-1}(v)}$. Then
$\s = st_\l \mu_{H^\alpha}\mu_v\o_\v= st_\l
\mu_{H^\alpha}\o_{\psi^{-1}(v)\v }$.  $\Box $

$\noindent $

{\bf The groups $\mA_n^*$ and $\mS_n^*$ and their centres}. Since
$\mS_n\subseteq \mA_n$, we have the inclusion $\mS_n^*\subseteq
\mA_n^*$.

\begin{theorem}\label{24Ma7}
{\rm (Theorem 4.4, \cite{Bav-Jacalg})}
\begin{enumerate}
\item $\mA_n^*= K^*\times (\CH_n\ltimes (1+\ga_n)^*)$ where
 the group $\CH_n$ is defined in (\ref{defH2}) and $(1+\ga_n)^*$
 is the group of units of $\mA_n$  the type $1+a$ for some element $a\in \ga_n$.
 \item The centre of the group $\mA_n^*$ is $K^*$.
\end{enumerate}
\end{theorem}
As a consequence of this theorem we obtain a description of the
group $\mS_n^*$ of units of the algebra $\mS_n$ (Corollary
\ref{a19Apr9}) which is used in finding explicit generators for
the group $G_n$ in \cite{Snaut}.

\begin{corollary}\label{a19Apr9}
\begin{enumerate}
\item $\mS_n^* = K^*\times (1+\ga_n)^*$ where the ideal $\ga_n$ of
the algebra $\mS_n$ is the sum of all the height 1 prime ideals of
the algebra $\mS_n$ (see (\ref{SnSn})).  \item  The centre of the
group $\mS_n^*$ is $K^*$, and the centre of the group
$(1+\ga_n)^*$ is $\{ 1 \}$. \item The map $(1+\ga_n)^*\ra \Inn
(\mS_n)$, $u\mapsto \o_u$, is a group  isomorphism.
\end{enumerate}
\end{corollary}

{\it Proof}. We denote, for a moment, the ideal $\ga_n$ of the
algebra $\mS_n$ by $\ga_n'$ in order to distinguish it from the
only maximal ideal $\ga_n$  of the algebra $\mA_n$. Since
$\ga_n'\subseteq \ga_n$, and, by (\ref{SnSn}),
$\mS_n/\ga_n'=L_n\subset \mA_n/\ga_n$, we have the equality
$\ga_n'=\mS_n\cap \ga_n$ which implies the equality $(1+\ga_n')^*
= \mS_n\cap (1+\ga_n)^*= \mS_n^*\cap (1+\ga_n)^*$. Since
$(\mS_n/\ga_n')^*\cap (\mA_n/\ga_n)^*\cap (K^*\times \xi
(\CH_n))=K^*$, we see that
\begin{eqnarray*}
 \mS_n^*&=& \mS_n^*\cap \mA_n^* = \mS_n^*\cap (K^*\times (\CH_n\ltimes (1+\ga_n)^*))= K^*\times (\mS_n^*\cap (1+\ga_n)^*) \\
 &=& K^*\times (1+\ga_n')^*. \;\;\;
\end{eqnarray*}
2. By statement 1, it suffices to show that $Z((1+\ga_n)^*)=\{
1\}$. For $n=1$, this was done in \cite{shrekaut}, Theorem 4.6.
So, let $n>1$, and we use induction on $n$. Let $u\in
Z((1+\ga_n)^*)$. For each integer $i=1, \ldots , n$, the element
$u+\gp_i$ belongs to the centre  of the algebra $ \mS_n/\gp_i$.
Note that
$$ \mS_n/\gp_i\simeq K[x_i, x_i^{-1}]\t \mS_{n-1, i}\subset K(x_i)
\t \mS_{n-1, i}\;\; {\rm where}\;\; \mS_{n-1,
i}:=\bigotimes_{j\neq i}^n\mS_1(j).$$ Therefore, $u+\gp_i\in
Z(K(x_i) \t \mS_{n-1, i})\cap K[x_i, x_i^{-1}]\t \mS_{n-1,
i}=K(x_i)^* \cap K[x_i, x_i^{-1}]\t \mS_{n-1, i}=K[x_i,
x_i^{-1}]^*$. Since $u\in (1+\ga_n)^*$, we see that  $u\equiv
1\mod \gp_i$ for all $i$, and so
 $u\in \bigcap_{i=1}^n (1+\gp_i)\cap (1+\ga_n)^*= (1+F_n)\cap
 (1+\ga_n)^* = (1+F_n)^*$. Since $(1+F_n)^*\simeq \GL_\infty (K)$ and
 $Z(\GL_\infty (K))=\{
1\}$, we see that $u=1$ since $u\in Z(\GL_\infty (K))$.

3. Statement 3 follows from statement 2.  $\Box $

$\noindent $

In the algebras $\mS_n$ and $\mA_n$, there are elements that are
invertible linear maps on $P_n$ but not units of the algebras
$\mS_n$ and $\mA_n$ as the following lemma shows.

\begin{lemma}\label{a20Apr9}
$\mS_n^* \subsetneqq \mS_n \cap \Aut_K(P_n)$ and  $\mA_n^*
\subsetneqq \mA_n \cap \Aut_K(P_n)$.
\end{lemma}

{\it Proof}. The element $u:=\prod_{i=1}^n (1-y_i)\in
\mS_n\subseteq \mA_n$ belongs to $\Aut_K(P_n)$, see below,  but
$u\not\in \mA_n^*$ (hence $u\not\in \mS_n^*$) since the element
$u+\ga_n$ is not a unit of the algebra $\CA_n = \mA_n/ \ga_n$. To
prove that the inclusion $u\in \Aut_K(P_n)$ holds  it suffices to
show this for  $n=1$ since $P_n=\bigotimes_{i=1}^n K[x_i]$. The
kernel of the linear map $u$ is equal to zero since $(1-y)*p=0$
for an element $p\in K[x]$ implies that $p=y*p= y^2*p= \cdots =
y^s*p=0$ for all $s\gg 0$ ($y$ is a locally nilpotent map). The
map $u$ is surjective since for each element $q\in K[x]$ there
exists a natural number, say $t$, such that $y^t*q=0$, and so
$q=(1-y^t)*q=u(1+y+\cdots +y^{t-1})*q$. Therefore, $u\in
\Aut_K(P_n)$. $\Box $

$\noindent $

{\bf An inversion formula for  $\s \in G_n$}. Recall that
$G_n\subseteq \mG_n$ (Proposition \ref{f21Mar9}). By Proposition
\ref{f21Mar9}.(1) and Corollary \ref{a19Apr9}, each element $\s
\in G_n$ is a unique product $\s = st_\l \o_\v$ where $s\in S_n$,
$t_\l \in \mT^n$, and $\v \in (1+\ga_n)^*$ ($\ga_n$ is from
Corollary \ref{a19Apr9}).

$\noindent $

{\it Definition}. The unique product $\s = st_\l \o_\v$ is called
the {\em canonical presentation} for $\s $. This is also the
canonical presentation for the automorphism $\s $ treated as an
element of the group $\mG_n$.

$\noindent $

By (\ref{mGinv}), the inversion formula for $\s^{-1}$ is given by
the rule: 
\begin{equation}\label{1mGinv}
\s^{-1} = s^{-1}t_{s(\l^{-1})}\o_{st_\l (\v^{-1})},
\end{equation}
this is the canonical presentation for the automorphism $\s^{-1}$
where $\v^{-1}$ is given by (\ref{1vf18Apr9}).

The next theorem shows that the group $\mG_n$ has trivial centre.

\begin{theorem}\label{16Apr9}
The centre of the group $\mG_n$ is $\{ e\}$.
\end{theorem}

{\it Proof}. The centre of the group $\mA_n^*$ is $K^*$ (Theorem
\ref{24Ma7}). For each element $u\in \mA_n^*$, let $\o_u$ be the
inner automorphism determined by the unit $u$. Let $\s \in
Z:=Z(\mG_n)$. We have to show that $\s =e$. For all elements $u\in
\mA_n^*$, $\o_u= \s \o_u \s^{-1}= \o_{\s (u)}$, and so $\s (u) =
\l u$ for some element $\l = \l_u\in K^*$. In particular, $\s
(H_i) = \l_iH_i$ for all $i=1, \ldots, n$. By Theorem
\ref{6Apr9}.(1), the automorphism $\s$ is a unique product $\s =
s\tau$ where $s\in S_n$ and $\tau \in (\mT^n\times \mU_n) \ltimes
\ker (\xi )$. By Corollary \ref{c18Apr9}, $\s (H_i) \equiv
H_{s(i)}\mod \ga_n$. Hence, $\s (H_i) = H_i$ for all $i$. Then $\s
\in \St_{\mG_n}(H_1, \ldots , H_n)= \mT^n\times \mU_n$ (Corollary
\ref{a11Apr9}.(1)). Note that $1+E_{\alpha\beta} \in \mA_n^*$ for
all $\alpha, \beta \in \N^n$. Hence
$$ 1+\s (E_{\alpha\beta}) = \s (1+E_{\alpha\beta}) =
\l_{\alpha\beta}(1+E_{\alpha\beta})$$ for a nonzero scalar
$\l_{\alpha\beta}$. By Lemma \ref{e21Mar9},  $\l_{\alpha\beta}=1$
and $\s (E_{\alpha\beta})=E_{\alpha\beta}$ for all $\alpha , \beta
\in \N^n$. By  Corollary \ref{yz16Apr9}, $\s =e$, as required.
$\Box $

$\noindent $

Let $H$ be a subgroup of a group $G$. The {\em centralizer}
$\Cen_G(H):=\{ g\in G\, | \, gh=hg\;\; {\rm for \; all}\; h\in
H\}$ of $H$ in $G$ is a subgroup of $G$. In the proof of Theorem
\ref{16Apr9}, we have used only inner derivations of the algebra
$\mA_n$. So, in fact, we have proved there  the next corollary.

\begin{corollary}\label{a16Apr9}
$\Cen_{\mG_n}(\Inn (\mA_n))=\{e\}$.
\end{corollary}

\begin{theorem}\label{A16Apr9}
\begin{enumerate}
\item The centre of the group $\ker (\xi )$ is $\{ e\}$.\item The
centre of the group $(1+\ga_n)^*$ is $\{ 1\}$.
\end{enumerate}
\end{theorem}

{\it Proof}.  1. For each $\alpha \in \N^n$, the element $u(\alpha
) :=1+E_{\alpha 0}$ belongs to the group $(1+\ga_n)^*$. Therefore,
the inner automorphism $\o_{u(\alpha )}$ is an element of $\ker
(\xi )$. Let $\s $ be an element of the centre of the group $\ker
(\xi )$. We have to show that $\s = e$. Note that $\o_{u(\alpha
)}=\s\o_{u(\alpha )}\s^{-1}= \o_{\s (u(\alpha ))}$, and so $\s
(u(\alpha ))= \l_\alpha u(\alpha )$ for a scalar $\l_\alpha \in
K^*$ since $Z(\mA_n^*)=K^*$. It follows from Lemma \ref{e21Mar9}
that $\l_\alpha =1$ and $\s (E_{\alpha 0})=E_{\alpha 0}$ for all
$\alpha \in \N^n$. By Corollary \ref{yz16Apr9}, $\s =e$.

2. By Corollary \ref{b18Apr9}.(2), the groups $\ker (\xi )$ and
$(1+\ga_n)^*$ are isomorphic. Therefore, the centre of the group
$(1+\ga_n)^*$ is $\{ 1\}$, by statement 1. $\Box$

$\noindent $

{\bf Non-embeddability of the proper prime factor algebras of
$\mA_n$ into $\mA_n$}. We are going to prove an `analogue' of
Theorem \ref{24Mar9} for the algebra $\mA_n$ (Theorem
\ref{A10Apr9}).

\begin{theorem}\label{A10Apr9}
No proper prime factor algebra of $\mA_n$ can be embedded into
$\mA_n$ (that is, for each nonzero prime ideal $\gp$ of the
algebra $\mA_n$, there is no algebra monomorphism from $\mA_n/\gp$
into $\mA_n$).
\end{theorem}

{\it Proof}. By  Corollary 3.5, \cite{Bav-Jacalg}, $\gp =
\gp_{i_1}+\cdots +\gp_{i_s}$.   Without loss of generality, we may
assume that $\gp = \gp_1+\cdots +\gp_s$. Suppose that there is a
monomorphism $f: \mA_n/\gp \ra \mA_n$, we seek a contradiction.
For each element $a\in \mA_n$, let $\oa := a+\gp$.  Note that, for
$i=1, \ldots , s$,  $\by_i\bx_i= 1$ and $\bx_i\by_i
=\overline{1-E_{00}(i)}=1$ since $E_{00}(i)\in \gp_i\subseteq
\gp$. The elements $\{ \bx_i, \by_i, \overline{H}_i^{\pm 1}\, | \,
i=1, \ldots , s\}$ are units of the algebra $\mA_n /\gp$, hence
their images under the map $f$ are units of the algebra $\mA_n$.
Let $\pi : \mA_n\ra \mA_n / \ga_n$, $a\mapsto a+\ga_n$. By Theorem
\ref{24Ma7}, the elements $\{ \pi f(\bx_i ), \pi f(\by_i), \pi
f(\overline{H}_i^{\pm 1})\, | \, i=1, \ldots , s\}$ belong to the
commutative algebra $\CL_n$, and so the image of the simple,
noncommutative algebra $\CA_s:= \mA_s/ \ga_s$ under the
compositions of  homomorphisms $\CA_s:= \mA_s/ \ga_s\ra \mA_n /
\gp \stackrel{f}{\ra} \mA_n \stackrel{\pi}{\ra} \mA_n / \ga_n$
belong to the commutative algebra $\CL_n$, a contradiction (the
algebra $\CA_s$ is generated by the elements $x_i+\ga_s$,
$y_i+\ga_s$, $H_i^{\pm 1}+\ga_s$, $i=1, \ldots , s$). $\Box $

$\noindent $

{\it Question. Is an algebra endomorphism of $\mA_n$ necessarily a
monomorphism (automorphism)?}


\section{Stabilizers in $\Aut_{K-{\rm alg}}(\mA_n)$ of the ideals of $\mA_n$}\label{STABAN}

In this section, for each nonzero  ideal $\ga $ of the algebra
$\mA_n$ its stabilizer $\St_{G_n}(\ga ) :=\{ \s \in \mG_n \, | \,
\s (\ga )=\ga \}$ is found (Theorem \ref{15Apr9}) and it is  shown
that the stabilizer $\St_{\mG_n}(\ga )$ has finite index in the
group $\mG_n$ (Corollary \ref{d15Apr9}). When the ideal $\ga$ is
either prime or generic, this result can be refined even further
(Corollary \ref{b15Apr9}, Corollary \ref{c15Apr9}).
 In particular, when $n>1$ the stabilizer
of each height 1 prime  of $\mA_n$ is a maximal subgroup of
$\mG_n$ of index $n$ (Corollary \ref{b15Apr9}.(1)). It is shown
that the ideal $\ga_n$ is the only nonzero, prime,
$\mG_n$-invariant ideal of the algebra $\mA_n$ (Corollary
\ref{b15Apr9}.(3)).

 An ideal $\ga$  of
 $\mA_n$ is called  a {\em proper}
ideal if  $\ga\neq 0, \mA_n$. For an ideal $\ga$ of the algebra
$\mA_n$, $\Min (\ga )$ denotes  the set of all the minimal primes
over $\ga$. Two ideals $\ga$ and $\gb$ are called {\em
incomparable} if neither $\ga \subseteq \gb$ nor $\gb \subseteq
\ga$. The
 ideals of the algebra $\mA_n$ are classified  in
\cite{Bav-Jacalg}. The next theorem shows that each ideal of the
algebra $\mA_n$ is completely determined by its minimal primes. We
use this theorem in the proof of Theorem \ref{15Apr9}.

\begin{theorem}\label{a15Apr9}
{\rm (Theorem 3.8, \cite{Bav-Jacalg})} Let $\ga$ be a proper
 ideal of the algebra $\mA_n$. Then $\Min (\ga )$ is a
finite non-empty set, and the ideal $\ga$ is a unique product and
a unique intersection of incomparable prime  ideals of $\mA_n$
(uniqueness is  up to permutation). Moreover,
$$ a= \prod_{\gp \in \Min (\ga )}\gp =  \bigcap_{\gp \in \Min (\ga )}\gp .$$
\end{theorem}

Let $\Sub_n$ be the set of all the subsets of the set $\{ 1,
\ldots, n \}$. $\Sub_n$ is a partially ordered set with respect to
`$\subseteq $'. Let $\SSub_n$ be the set of all the  subsets of
$\Sub_n$. An element $\{ X_1, \ldots , X_s\}$ of $\SSub_n$ is
called {\em incomparable} if for all $i\ne j$ such that $1\leq
i,j\leq s$ neither $X_i\subseteq X_j$ nor $X_i\supseteq X_j$. An
empty set and one element set are called incomparable by
definition. Let $\Inc_n$ be the subset of $\SSub_n$ of all the
incomparable elements of $\SSub_n$. The symmetric group   $S_n$
acts in the obvious way on the sets $\SSub_n$ and $\Inc_n$ ($\s
\cdot \{ X_1, \ldots ,X_s\}= \{ \s (X_1), \ldots , \s (X_s)\}$).

\begin{theorem}\label{15Apr9}
Let $\ga $ be a proper  ideal of the algebra $\mA_n$. Then
$$ \St_{\mG_n}(\ga ) = \St_{S_n}(\Min (\ga ))\ltimes (\mT^n\times \mU_n)\ltimes
\ker (\xi )= \St_{S_n}(\Min (\ga ))\ltimes (\mT^n\times
\Xi_n)\ltimes \Inn (\mA_n)
$$ where $\St_{S_n}(\Min (\ga )):=\{ \s \in S_n \, | \,
\s (\gq ) \in \Min (\ga )$ for all $\gq \in \Min (\ga )\}$.
Moreover, if $\Min (\ga ) = \{ \gq_1, \ldots , \gq_s\}$ and, for
each number $t=1, \ldots , s$, $\gq_t=\sum_{i\in I_t}\gp_i$ for
some subset $I_t$ of $\{ 1, \ldots , n\}$ then the group
$\St_{S_n}(\Min (\ga ))$ is the stabilizer in the group $S_n$  of
the element $\{ I_1, \ldots , I_s\}$ of $\SSub_n$.
\end{theorem}

{\it Remark}. Note that the group $$\St_{\mG_n}(\Min (\ga ))=
\St_{S_n}(\{ I_1, \ldots , I_s\}):=\{ \s \in S_n \, | \, \{ \s
(I_1), \ldots , \s (I_s)\} = \{ I_1, \ldots , I_s\}\}$$ (and also
the group $\St_{\mG_n}(\ga )$) can be effectively computed in
finitely many steps.

{\it Proof}. Recall that each nonzero prime ideal of the algebra
 $\mA_n$ is a unique sum of height one  prime ideals of the algebra
 $\mA_n$. By Theorem \ref{a15Apr9}, $\St_{\mG_n}(\ga )\supseteq \St_{\mG_n}(\CH_1
 )= (\mT^n\times \mU_n)\ltimes
\ker (\xi )$ (Corollary \ref{a12Apr9}). Since $\mG_n=S_n\ltimes
(\mT^n\times \mU_n)\ltimes \ker (\xi )$ (Theorem \ref{6Apr9}.(1))
and $\Inn (\mA_n)=\mU^0_n\ltimes \ker (\xi )$,
\begin{eqnarray*}
\St_{\mG_n}(\ga )& =&(\St_{\mG_n}(\ga )\cap S_n)\ltimes
(\mT^n\times \mU_n) \ltimes \ker (\xi )=\St_{S_n}(\ga )\ltimes
(\mT^n\times \mU_n)\ltimes\ker (\xi )\\
&=&\St_{S_n}(\ga )\ltimes (\mT^n\times \Xi_n)\ltimes \Inn (\mA_n).
\end{eqnarray*}
By Theorem \ref{a15Apr9}, $\St_{S_n}(\ga )=\St_{S_n}(\Min (\ga ))=
\St_{S_n}(\{ I_1, \ldots , I_s\} )$, and the statement follows.
$\Box$

$\noindent $

The {\em index} of a subgroup $H$ in a group $G$ is denoted by
$[G:H]$.
\begin{corollary}\label{d15Apr9}
Let $\ga $ be a proper  ideal of $\mA_n$. Then $
[\mG_n:\St_{\mG_n}(\ga )] = |S_n:\St_{S_n}(\Min (\ga ))|<\infty$.
\end{corollary}

{\it Proof}. This follows from Theorem \ref{6Apr9}.(1) and Theorem
\ref{15Apr9}.  $\Box $

\begin{corollary}\label{b15Apr9}
\begin{enumerate}
\item $\St_{\mG_n}(\gp_i) \simeq S_{n-1}\ltimes (\mT^n\times
\mU_n)\ltimes \ker (\xi )\simeq S_{n-1}\ltimes (\mT^n\times
\Xi_n)\ltimes \Inn (\mA_n)$, for $i=1, \ldots , n$. Moreover, if
$n>1$ then the groups $\St_{\mG_n}(\gp_i)$ are maximal subgroups
of $\mG_n$ with $[\mG_n :\St_{\mG_n}(\gp_i)]=n$ (if $n=1$ then
$\St_{\mG_1}(\gp_1)=\mG_1$, see statement 3). \item Let $\gp$ be a
nonzero  prime ideal of the algebra $\mA_n$ and $h= \hht (\gp )$
be its height. Then $\St_{\mG_n}(\gp ) \simeq (S_h\times
S_{n-h})\ltimes  (\mT^n\times \mU_n)\ltimes \ker (\xi )\simeq
(S_h\times S_{n-h}) \ltimes (\mT^n\times \Xi_n)\ltimes \Inn
(\mA_n) $.\item The ideal $\ga_n$ is the only nonzero, prime,
$\mG_n$-invariant ideal of the algebra $\mA_n$. \item Suppose that
$n>1$. Let $\gp$ be a nonzero prime ideal of the algebra $\mA_n$.
Then its stabilizer $\St_{\mG_n}(\gp )$ is a maximal subgroup of
$\mG_n$ iff the ideal $\gp $ is of height one.
\end{enumerate}
\end{corollary}

{\it Proof}. 1. Clearly, $\St_{\mG_n}(\gp_i) \cap S_n=\{ \tau \in
S_n\, | \, \tau (\gp_i) = \gp_i\} \simeq S_{n-1}$. By Theorem
\ref{15Apr9}, $\St_{\mG_n}(\gp_i)=S_{n-1}\ltimes (\mT^n\times
\mU_n) \ltimes \ker (\xi )$. When $n>1$, the group
$\St_{\mG_n}(\gp_i)$ is a maximal subgroup of $\mG_n$ since
$$S_{n-1}\simeq \St_{G_n}(\gp_i)/(\mT^n\times \mU_n) \ltimes \ker
(\xi ) \subseteq G_n/ (\mT^n\times \mU_n) \ltimes \ker (\xi )
\simeq S_n$$ and $S_{n-1}=\{ \s \in S_n \, | \, \s (i)=i \}$ is a
maximal subgroup of $S_n$. Clearly, $[\mG_n
:\St_{\mG_n}(\gp_i)]=[S_n:S_{n-1}]=n$.

2. By Corollary 3.5, \cite{Bav-Jacalg},  $\gp=\gp_{i_1}+\cdots +
\gp_{i_h}$ for some distinct indices $i_1, \ldots , i_h\in \{ 1,
\ldots , n\}$. Let $I=\{ i_1, \ldots , i_h\}$ and $CI$ be its
complement. Statement 2 follows from Theorem \ref{15Apr9} and the
fact that
$$\St_{\mG_n} (\gp ) \cap S_n = \{ \s \in S_n \, | \, \s (I) = I, \s
(CI) = CI\}\simeq S_h\times S_{n-h}.$$

3. Since $\ga_n = \gp_1+\cdots +\gp_n$, statement 3 follows from
statement 2.

4. Statement 4 follows from statement 1 and 2. $\Box $

$\noindent $

We are going to apply Theorem \ref{15Apr9} to find the stabilizers
of the generic  ideals (see Corollary \ref{c15Apr9}) but first we
recall the definition of the {\em wreath product} $A\wr B$ of
finite groups $A$ and $B$. The set $\Fun (B,A)$ of all functions
$f: B\ra A$ is a group: $(fg) (b) := f(b) g(b)$ for all $b\in B$
where $g\in \Fun (B,A)$. There is a group homomorphism
$$ B\ra \Aut (\Fun (B,A)), \; b_1\mapsto (f\mapsto b_1(f):b\mapsto
f(b_1^{-1}b)).$$ Then the semidirect product $\Fun (B,A) \rtimes
B$ is called the {\em wreath product} of the groups $A$ and $B$
denoted $A\wr B$, and so the product in $A\wr B$ is given by the
rule:
$$f_1b_1\cdot f_2b_2= f_1b_1(f_2) b_1b_2, \;\; {\rm where}\;\;
f_1, f_2\in \Fun (B,A) , \;\; b_1,b_2\in B.$$ Recall that each
nonzero prime ideal $\gp$ of the algebra $\mA_n$ is a unique sum
$\gp = \sum_{i\in I} \gp_i$ of height one prime ideals. The set
$\Supp (\gp ):= \{ \gp_i\, | \, i\in I\}$ is called the {\em
support} of $\gp$.

$\noindent $

{\it Definition}. We say that a proper  ideal $\ga$ of $\mA_n$ is
{\em generic} if $\Supp (\gp ) \cap \Supp (\gq )=\emptyset$ for
all $\gp , \gq \in \Min (\ga )$ such that $\gp \neq \gq$.

\begin{corollary}\label{c15Apr9}
Let $\ga$ be a generic  ideal of the algebra $\mA_n$. The set
$\Min (\ga )$ of minimal primes over $\ga$ is the disjoint union
of its non-empty subsets, $\Min_{h_1}(\ga ) \bigcup \cdots \bigcup
\Min_{h_t}(\ga )$, where $1\leq h_1<\cdots < h_t\leq n$ and the
set $\Min_{h_i}(\ga )$ contains all the minimal primes over $\ga$
of height $h_i$. Let $n_i:= |\Min_{h_i}(\ga )|$. Then
$$ \St_{G_n}(\ga )= (S_m\times \prod_{i=1}^t(S_{h_i}\wr
S_{n_i}))\ltimes (\mT^n\times \mU_n)\ltimes \ker (\xi )\simeq
(S_m\times \prod_{i=1}^t(S_{h_i}\wr S_{n_i})) \ltimes (\mT^n\times
\Xi_n)\ltimes \Inn (\mA_n) $$ where $m= n-\sum_{i=1}^t n_ih_i$.
\end{corollary}

{\it Proof}. Suppose that $\Min (\ga ) = \{ \gq_1, \ldots ,
\gq_s\}$ and the sets $I_1, \ldots , I_s$ are defined in Theorem
\ref{15Apr9}. Since the ideal $\ga$ is generic, the sets $I_1,
\ldots , I_s$ are disjoint.  By Theorem \ref{15Apr9}, we have to
show that 
\begin{equation}\label{Smm}
\St_{S_m}(\{ I_1, \ldots , I_s\}) \simeq S_m\times
\prod_{i=1}^t(S_{h_i}\wr S_{n_i}).
\end{equation}
The ideal $\ga$ is generic, and so the set $\{ 1, \ldots , n\}$ is
the disjoint union $\bigcup_{i=0}^t M_i$ of its subsets where
$M_i:= \bigcup_{|I_j|=h_i}I_j$, $i=1, \ldots , t$, and $M_0$ is
the complement of the set $\bigcup_{i=1}^tM_i$. Let $S(M_i)$ be
the symmetric group corresponding to the set $M_i$ (i.e. the set
of all bijections $M_i\ra M_i$). Then each element $\s \in
\St_{G_n}(\{ I_1, \ldots , I_s\} )$ is a unique product $\s =
\s_0\s_1\cdots \s_t$ where $\s_i\in S(M_i)$. Moreover, $\s_0$ can
be an arbitrary element of $S(M_0) \simeq S_m$, and, for $i\neq
0$, the element $\s_i$ permutes the sets $\{ I_j\, | \,
|I_j|=h_i\}$ and simultaneously permutes the elements inside each
of the sets $I_j$, i.e. $\s_i\in S_{h_i}\wr S_{n_i}$. Now,
(\ref{Smm}) is obvious. $\Box $

\begin{corollary}\label{e15Apr9}
For each number $s=1, \ldots, n$, let
$\gb_s:=\prod_{|I|=s}(\sum_{i\in I} \gp_i)$ where $I$ runs through
all the subsets of the  set  $\{ 1, \ldots , n\}$ that contain
exactly $s$ elements.  The ideals $\gb_s$ are the only proper,
$\mG_n$-invariant ideals of the algebra $\mA_n$.
\end{corollary}

{\it Proof}. By Theorem \ref{15Apr9}, the ideals  $\gb_s$ are
$\mG_n$-invariant, and they are proper. The converse follows at
once from the classification of  ideals for the algebra $\mA_n$
(Theorem \ref{a15Apr9}) and Theorem \ref{15Apr9}. $\Box $

$\noindent $

$\noindent $

Department of Pure Mathematics

University of Sheffield

Hicks Building

Sheffield S3 7RH

UK

email: v.bavula@sheffield.ac.uk

\end{document}